\documentclass[10pt]{article}

\usepackage{amsmath}
\usepackage{amsmath, amsthm, amsfonts, amssymb, color}
\usepackage{mathrsfs}
\usepackage{latexsym,bm}
\usepackage{graphicx}
\usepackage[toc,page]{appendix}
\usepackage{mathtools}
\mathtoolsset{showonlyrefs}
\usepackage{comment}
\usepackage{xcolor}
\usepackage[numbers, sort&compress]{natbib}

\newtheorem{theorem}{Theorem}[section]

\newtheorem{lemma}[theorem]{Lemma}

\newtheorem{prop}[theorem]{Proposition}
\newtheorem{remark}[theorem]{Remark}
\newtheorem{example}[theorem]{Example}

\definecolor{wco}{rgb}{0.5,0.2,0.3}

\newcommand{\pp}{\mathbb{P}}

\newcommand{\mf}{\mathcal{M}(E)}
\newcommand{\mfo}{\mathcal{M}^{0}(E)}

\newcommand{\ud}{\mathrm{d}}

\newcommand{\R}{\mathbb{R}}
\newcommand{\e}{\mathrm{e}}

\newcommand{\cc}{\mathcal{C}}
\newcommand{\re}{\mathrm{Re}}
\newcommand{\im}{\mathrm{Im}}
\numberwithin{equation}{section} 
\numberwithin{theorem}{section} 

\begin{document}
	
	\allowdisplaybreaks

	\title{\bf	Law of the iterated logarithm for supercritical non-local  spatial branching processes
	}
	\author{  \bf   Haojie Hou\footnote{School of Statistics and Data Science, Nankai University,  Tianjin, 300071, P.R. China.   Email:  houhaojie@nankai.edu.cn.
	 } \hspace{1mm}\hspace{1mm}
		Ting Yang\footnote{Corresponding author. School of Mathematics and Statistics, Beijing Institute of Technology, Beijing, 100081, P.R.China.
			Email: yangt@bit.edu.cn
			The research of this author is supported by NSFC (Grant Nos. 12271374 and 12371143).\hspace{1mm} } \hspace{1mm}\hspace{1mm}
	}
	
	\date{}
	\maketitle

	\begin{abstract}
		Suppose that $X=(X_{t})_{t\ge 0}$ is either a general supercritical non-local branching Markov process, or a general supercritical non-local superprocess, on a Luzin space.
		Here, by ``supercritical" we mean that the mean semigroup of $X$ exhibits a Perron-Frobenius type behaviour with a positive principal eigenvalue.
		In this paper, we study the almost sure behaviour of a family of martingales naturally associated with the real or complex-valued eigenpairs of the mean semigroup.
		Under a fourth-moment condition, we establish
		limit theorems of the iterated logarithm type for these martingales.
		In particular, we discover three regimes, each resulting in different scaling factors and limits.
		Furthermore, we obtain a law of the iterated logarithm for the linear functional $\langle \re(f),X_{t}\rangle$ where $f$ is a sum of finite terms of eigenfunctions and $\re(f)$ denotes its real part.
		In the context of branching Markov processes, our results improve on existing literature by complementing
		the known results
		for multitype branching processes in
		Asmussen [Trans. Amer. Math. Soc. 231 (1) (1977) 233--248]
		and generalizing the recent work
		of
		Hou, Ren and Song [arXiv: 2505.12691v1, arXiv: 2505.12691v2] to allow for non-local branching mechanisms.
		For superprocesses, as far as we know, our results are new.
	\end{abstract}

	\medskip
	

	\medskip


	\section{Introduction}
	This paper investigates the long-time behaviour of a broad class of spatial branching processes that are defined in terms of a spatial motion and a branching mechanism, including both branching Markov processes and superprocesses.
	We mainly focus on the supercritical case where the leading eigenvalue of the mean semigroup is positive.
	
	Recent decades have seen significant developments in the study of the limit theorems for such processes, including (weak and strong) laws of large numbers (LLNs), e.g.,
	\cite{CRW08,CS07, CRY19, EKW,
		E09,EHK10,ET02, EW06,LRS13,W10}
	and the references therein,
	and central limit theorems (CLTs), e.g.,
	\cite{AM2015,Janson04,M18,RSSZ,RSZ,RSZ2015,RSZ2017}
	and the references therein.
	These	results have been extended to the more general framework of non-local branching mechanisms.
	For example, LLNs and CLTs are established in
	\cite{DH2025,CZ2025, HK}
	for general non-local branching Markov processes, and in
	\cite{PY, Y25}
	for general non-local superprocesses.
	
	A natural and important question to ask is, whether there are, in addition, limit theorems of the iterated logarithm type. This question has already been addressed for Galton--Watson processes:
	Let $\big((Z_n)_{ n\geq 0}, \mathbb{P}\big)$ be a supercritical Galton--Watson
	process with $Z_0=1$ and $m:= \mathbb{E}[Z_1]\in (1,\infty)$. It is well-known that $W_n(Z):=  Z_n/m^n$ is a nonnegative martingale, and thus converges almost surely to a finite limit $W_{\infty}(Z)$. Assuming $\mathbb{E}[Z_1^2]<\infty$, Heyde \cite{Heyde1970, Heyde1971} proved that $m^{n/2}(W_{n}(Z)-W_{\infty}(Z))$ converges in distribution to a mixed Gaussian law. It was further proved by Heyde \cite{Heyde1971-2} and Heyde and Leslie \cite{HL1971} that under the same second moment condition, $\mathbb{P}$-a.s. on the non-extinction event,
	\begin{align}\label{Ref1}
		\limsup_{n\to+\infty} \slash \liminf_{n\to+\infty} \frac{m^{n/2}(W_{n}(Z)-W_{\infty}(Z))}{\sqrt{2\log n}} = (+\slash -) \sqrt{\frac{\mathrm{Var}[Z_1]}{m^2-m}W_{\infty}(Z)}.
	\end{align}
	Since $\log n \sim \log \log Z_n$ as $n\to+\infty$, \eqref{Ref1} is often referred to as the \textit{law of the iterated logarithm} (LIL) for the Galton--Watson processes.
	
	While the LIL is well-established for Galton--Watson processes and has been extended to multitype branching processes in \cite{Asmussen1977}, the corresponding result remains largely unexplored in the setting of spatial branching processes. Recently,
	Hou et al. \cite{HRS2025,HRS20252}
	studied the asymptotic behaviour of branching Markov processes where each particle evolves as a
	Hunt process
	and the reproduction occurs locally. In this setting and under restrictive spectral assumptions, LIL results are established for a large class of linear functionals.
	Additionally, we mention two related works \cite{IK} and \cite{IKM}, where the LIL is established, respectively, for the
	Biggins
	martingale in branching random walk and the Nerman's martingale in Crump-Mode-Jagers branching processes.
	However, no analogous results currently exist for superprocesses.

	In this paper, we establish LIL type limit theorems
	for
	both supercritical branching Markov processes and supercritical superprocesses. We work under the assumptions that the mean semigroup displays a natural Perron-Frobenius type behaviour and
	the fourth-moments exist.
	This setting is	notably	general, even allowing for non-local branching mechanisms. Our results thus
	extend the established results
	for multitype branching processes in \cite{Asmussen1977} and branching Markov processes in \cite{HRS2025,HRS20252},
	as well as provide	completely
	new results for superprocesses. Moreover, we are able to develop a unified approach applicable for both branching Markov processes and superprocesses.
	A more detailed discussion on the related literature and the ideas of proof is differed to Sections \ref{sec1.4} and \ref{sec1.5}.
	Before that, we introduce the model in Section \ref{sec1.1}, and state our assumptions and main results in Sections \ref{sec1.2} and \ref{sec1.3}.

	\subsection{Non-local spatial branching processes}
	\label{sec1.1}
	In this paper, we use ``:=" as a way of definition.
	Notations
	$\R$ and $\mathbb{C}$ stand for the sets of real and complex numbers, respectively.
	The letter $c$, with or without subscript, denotes a finite positive constant whose value may vary from place to place.
	Let $E$ be a Luzin topological space
	with Borel $\sigma$-algebra $\mathcal{B}(E)$, and
	$E_{\partial}:=E\cup\{\partial\}$ be the one-point compactification of $E$. Any function on $E$ will be automatically extended to $E_{\partial}$ by setting $f(\partial)=0$.
	We use $\mathcal{B}_{b}(E)$ (resp. $\mathcal{B}^{+}(E)$) to denote the space of bounded (resp. non-negative) measurable functions on $(E,\mathcal{B}(E))$.
	We write $\langle f,\mu\rangle :=\int_{E}f(x)\mu(\mathrm{d}x)$ for a measure $\mu$ on $E$ and a function $f\in\mathcal{B}^{+}(E)$.
	For $k\ge 1$ and a probability space $(\Omega,\mathcal{F},\mathrm{P})$ we use $L^{k}(\mathrm{P})$ to denote the space of $\mathbb{C}$-valued random variables $\eta$ on $\Omega$ satisfying that $\int_{\Omega}|\eta|^{k}\mathrm{d}\mathrm{P}<+\infty$. For a $\mathbb{C}$-valued function $f$, we use $\bar{f}$ for its complex-conjugate, $\re f$ and $\im f$ for its real and imaginary parts, respectively.
	For any two positive functions $f$ and $g$, we use $f\stackrel{c}{\lesssim}g$ (resp., $f\stackrel{c}{\asymp}g$) to denote that, there is a positive constant $c$ such that $f\le c g$ (resp., $c^{-1}g\le f\le c g$) on their common domain of definition. We also write ``$\lesssim$" (resp. ``$\asymp$") if $c$ is unimportant or understood.

	We suppose that
	the spatial motion $\xi=(\Omega,\mathcal{H},\mathcal{H}_{t},\xi_{t},\Pi_{x})$ is a Borel right process in $E$ with transition semigroup $(P_{t})_{t\ge 0}$.
	Here $\{\mathcal{H}_{t}:t\ge 0\}$ is the minimal admissible filtration
	and $\Pi_{x}$ is the law of the process starting from $x$.
	With a slight abuse of notation, we use the same notation $\Pi_x$ to denote the associated expectation operator.
	The Borel right process $\xi$ generates an intrinsic topology, called the fine topology, which is the smallest topology on $E$ that makes all superharmonic functions relative to $\xi$ continuous; see, e.g.\cite[Appendix A.3]{Li}. It is known that a bounded measurable function $f$ is finely continuous relative to $\xi$ if and only if $t\mapsto f(\xi_{t})$ is a.s. right continuous on $[0,+\infty)$.

	\subsubsection{Non-local branching Markov process}
	
	We begin by defining the concept of a non-local branching Markov process (non-local BMP), following the presentation in \cite{HK}.
	A non-local BMP is a collection of particles that evolves according to the following rules: Given their point of creation, particles move independently according to the spatial motion $\xi$. In an event, which we call ``branching", particles positioned at $x\in E$ die at rate $\beta(x)$, where $\beta\in\mathcal{B}^{+}_{b}(E)$, and instantaneously, new particles are created in $E$ according to a point process. The configuration of these offspring are described by the random counting measure
	$$\mathcal{Z}(A)=\sum_{i=1}^{N}\delta_{x_{i}}(A)$$
	for Borel sets $A$ of $E$. The law of the aforementioned point process can depend on $x$, the point of death of the parent, and we denote it by $\mathcal{P}_{x}$, with associated expectation operator given by $\mathcal{E}_{x}$. We assume that $\sup_{x\in E}\mathcal{E}_{x}[N]<+\infty$. The information on the branching event is captured in the so-called branching mechanism
	$$G[f](x):=\beta(x)\mathcal{E}_{x}\left[\prod_{i=1}^{N}f(x_{i})-f(x)\right]\quad\forall x\in E,\ f\in\mathcal{B}^{+}_{1}(E),$$
	where
	$\mathcal{B}^{+}_{1}(E)$ denotes the set of nonnegative measurable functions on $E$ that are bounded by $1$.
	In particular, a BMP is said to be
	``local",
	if all offspring are positioned at their parent's point of death. In that case, the branching mechanism reduces to
	$$G[s](x)=\beta(x)\left(\sum_{k=0}^{+\infty}p_{k}(x)s^{k}-s\right)\quad x\in E,\ s\in [0,1],$$
	where $\{p_{k}(x):k\ge 0\}$ is the offspring distribution when a parent branches at site $x\in E$.
	Without loss of generality, we assume $\mathcal{P}_{x}(N=1)=0$ for all $x\in E$, by viewing a branching event with one offspring as an extra jump in the spatial motion.
	
	Let $\mathcal{N}(E)$ be the space of finite counting measures on $E$. If the configuration of particles at time $t$ is denoted by $\{x_{1}(t),\cdots,x_{N_{t}}(t)\}$, then on the event that the process has not
	become extinct or exploded,
	the BMP can be described as the coordinate process $X:=(X_{t})_{t\ge 0}$ given by
	$$X_{t}(\cdot):=\sum_{i=1}^{N_{t}}\delta_{x_{i}(t)}(\cdot)\quad t\ge 0,$$
	evolving in $\mathcal{N}(E)$.
	In particular, $X$ is a Markov process in $\mathcal{N}(E)$, and its probabilities will be denoted by $\pp:=\{\pp_{\mu}:\ \mu\in\mathcal{N}(E)\}$. Henceforth, we refer to this process
	as a $(P_{t},G)$-BMP.
	With a slightly abuse of notation, we also use the same notation $\mathbb{P}_\mu$ to denote the corresponding expectation operator.

		It is worth noting that the independence in the definition of branching events and movement implies that, if we define
		$$\e^{-U_{t}f(x)}:=\pp_{\delta_{x}}\big[\e^{-\langle f,X_{t}\rangle}\big]\quad\forall t\ge 0,\ x\in E,\ f\in\mathcal{B}^{+}_{b}(E),$$
		then for $\mu\in\mathcal{N}(E)$, we have
		$$\pp_{\mu}\big[\e^{-\langle f,X_{t}\rangle}\big]=\e^{-\langle U_{t}f,\mu\rangle}.$$
		Moreover, for $f\in\mathcal{B}^{+}_{b}(E)$, $t\ge 0$ and $x\in E$,
		\begin{equation}\label{eq:integral equation for U_t}
			\e^{-U_{t}f(x)}=P_{t}\big(\e^{-f}\big)(x)+\int_{0}^{t}P_{s}\big(G\big[\e^{-U_{t-s}f}\big]\big)(x)\mathrm{d}s.
		\end{equation}
		To derive the above equation, it suffices to split the expectation on the first branching event and apply standard reasoning for semigroup integral equations, see e.g. \cite[Chapter 8]{HK}.

	\subsubsection{Non-local superprocess}
		Superprocesses, which are an important class of measure-valued processes, arise as the short lifetime and high density limits of branching Markov processes. There is an extensive body of literature on superprocesses with the so-called local branching mechanisms, such as \cite{Dawson,Dynkin,Eth,Li}. In recent years, research has been extended to the more general framework of non-local branching mechanisms, as seen in \cite{DGL,Li,PY,HKMPR} and the references therein. We introduce the definition of non-local superprocesses in this subsection and refer to \cite{Li} as a standard reference.

	Let $\mf$ denote the space of finite Borel measures on $E$ endowed with
	the weak topology.
	Consider a superprocess $X:=(X_{t})_{t\ge 0}$ with spatial motion $\xi$ and a non-local branching mechanism $\psi$ given by
	\begin{align}
		\psi(x,f)&:=a(x)f(x)+b(x)f(x)^{2}-\eta(x,f)\nonumber\\
		&\quad+\int_{\mfo}\left(\e^{-\nu(f)}-1+\nu(\{x\})f(x)\right)H(x,\mathrm{d}\nu)\quad\forall x\in E,\ f\in\mathcal{B}^{+}(E),\label{def:branching mechanism}
	\end{align}
	where $a\in\mathcal{B}_{b}(E)$, $b\in\mathcal{B}^{+}_{b}(E)$, $\eta(x,\mathrm{d}y)$ is a bounded kernel on $E$,
	and $\left(\nu(1)\wedge \nu(1)^{2}+\nu_{x}(1)\right)H(x,\mathrm{d}\nu)$ is a bounded kernel from $E$ to $\mfo:=\mf\setminus \{0\}$ with $\nu_{x}(\mathrm{d}y)$ defined by $\nu_x(A):=\nu(A\setminus \{x\})$ for all $A\in \mathcal{B}(E)$.
	To be specific, $X$ is an $\mf$-valued Markov process with probabilities $\pp:=\{\pp_{\mu}:\mu\in \mf\}$, satisfying that for every $f\in\mathcal{B}^{+}_{b}(E)$ and $\mu\in\mf$,
	\begin{equation}\label{eq1}
		\pp_{\mu}\left[\e^{-\langle f,X_{t}\rangle}\right]=\e^{-\langle V_{t}f,\mu\rangle},
	\end{equation}
	where $V_{t}f(x):=-\log\pp_{\delta_{x}}\left[\e^{-\langle f,X_{t}\rangle}\right]$ is the unique nonnegative locally bounded solution to the integral equation
	\begin{equation}\label{eq:integral equation for V_{t}}
		V_{t}f(x)=P_{t}f(x)-\Pi_{x}\left[\int_{0}^{t}\psi(\xi_{s},V_{t-s}f)\mathrm{d}s\right].
	\end{equation}
Here the same notation $\mathbb{P}_\mu$ is also used for the corresponding expectation operator.
	By ``locally bounded", we mean that
$$\sup_{x\in E,t\in [0,T]}|V_{t}f(x)|<+\infty \mbox{ for every }T\in (0,+\infty).$$
	Such a process is defined in \cite{Li} via its log-Laplace functional and is referred to as the $(P_{t},\psi)$-superprocess.
	The branching mechanism given in \eqref{def:branching mechanism} is quite general.
	For example, let
	\begin{equation}\nonumber
		\psi^{L}(x, r ):=\hat{a}(x) r +b(x) r ^{2}+\int_{(0,+\infty)}\left(e^{- r  u}-1+ r  u\right)\pi(x,\mathrm{d} u)
	\end{equation}
	for $x\in E$ and $ r \ge 0$,  and
	\begin{equation*}
		\psi^{NL}(x,f):=-\eta(x,f)+\int_{\mfo}\left(\e^{-\nu(f)}-1\right)\Gamma(x ,\mathrm{d} \nu)
	\end{equation*}
	for $x\in E$ and $f\in\mathcal{B}^{+}(E)$, where $b,\eta$ are as in \eqref{def:branching mechanism}, $\hat{a}\in\mathcal{B}_{b}(E)$,  $(u\wedge u^{2})\pi(x,\mathrm{d} u)$ is a bounded kernel from $E$ to $(0,+\infty)$, $\nu(1)\Gamma(x,\mathrm{d} \nu)$ is a bounded kernel from $E$ to $\mfo$. Then \ $(x,f)\mapsto \psi^{L}(x,f(x))+\psi^{NL}(x,f)$ is a branching mechanism that can be represented in the form of \eqref{def:branching mechanism} with $a(x)=\hat{a}(x)-\int_{\mfo}\nu(\{x\})\Gamma(x,\mathrm{d} \nu)$ and $H(x,\mathrm{d}\nu)=\Gamma(x,\mathrm{d}\nu)+\int_{(0,+\infty)}\delta_{r\delta_{x}}(\mathrm{d}\nu)\pi(x,\mathrm{d}r)$.
	A branching mechanism of this type is said to be \textit{decomposable} with local part $\psi^{L}$ and non-local part $\psi^{NL}$.
	In particular, if the non-local part equals $0$, we call such a branching mechanism \textit{purely local}. We emphasize that the branching mechanism considered in this paper is allowed to be non-local and non-decomposable.

	By \cite[Theorem 5.13]{Li}, a $(P_{t},\psi)$-superprocess $X$ has a Borel right realization in $\mf$. We denote by $\mathcal{W}^{+}_{0}$ the space of c\`{a}dl\`{a}g paths from $[0,+\infty)$ to $\mf$ having zero as a trap. Here, we assume that $X$ is the coordinate process in $\mathcal{W}^{+}_{0}$.

	\subsection{Basic assumptions}\label{sec1.2}
	In what follows, $(X,\mathbb{P})$ is taken as either a $(P_{t},G)$-BMP or a $(P_{t},\psi)$-superprocess, and  $(\mathcal{F}_{t})_{t\in[0,\infty]}$ is the filtration generated by this process, which is completed  with the class of $\pp_{\mu}$-negligible measurable sets for
	every initial measure $\mu$.
	
	Define the mean semigroup as
	$$T_{t}f(x):=\pp_{\delta_{x}}\big[\langle f,X_{t}\rangle\big]\quad\forall x\in E,\ f\in\mathcal{B}_{b}(E)\mbox{ and }t\ge 0.$$
	If $(X,\mathbb{P})$ is a $(P_{t},G)$-BMP (resp., a $(P_{t},\psi)$-superprocess), then by \cite[Lemma 8.1]{HK} (resp., \cite[Proposition 2.27]{Li}),
	$$\pp_{\mu}\big[\langle f,X_{t}\rangle\big]=\langle T_{t}f,\mu\rangle\quad\forall \mu\in \mathcal{N}(E)\quad\mbox{(resp., $\forall \mu\in\mf$).}$$
	Moreover,
	in the BMP setting, $T_{t}f(x)$ is the unique locally bounded solution to the integral equation
	\begin{align*}
		T_{t}f(x)=&P_{t}f(x)-\Pi_{x}\left[\int_{0}^{t}\beta(\xi_{s})T_{t-s}f(\xi_{s})\mathrm{d}s\right]+\Pi_{x}\left[\int_{0}^{t}m[T_{t-s}f](\xi_{s})\mathrm{d}s\right],
	\end{align*}
	where $m[f](x):=\beta(x)\mathcal{E}_{x}\big[\sum_{i=1}^{N}f(x_{i})\big]$, while in the superprocess setting,
	\begin{align*}
		T_{t}f(x)=&P_{t}f(x)-\Pi_{x}\left[\int_{0}^{t}a(\xi_{s})T_{t-s}f(\xi_{s})\mathrm{d}s\right]+\Pi_{x}\left[\int_{0}^{t}m[T_{t-s}f](\xi_{s})\mathrm{d}s\right],
	\end{align*}
	where $m[f](x):=\eta(x,f)+\int_{\mfo}\nu_{x}(f)H(x,\mathrm{d}\nu)$.
		These integral equations can be derived by differentiating the non-linear semigroup equations \eqref{eq:integral equation for U_t} for BMPs and \eqref{eq:integral equation for V_{t}} for superprocesses.
		Using the above equations and the Gronwall's inequality, one can easily show that there exists a constant $c_{0}>0$ such that
		\begin{equation}\label{new4}
			\|T_{t}f\|_{\infty}\le \e^{c_{0}t}\|f\|_{\infty}\quad\forall t\ge 0,\ f\in\mathcal{B}_{b}(E).
		\end{equation}

	Throughout this paper we impose
	the following fundamental assumptions.
	\begin{description}
		\item{(H1)} There exist an eigenvalue $\lambda_{1}>0$, a corresponding positive bounded function $\varphi$ and a finite left eigenmeasure $\widetilde{\varphi}$, with normalization $\langle \varphi,\widetilde{\varphi}\rangle=1$, such that for $f\in\mathcal{B}^{+}_{b}(E)$, $\mu\in\mf$ and $t\ge 0$,
		\begin{equation}\label{H1.1}
			\langle T_{t}\varphi,\mu\rangle=\e^{\lambda_{1} t}\langle \varphi,\mu\rangle\mbox{ and }\langle T_{t}f,\widetilde{\varphi}\rangle=\e^{\lambda_{1} t}\langle f,\widetilde{\varphi}\rangle.
		\end{equation}
		Further, we define
		$$
		{  \triangle_{t}:=
		\sup_{x\in E,f\in\mathcal{B}^{+}_{1}(E)}\left| \e^{-\lambda_{1} t}T_{t}f(x)-\langle f,\widetilde{\varphi}\rangle\varphi(x)\right|\quad\forall t\ge 0.}
		$$
		We assume that
		{
			\begin{equation}\label{H1.2}
				\lim_{t\to+\infty}t^{\rho}\triangle_{t}=0\mbox{ for some }\rho>1/2.
			\end{equation}
		}
		
		\item{(H2)} If $(X,\pp)$ is a
		$(P_{t},G)$-BMP, then
		$$\sup_{x\in E}\mathcal{E}_{x}[N^{4}]<+\infty,$$
		and if $(X,\pp)$ is a
		$(P_{t},\psi)$-superprocess, then
		$$\sup_{x\in E}\int_{\{\nu(1)>1\}}
		\nu(1)^{4}H(x,\mathrm{d}\nu)<+\infty.$$
		\item{(H3)} If $(X,\pp)$ is a $(P_{t},\psi)$-superprocess, we additionally assume that, the operators $f\mapsto\psi(\cdot,f)$ and $f\mapsto -a f+m[f]$ preserve $C^{\xi}_{b}(E)$. Here $C^{\xi}_{b}(E)$ denotes the space of bounded measurable functions on $E$ that are finely continuous with respect to $\xi$.
	\end{description}

	Let us give some remarks on these assumptions.
	Assumption (H1) is a Perron-Frobenius type assumption that guarantees the existence of a principal eigenvalue and its corresponding eigenfunction. The spatial branching process is known as {\it supercritical} when $\lambda_{1}>0$.
	Following the terminology in \cite{HK}, processes satisfying (H1), without restriction on the sign of $\lambda_{1}$, are now referred to as the Asmussen--Hering class, acknowledging the fundamental treatment for this class in \cite{AH}.

	Although (H1) may seem like a strong assumption, it can be shown to hold for a large class of spatial branching processes.
	Here we give some concrete examples to illustrate the generality of this setting.

	\begin{example}\rm
		A (continuous-time) Galton--Watson process may be considered as
		the simplest BMP example.
		In this setting, we can take $E=\{0\}$, $\xi$ to be the Markov process which remains at $0$, and a spatially-independent branching mechanism. (H1) trivially holds if we take $\lambda_{1}=\beta(\mathcal{E}[N]-1)$, $\varphi=1$ and $\widetilde{\varphi}=\delta_{\{0\}}$.
		In the context of superprocesses, continuous-state branching process is a natural continuous-mass analogue of Galton--Watson process, and can be shown to satisfy (H1) similarly.
	\end{example}
	
	\begin{example}\rm	
		Multitype branching process is usually considered as a basic example of a non-local BMP. In this setting, $E$ is a finite set. Particles do not move during their lifetime, which is independent and exponentially distributed depending on the particle's state, and reproduce on their death, sending their offspring both to their current site as well as other sites in $E$.
		Then (H1) follows from the Perron-Frobenius theorem applied to the mean semigroup matrix; see e.g. \cite{Janson04}.
	\end{example}

	{ 	\begin{example}\rm
			We observe that $\triangle_{t}$ can be rewritten as
			$$\triangle_{t}=\sup_{x\in E,\,f\in\mathcal{B}^{+}_{1}(E)}\varphi(x)\big|\varphi(x)^{-1}\e^{-\lambda_{1}t}T_{t}f(x)-\langle f,\widetilde{\varphi}\rangle\big|.$$
			Since $\varphi$ is bounded and positive, \eqref{H1.2} is implied by the stronger condition
			\begin{equation}
				\sup_{x\in E,\,f\in\mathcal{B}^{+}_{1}(E)}\big|\varphi(x)^{-1}\e^{-\lambda_{1}t}T_{t}f(x)-\langle f,\widetilde{\varphi}\rangle\big|=O(\e^{-\epsilon t})\mbox{ as }t\to+\infty\mbox{ for some }\epsilon>0.
			\end{equation}
			This latter condition is one of the fundamental assumptions imposed in \cite{HHKW,HKMPR}, and it has been verified for a wide class of non-local spatial branching processes, including the neutron branching process, branching diffusion on a compact domain, multitype continuous-state branching process, and others. For a detailed discussion of these models, we refer the reader to \cite{HHKW,HKMPR} and the references therein.
		\end{example}}

	Assumption (H2) is obviously a fourth-moment condition, which ensures that, if 	$f\in\mathcal{B}_{b}(E)$
	and $t>0$, then $\pp_{\mu}\left[\langle f,X_{t}\rangle^{4}\right]$ is finite for all $\mu\in\mathcal{N}(E)$ in the BMP setting and all $\mu\in\mf$ in the superprocess setting.
	Assumption (H2) allows us to define for any $\mathbb{C}$-valued bounded functions $f$ and $g$,
	$$\vartheta[f,g](x):=\beta(x)\mathcal{E}_{x}\Big[\sum_{i,j=1,i\not=j}^{N}f(x_{i})g(x_{j})\Big]\quad\forall x\in E,$$
	if $(X,\pp)$ is a $(P_{t},G)$-BMP, or,
	$$\vartheta[f,g](x):=2b(x)f(x)g(x)+\int_{\mfo }\nu(f)\nu(g)H(x,{ \ud \nu})\quad\forall x\in E,$$
	if $(X,\pp)$ is a $(P_{t},\psi)$-superprocess. For convenience, we write $\vartheta[f](x)$ for $\vartheta[f,f](x)$.

	In the superprocess setting, it follows from \cite[Proposition 2.38]{Li}) that, for every $t\ge 0$, $f\in\mathcal{B}_{b}(E)$ and $\mu\in\mf$,
	\begin{equation}\label{Second-moment-Sup}
		\pp_{\mu}\left[\langle f,X_{t}\rangle^{2}\right]=\langle T_{t}f,\mu\rangle^{2}+\langle\int_{0}^{t}T_{t-s}(\vartheta[T_{s}f]){ \ud s},\mu\rangle.
	\end{equation}
	This identity implies, in particular, that $\mathrm{Var}_{\mu}[\langle f,X_{t}\rangle]=\langle\mathrm{Var}_{\delta_{\cdot}}\left[\langle f,X_{t}\rangle\right],\mu\rangle$ where
	\begin{equation*}
		\mathrm{Var}_{\delta_{x}}\left[\langle f,X_{t}\rangle\right]=\int_{0}^{t}T_{t-s}(\vartheta[T_{s}f])(x){ \ud s}.
	\end{equation*}
	On the other hand, for BMPs, \cite[Proposition 9.1]{HK} gives
	\begin{equation}\label{Second-moment-BMP}
		\pp_{\delta_{x}}\left[\langle f,X_{t}\rangle^{2}\right]=T_{t}(f^{2})(x)+\int_{0}^{t}T_{t-s}(\vartheta[T_{s}f])(x){ \ud s}
	\end{equation}
	for every $t\ge 0$, $f\in\mathcal{B}_{b}(E)$ and $x\in E$. By the branching property, it follows that for every $\mu\in\mathcal{N}(E)$, $\mathrm{Var}_{\mu}[\langle f,X_{t}\rangle]=\langle\mathrm{Var}_{\delta_{\cdot}}\left[\langle f,X_{t}\rangle\right],\mu\rangle$
	with
	\begin{equation*}
		\mathrm{Var}_{\delta_{x}}\left[\langle f,X_{t}\rangle\right]=T_{t}(f^{2})(x)-(T_{t}f(x))^{2}+\int_{0}^{t}T_{t-s}(\vartheta[T_{s}f])(x){ \ud s}.
	\end{equation*}
	To unify the notation, we define for any $\mathbb{C}$-valued bounded functions $f$ and $g$ and for $x\in E$,
	\begin{equation}\nonumber
		\cc[f,g](x)
		:=\begin{cases}
			f(x)g(x)&\mbox{ if $X$ is a BMP,}\\
			0&\mbox{ if $X$ is a superprocess,}
		\end{cases}
	\end{equation}
	and write $\cc[f]$ for $\cc[f,f]$. Then in both cases, the variance $\mathrm{Var}_{\delta_{x}}[\langle f,X_{t}\rangle]$ can be expressed as
	\begin{equation}\label{eq:variance}
		\mathrm{Var}_{\delta_{x}}[\langle f,X_{t}\rangle]=T_{t}(\cc[f])(x)-\cc[T_{t}f](x)+\int_{0}^{t}T_{t-s}(\vartheta[T_{s}f])(x)\ud{s}.
	\end{equation}

	Assumption (H3), inherited from \cite{Y25}, serves as
	a regularity condition on the branching mechanism. This condition is only required to ensure the existence of a stochastic integral representation for the non-local superprocesses,
which in turn is used to establish an upper bound for the fourth moment of the superprocess (see, e.g.\cite[Lemma A.3]{Y25}). This bound will be employed later in our argument.
So, condition (H3) is not
	strictly necessary
	and can be replaced by some other conditions,
	such as, Conditions 7.1 and 7.2 of \cite[Chapter 7]{Li}.

	\subsection{Main results}\label{sec1.3}
	Before we state the main results of this paper, let us define the martingales with respect to the spatial branching processes.
	We say $(\lambda,g)$ is an eigenpair with respect to the mean semigroup $T_{t}$, if $\lambda\in\mathbb{C}$, and $g$ is a non-trivial $\mathbb{C}$-valued bounded measurable function satisfying that
	$$T_{t}g(x)=\e^{\lambda t}g(x)\quad\forall x\in E,\ t\ge 0.$$
	In that case, define
	\begin{equation}\label{def:martingale}
		W_{t}(\lambda,g):=\e^{-\lambda t}\langle g,X_{t}\rangle \quad\forall t\ge 0.
	\end{equation}
	One can easily show by the Markov property that $(W_{t}(\lambda,g))_{t\ge 0}$
	forms a complex-valued
	$\pp_{\mu}$-martingale with filtration $(\mathcal{F}_{t})_{t\ge 0}$.
	In particular, $W_{t}(\lambda_{1},\varphi)$ is a non-negative martingale, and automatically converges almost surely to a finite limit.
	For notational simplicity, we use $W^{\varphi}_{t}$ to denote
	$W_{t}(\lambda_{1},\varphi)$,
	and $W^{\varphi}_{\infty}$ to denote its limit.

	Martingales of the form \eqref{def:martingale} play an important role in the limit theory of spatial branching processes,
	as evidenced by \cite{KS2, A1969,Asmussen1977,Janson04} for multitype branching processes and by recent works
	\cite{RSZ2017,DH2025,HRS2025,HRS20252}
	for BMPs.
	An analogous object
	in branching random walk
	is the so-called Biggins martingale (with complex parameter),
	whose	fluctuations have been extensively investigated
	in \cite{Biggins92,KM2017,IKM2,ILL2019} and the reference therein.
	There are other related results
	on complex-valued martingales defined in
	various branching models. We refrain from providing references here, and refer to a discussion in \cite{MP} instead.

	\begin{remark}\rm
		(1)
		Suppose $(\lambda,g)$ is an eigenpair with respect to $T_{t}$.
		Since $\e^{\lambda t}\langle g,\widetilde{\varphi}\rangle=\langle T_{t}g,\widetilde{\varphi}\rangle=\e^{\lambda_{1}t}\langle g,\widetilde{\varphi}\rangle$, we have
		$$\langle g,\widetilde{\varphi}\rangle=0\quad \mbox{ if }\quad \lambda\not=\lambda_{1}.$$
		Consequently,
		\begin{equation*}
			\e^{-\lambda_{1}t}T_{t}g(x)-\langle g,\widetilde{\varphi}\rangle \varphi(x)=\begin{cases}
				\e^{(\lambda-\lambda_{1})t}g(x),& \, \lambda\not=\lambda_{1};\\
				g(x)-\langle g,\widetilde{\varphi}\rangle\varphi(x),&\, \lambda=\lambda_{1}.
			\end{cases}
		\end{equation*}
		Assumption (H1) implies that
		$|\e^{-\lambda_{1}t}T_{t}g(x)-\langle g,\widetilde{\varphi}\rangle \varphi(x) |\to 0$
		as $t\to+\infty$. Hence, $\lambda_{1}$ is the largest real eigenvalue such that every other eigenvalue $\lambda$ satisfies that $\re\lambda<\lambda_{1}$. Moreover, $\lambda_{1}$ is a simple eigenvalue, and $\varphi$ is the unique (up to a constant multiple) eigenfunction.
		
		(2) We have for all $t\ge 0$ and $x\in E$,
		\begin{align*}
			T_{t}(\re g)(x)&=\re (T_{t}g(x))=\e^{t\re\lambda}\big(\cos(t\im \lambda)\re g(x)-\sin(t\im \lambda)\im g(x)\big),\\
			T_{t}(\im g)(x)&=\im (T_{t}g(x))=\e^{t\re \lambda}\big(\cos(t\im \lambda)\im g(x)+\sin(t\im \lambda)\re g(x)\big).
		\end{align*}
		If $\im \lambda =0$, then both $\re g$ and $\im g$ are real-valued eigenfunctions corresponding to $\lambda$.
		On the other hand, if $\im\lambda\not=0$, we must have $	\re g	\not\equiv 0$. Otherwise, we get by the above equalities that
		$\sin(t\im \lambda)\im g(x)=0$ for all $t\ge 0$ and $x\in E$, and thus $\im g\equiv 0$,
		which contradicts the fact that $g$ is non-trivial. Similarly, we have $\im g\not\equiv 0$ if $\im \lambda \not=0$.
	\end{remark}
	
	In view of the above remark, we may and do assume that, whenever $\lambda$ is a real eigenvalue, the corresponding eigenfunction $g$ is real-valued.
	
	We are now ready to state our main results, which are derived under the assumptions (H1)-(H3).

	\begin{theorem}\label{them0}
		If $\re\lambda>\lambda_{1}/2$, then for every $x\in E$, the martingale $W_{t}(\lambda,g)$ converges $\pp_{\delta_{x}}$-a.s. and in $L^{2}(\pp_{\delta_{x}})$ to a random variable $W_{\infty}(\lambda,g)$.
	\end{theorem}

	The condition $\re\lambda>\lambda_{1}/2$ is not only sufficient, as shown in Theorem \ref{them0}, but also necessary for the convergence of $W_{t}(\lambda,g)$, as will been seen from Theorem \ref{them1} below. Indeed, Theorem \ref{them1} further describes the decay rate of $W_{t}(\lambda,g)$ to its limit when $\re\lambda>\lambda_{1}/2$, and the growth rate of $W_{t}(\lambda,g)$ when $\re\lambda\le \lambda_{1}/2$.
	The reader will note that our results are stated only for the real parts and the modulus of the martingales.
	The results for the imaginary parts
	can be immediately derived by the identity
	$\im W_{t}(\lambda,g)=\re W_{t}(\lambda,-\mathrm{i}g)$.

	\begin{theorem}\label{them1}
		Suppose $(\lambda,g)$ is an eigenpair.
		 Define
			\begin{align}
				\Sigma(\lambda, g):= 	
				\left\{
				\begin{aligned}
					&\frac{\langle \vartheta[g,\bar{g}],\widetilde{\varphi}\rangle}{\lambda_{1}-2\re\lambda}+\langle\cc[|g|],\widetilde{\varphi}\rangle+\Big| \frac{\langle\vartheta[g],\widetilde{\varphi}\rangle }{\lambda_{1}-2\lambda}+\langle \cc[g],\widetilde{\varphi}\rangle \Big|, & \mbox{if } \re\lambda<\lambda_{1}/2,\\
					&(1+1_{\{\im\lambda=0\}})\langle\vartheta[g,\bar{g}],\widetilde{\varphi} \rangle , & \mbox{if }
					\re\lambda=\lambda_{1}/2,\\
					&\frac{\langle\vartheta[g,\bar{g}],\widetilde{\varphi}\rangle}{2\re\lambda	-\lambda_{1}}-\langle \cc[|g|],\widetilde{\varphi}\rangle+\Big|\frac{\langle\vartheta[g],\widetilde{\varphi}\rangle}{2\lambda-\lambda_{1}}-\langle\cc[g],\widetilde{\varphi}\rangle\Big|, &  \mbox{if } \re\lambda>\lambda_{1}/2.
				\end{aligned}
				\right.
			\end{align}
		The following hold
		$\pp_{\delta_{x}}$-a.s. for every $x\in E$.
		\begin{description}
			\item[(i)] If $\re\lambda<\lambda_{1}/2$, then
			\begin{align*}
				\limsup_{t\to+\infty} \slash \liminf_{t\to+\infty}\frac{\e^{(\re\lambda-\frac{\lambda_{1}}{2})t}\re W_{t}(\lambda,g)}{\sqrt{\log t}}
				=&(+\slash -)\sqrt{
				\Sigma(\lambda, g) W^{\varphi}_\infty
					},
			\end{align*}
			and
			$$\limsup_{t\to+\infty}\frac{\e^{(\re\lambda-\frac{\lambda_{1}}{2})t}\big|  W_{t}(\lambda,g)\big| }{\sqrt{\log t}}=\sqrt{
			\Sigma(\lambda, g) W^{\varphi}_\infty
				}.$$
			\item[(ii)] If $\re\lambda=\lambda_{1}/2$, then
			$$\limsup_{t\to+\infty} \slash \liminf_{t\to+\infty}\frac{\re W_{t}(\lambda,g)}{\sqrt{t\log\log t}}
			=(+\slash -)\sqrt{
			\Sigma(\lambda, g)
		       W^{\varphi}_{\infty}},$$
			and
			$$\limsup_{t\to+\infty} \frac{\big|  W_{t}(\lambda,g)\big| }{\sqrt{t\log\log t}}=
			\sqrt{
			\Sigma(\lambda, g)
				 W^{\varphi}_{\infty}}.$$
			\item[(iii)] If $\re\lambda>\lambda_{1}/2$, then
			\begin{align*}
				&\limsup_{t\to+\infty} \slash \liminf_{t\to+\infty}\frac{\e^{(\re\lambda-\frac{\lambda_{1}}{2})t}\big(\re W_{\infty}(\lambda,g)-\re W_{t}(\lambda,g)\big)}{\sqrt{\log t}}
				 =
				(+\slash -)\sqrt{
			\Sigma(\lambda, g) W^{\varphi}_\infty
					},
			\end{align*}
			and
			$$\limsup_{t\to+\infty} \frac{\e^{(\re\lambda-\frac{\lambda_{1}}{2})t}\big| W_{\infty}(\lambda,g)-W_{t}(\lambda,g)\big|}{\sqrt{\log t}}=\sqrt{
			\Sigma(\lambda, g) W^{\varphi}_\infty
				}.$$
			In particular,
			$$\limsup_{t\to+\infty}\slash \liminf_{t\to+\infty}\frac{\e^{\lambda_{1}t/2}\big(W^{\varphi}_{\infty}-W^{\varphi}_{t}\big)}{\sqrt{\log t}}=(+\slash -)
			\sqrt{\Sigma(\lambda_{1},\varphi)
				W^{\varphi}_{\infty}}.$$
		\end{description}
	\end{theorem}

	\begin{remark}\rm
		\begin{description}
				\item{(i)} The LIL for martingales exhibits three qualitatively distinct regimes, depending on whether $\re\lambda>$, $=$ or $<\lambda_{1}/2$. Equivalently, this trichotomy can be stated in terms of the spectral gap $\lambda_{1}-\re\lambda$ being less than, equal to or greater than $\lambda_{1}/2$. Indeed, the size of the spectral gap determines the growth rate of the martingale variance (see Lemma \ref{lemn2} below), and hence determines how the martingale fluctuates around its mean.

\item{(ii)} We mention that there are two facts standing out in the critical case ($\re\lambda=\lambda_{1}/2$). One is the scaling $\sqrt{t\log\log t}$, which is a delicate conclusion from a martingale analogue of Kolmogorov's law of the iterated logarithm (see, Stout \cite{Stout}).
			A more detailed discussion is provided in Section \ref{sec1.5}.
			The other fact is that, there is an
			additional	factor $2$ in the limit when $\im\lambda=0$.
			This is intuitively reasonable,
			given that	the limit
			of the normalized variance of the martingale depends on
			whether $\im\lambda$ vanishes
			(see, Lemma \ref{lemn2}(ii) below).
			
 \item{(iii)} The proof for the almost sure behaviour of the modulus of the martingales is non-trivial: Based on a key observation that the limit in Theorem \ref{them1}
remains the same when replacing $g$ by $\theta g$ for any $\theta\in\mathbb{C}$ such that $|\theta|=1$, we can prove the LIL for the modulus of the martingale with the help of Lemma \ref{lem:modula}.
		\end{description}	
	\end{remark}

	Let $\mathcal{I}$ be the set of all the eigenpairs $(\lambda, g)$.
	Define
	\begin{align*}
		\mathbb{T}
		:= \left\{ \sum_{k=1}^m g_k:\  m\in \mathbb{N},  (\gamma_k, g_k)\in \mathcal{I},\ \re(\gamma_1)\geq ...\geq \re(\gamma_m)\  \mbox{and}\ \gamma_i\neq \gamma_j,\ \forall  i\neq j  \right\}.
	\end{align*}
	Our next objective is to establish
	an LIL
	for the real part of $\langle h,X_{t}\rangle$ when the test function $h$	belongs to
	$\mathbb{T}$. Noticing that if $(\gamma_{i},g_{i})$ and $(\gamma_{j},g_{j})$ are eigenpairs with $\gamma_{i}=\bar{\gamma_{j}}$, then $\re\big(\langle g_{i}+g_{j},X_{t}\rangle\big)=\re\big(\langle g_{i}+\bar{g}_{j},X_{t}\rangle\big)$, where $g_{i}+\bar{g}_{j}$ is an eigenfunction corresponding to $\gamma_{i}$.
	This implies that, every function in $\mathbb{T}$ can be reduced to a function in the subclass $\widehat{\mathbb{T}}$, where
	\begin{align*}
		\widehat{\mathbb{T}}
		:=  \left\{ \sum_{k=1}^m g_k:\  m\in \mathbb{N},  (\gamma_k, g_k)\in \mathcal{I},\ \re(\gamma_1)\geq ...\geq \re(\gamma_m)\  \mbox{and}\ \gamma_i\neq \gamma_j, \ \gamma_i\neq \bar{\gamma}_j,\ \forall  i\neq j  \right\}.
	\end{align*}
	Therefore, without loss of generality, we only need to consider the subclass $\widehat{\mathbb{T}}$ of test functions.
	Here are our main results for $\langle h, X_t\rangle$.

	\begin{theorem}\label{them4}
		Suppose $h=\sum_{k=1}^m g_k\in
		\widehat{\mathbb{T}}$
		and $x\in E$.
		\begin{description}
			\item[(i)] If $\{\gamma_k: \re(\gamma_{k})=
			\lambda_{1}/2
			\}=\emptyset$, then $\pp_{\delta_{x}}$-a.s.,
			\begin{align*}
				&\limsup_{t\to+\infty}\slash \liminf_{t\to+\infty} \frac{ e^{-\lambda_1 t/2} \re\Big( \langle h, X_{t}\rangle      - \sum_{\re(\gamma_{k})> \frac{\lambda_1}{2}}   e^{\gamma_k t} W_{\infty}(\gamma_k, g_k) \Big)  }{\sqrt{2\log t}}\nonumber\\
				& = (+\slash -) \sqrt{ \big(H^{sm}_{\infty}(h)+H^{la}_{\infty}(h)\big) W_\infty^\varphi},
			\end{align*}
			where
			$$H^{sm}_{\infty}(h):=\int_{0}^{+\infty}\e^{-\lambda_{1}s}\langle\vartheta\big[\sum_{\re(\gamma_{k})<\lambda_{1}/2}\re\big(\e^{\gamma_{k}s}g_{k}\big)\big],\widetilde{\varphi}\rangle \mathrm{d}s
			+\langle\cc\big[\sum_{\re(\gamma_{k})<\lambda_{1}/2}\re(g_{k})\big],\widetilde{\varphi}\rangle
			,$$
			and $$H^{la}_{\infty}(h):=\langle \mathrm{Var}_{\delta_{\cdot}}\big[\sum_{\re(\gamma_{k})>\lambda_{1}/2}\re W_{\infty}(\gamma_{k},g_{k})\big],\widetilde{\varphi}\rangle.$$

			In particular, if $\re(\gamma_{1})<\lambda_{1}/2$, then $\pp_{\delta_{x}}$-a.s.,
			$$\limsup_{t\to+\infty}\slash \liminf_{t\to+\infty} \frac{ e^{-\lambda_1 t/2} \re\big( \langle h, X_{t}\rangle\big)}{\sqrt{2\log t}}= (+\slash -) \sqrt{ H^{sm}_{\infty}(h) W_\infty^\varphi}.$$
			
			\item[(ii)]
			If $\{\gamma_k: \re(\gamma_{k})=
			\lambda_{1}/2
			\}\neq \emptyset$, then $\pp_{\delta_{x}}$-a.s.,
			\begin{equation*}
				\limsup_{t\to+\infty}\slash \liminf_{t\to+\infty} \frac{ e^{-\lambda_1 t/2} \re\Big( \langle h, X_{t}\rangle      - \sum_{\re(\gamma_{k})> \frac{\lambda_1}{2}}   e^{\gamma_k t} W_{\infty}(\gamma_k, g_k) \Big)  }{\sqrt{t\log \log t}}
				= (+\slash -)\sqrt{
					\sum_{\re(\gamma_{k})=\frac{\lambda_1}{2}}
				 \Sigma (\gamma_k, g_k)
				 W^{\varphi}_{\infty} }.
			\end{equation*}
			In particular, if $\re(\gamma_{1})=\lambda_{1}/2$, then $\pp_{\delta_{x}}$-a.s.,
			$$\limsup_{t\to+\infty}\slash \liminf_{t\to+\infty} \frac{ e^{-\lambda_1 t/2} \re\big( \langle h, X_{t}\rangle\big)  }{\sqrt{t\log \log t}}
			= (+\slash -)\sqrt{
				\sum_{\re(\gamma_{k})=\frac{\lambda_1}{2}}
			 \Sigma (\gamma_k, g_k)
			W^{\varphi}_{\infty} }.$$
			Moreover, $\pp_{\delta_{x}}$-a.s.
				\begin{align*}
					\limsup_{t\to+\infty} \frac{ e^{-\lambda_1 t/2} \Big| \langle h, X_{t}\rangle      - \sum_{\re(\gamma_{k})> \frac{\lambda_1}{2}}   e^{\gamma_k t} W_{\infty}(\gamma_k, g_k) \Big|  }{\sqrt{t\log \log t}} = \sqrt{
						\sum_{\re(\gamma_{k})=\frac{\lambda_1}{2}}
					 \Sigma (\gamma_k, g_k)
					W^{\varphi}_{\infty} }.
				\end{align*}
		\end{description}
	\end{theorem}
	
	\begin{remark}\label{Rem}
		\rm	
		\begin{description}		
			\item{(i)} In the above theorem, we can write $e^{-\lambda_1 t/2} \re\Big( \langle h, X_{t}\rangle-\sum_{\re(\gamma_{k})> \frac{\lambda_1}{2}}   e^{\gamma_k t} W_{\infty}(\gamma_k, g_k) \Big)$ as
			\begin{align*}
				&\sum_{\re(\gamma_{k})<\lambda_{1}/2}\e^{(\re(\gamma_{k})-\frac{\lambda_{1}}{2})t}\re\big(\e^{\mathrm{i}t\im(\gamma_{k})}W_{t}(\gamma_{k},g_{k})\big)
				+\sum_{\re(\gamma_{k})=\lambda_{1}/2}\re\big(\e^{\mathrm{i}t\im(\gamma_{k})}W_{t}(\gamma_{k},g_{k})\big)\\
				&+\sum_{\re(\gamma_{k})>\lambda_{1}/2}
				\e^{(\re(\gamma_{k})-\frac{\lambda_{1}}{2})t}\re\Big(\e^{\mathrm{i}t\im(\gamma_{k})}\big(W_{t}(\gamma_{k},g_{k})-W_{\infty}(\gamma_{k},g_{k})\big)\Big).
			\end{align*}
			From this representation, the results of Theorem \ref{them4} can be viewed as an extension of Theorem \ref{them1} to a special class of linear combinations of the martingales.
			Moreover, the above representation together with Theorem \ref{them1} enables us to identify the dominant terms in the asymptotic behaviour.
			Nevertheless, neither the statement nor the proof of Theorem \ref{them4} is a straightforward adaption of Theorem \ref{them1}. In particular, when $h$ is chosen to be a single eigenfunction in Theorem \ref{them4}, the limit in the critical case coincides with that in Theorem \ref{them1}, while in the non-critical case this coincidence occurs if and only if the corresponding eigenvalue is real. This observation reveals, especially in the non-critical case, the oscillatory effects arising from the imaginary parts of the eigenvalues. In addition, the presence of the factor $\e^{\mathrm{i}t\im(\gamma_{k})}$ ($\im(\gamma_{k})$ may not be $0$) brings further technical difficulties in the proof, particularly in the critical case.
			To achieve the main results,
			we 	employ a fundamental theorem
			in the theory of dynamical systems.
			A more detailed discussion is deferred to Section \ref{sec1.5}.
			\item{(ii)} In the non-critical case where $\{\gamma_k: \re(\gamma_{k})=
				\lambda_{1}/2
				\}=\emptyset$, Lemma \ref{lem:modula} is not applicable for obtaining the exact value of the limit superior of the modulus. Nevertheless, we can establish a sharp two-sided estimate as follows.
				For $h=\sum_{k=1}^{m}g_{k}$, define
				$a_{1}(h):=H^{sm}_{\infty}(h)+H^{la}_{\infty}(h)$, $a_{3}(h):=H^{sm}_{\infty}(-\mathrm{i}h)+H^{la}_{\infty}(-\mathrm{i}h)$, and
				\begin{align*}
					a_{2}(h):=&\int_{0}^{+\infty}\e^{-\lambda_{1}s}\langle \vartheta\big[\sum_{\re(\gamma_{k})<\lambda_{1}/2}\re\big(\e^{\gamma_{k}s}g_{k}\big),\sum_{\re(\gamma_{k})<\lambda_{1}/2}\im\big(\e^{\gamma_{k}s}g_{k}\big)\big],\widetilde{\varphi}\rangle\mathrm{d}s\\
					&+\langle\mathcal{C}\big[\sum_{\re(\gamma_{k})<\lambda_{1}/2}\re(g_{k}),\sum_{\re(\gamma_{k})<\lambda_{1}/2}\im(g_{k})\big],\widetilde{\varphi}\rangle\\
					&+\langle\mathrm{Cov}_{\delta_{\cdot}}\big[\sum_{\re(\gamma_{k})>\lambda_{1}/2}\re W_{\infty}(\gamma_{k},g_{k}),\sum_{\re(\gamma_{k})>\lambda_{1}/2}\im W_{\infty}(\gamma_{k},g_{k})\big],\widetilde{\varphi}\rangle.
				\end{align*}
				Set $c_{1}(h):=(a_{1}(h)+a_{3}(h)+\sqrt{(a_{1}(h)-a_{3}(h))^{2}+4a_{2}(h)^{2}})/2$ and $c_{2}(h):=a_{1}(h)+a_{3}(h)$. Then, $\pp_{\delta_{x}}$-a.s.,
				\begin{equation}\label{new5}
					\sqrt{c_{1}(h)W^{\varphi}_{\infty}}\le\limsup_{t\to+\infty}\frac{ e^{-\lambda_1 t/2} \Big|\langle h, X_{t}\rangle      - \sum_{\re(\gamma_{k})> \frac{\lambda_1}{2}}   e^{\gamma_k t} W_{\infty}(\gamma_k, g_k) \Big|}{\sqrt{2\log t}}\le \sqrt{c_{2}(h)W^{\varphi}_{\infty}}.
				\end{equation}
				In particular, if $a_{2}(h)^{2}=a_{1}(h)a_{3}(h)$, then $c_{1}(h)=c_{2}(h)$, and the equality holds in \eqref{new5}. The proof of \eqref{new5} is deferred to the Appendix
				 \ref{Appen3}.
		\end{description}
	\end{remark}

	\subsection{Discussion}\label{sec1.4}	
	As mentioned above,
	there is a variety of branching models satisfying (H1).
	Consequently, under the moment and regularity assumptions on the branching mechanisms, our main theorems apply to these models.

	We remark here that the fourth-moment condition (H2) can be partially relaxed for BMPs. Specifically, in the BMP setting, condition (H2) may be weakened to a second moment condition of the form $\sup_{x\in E} \mathcal{E}_{x}\big[N^{2}\big]<	+\infty$, as demonstrated in \cite{HRS2025} through an auxiliary
	truncation procedure.
	However, this approach fails for superprocesses due to
	the absence of a discrete particle representation.  To unify the treatment of both cases, we retain the stronger assumption (H2) throughout this paper.

	Let us end this section with a short discussion on how our results relate to the existing literature.
	As mentioned earlier, an analogue of Theorem \ref{them1}(iii)
		in the form of $W_{\infty}(Z)-W_{n}(Z)$ is well-known for supercritical Galton--Watson processes.
		If $(\lambda,g)$ is a real-valued eigenpair, then, for multitype branching processes,
		results corresponding to Theorem \ref{them1}(i) and (ii)
		can be derived from \cite[Theorem 2]{Asmussen1977}, while \cite[Theorem 4]{Asmussen1977} provides a discrete-time analogue of Theorem \ref{them1}(iii).
	Beyond these references, however, we have not found any analogous results for
	genuine (non-real) complex-valued eigenfunctions
	in the setting of multitype branching processes.
	Therefore, we believe that Theorem \ref{them1} is a new result in this setting.

	With regard to Theorem \ref{them4}, for multitype branching processes, the limits in the cases $\re(\gamma_{1})\le \lambda_{1}/2$ can be derived from
	\cite[Theorem 2]{Asmussen1977},
	but, the results we provide for $\re(\gamma_{1})>\lambda_{1}/2$ appear to be new.	
	For general BMPs, the most recent and comparable results we are aware of are in \cite{HRS2025}, which inspired this paper,
	but only deals with a special class of BMPs.
		In \cite{HRS2025}, the spatial motion is taken to be a symmetric Hunt process whose transition semigroup admits a density with respect to some positive Radon measure $m$ on $E$, and the branching mechanism is local. In this setting and under additional assumptions on the transition density, it is implied that the infinitesimal generator of the mean semigroup is a compact self-adjoint operator on $L^{2}(E,m)$. Its spectrum consists of countably many real eigenvalues $\lambda_{1}>\lambda_{2}\ge \lambda_{3}\ge\cdots$ with corresponding real-valued bounded eigenfunctions $\{\varphi_{k}:\ k\ge 1\}$ that forms a complete orthonormal basis of $L^{2}(E,m)$. The principal eigenpair $(\lambda_{1},\varphi_{1})$ is chosen to be positive and normalized by $\int_{E}\varphi^{2}_{1}\mathrm{d}m=1$. Furthermore, the mean semigroup $T_{t}$ is \textit{intrinsically ultracontractive} in the following sense: there exists a family of positive symmetric functions $\{p_{t}(x,y):\ t>0\}$ on $E\times E$ such that
		$$T_{t}f(x)=\int_{E}p_{t}(x,y)f(y)m(\mathrm{d}y)\quad\forall t>0,\ x\in E,\ f\in\mathcal{B}^{+}(E),$$
		and
\begin{equation}\label{eq:IU}
		\Big|\frac{\e^{-\lambda_{1}t}p_{t}(x,y)}{\varphi_{1}(x)\varphi_{1}(y)}-1\Big|\le c_{1}\e^{-c_{2}t}\quad\forall (t,x,y)\in (1,+\infty)\times E\times E
\end{equation}
		for some $c_{1},c_{2}>0$.
Clearly, the identities in \eqref{H1.1} are satisfied with with $\varphi=\varphi_{1}$ and $\widetilde{\varphi}=\varphi_{1}m$. However, the assumptions in \cite{HRS2025} do not guarantee that $\widetilde{\varphi}$ is a finite measure. For (H1) to hold, we additionally assume this is the case. This extra condition is satisfied, for instance, if $m$ is a finite measure on $E$. In this setting, by \eqref{eq:IU}, we have for every $x\in E$, $t>1$ and $f\in\mathcal{B}^{+}_{1}(E)$,
		\begin{align*}
			 & \big|\e^{-\lambda_{1}t}T_{t}f(x)-\varphi_{1}(x)\int_{E}f(y)\varphi_{1}(y)m(\mathrm{d}y)\big| \nonumber\\
			 &\le\|\varphi_{1}\|_{\infty}\int_{E}\Big|\frac{\e^{-\lambda_{1}t}p_{t}(x,y)}{\varphi_{1}(x)\varphi_{1}(y)}-1\Big|
			\cdot	f(y)\varphi_{1}(y)
			m(\mathrm{d}y)\le c_{1}\|\varphi_1\|_{\infty}\e^{-c_{2}t} \int_E \varphi_{1}(y)
			m(\mathrm{d}y).
		\end{align*}
		Thus, \eqref{H1.2} of (H1) is also satisfied.
		In this setting, Theorems 2.1-2.6 of \cite{HRS2025} follow as special cases of our Theorem \ref{them4} applied to test functions of the form $f=\sum_{k=1}^{m}a_{k}\varphi_{k}$ with $m\in\mathbb{N}$ and $a_{k}\in\mathbb{R}$.

	For superprocesses, we have been unable to find any comparable results. We believe our results to be the first of their kind
	in the literature.
	
	\subsection{Outline of the paper and strategy for proofs}\label{sec1.5}

	There are three sections in the remainder of this paper.
	
	Section \ref{S2} provides some
	preliminary convergence results.
	In Section \ref{S2.1},
	we first investigate the long time behaviour of the variance of martingales, and then establish a law of large numbers for spatial branching processes.
	In Section \ref{S2.2}, we
	establish an almost sure convergence result for
	the martingale difference $W_{t+s}(\lambda,g)-W_{t}(\lambda,g)$ as $t$ goes to infinity,
	as stated in
	Proposition \ref{lem3.5}.
	The proof contains two main steps: The
	first
	step is to prove the convergence along lattice times (Lemma \ref{lem3.4}). Key to the proof is a refined bound of the maximum difference between the distribution of the normalized martingale difference and a standard normal distribution.
	In the proof of Lemma \ref{lem3.4}, the reader will see
	a fundamental difference between the approaches we take for superprocesses and BMPs.
	This is due to the way in which we obtain the desired bound.
	For a BMP, by the branching property we can write $X_{t+s}=\sum_{i=1}^{N_{t}}X^{(i),t}_{s}$, where $X^{(i),t}_{s}$ denotes the BMP starting from the $i$-th particle alive at time $t$. Thus, conditionally on $X_{t}$, $W_{t+s}(\lambda,g)-W_{t}(\lambda,g)$ can be viewed as a sum of a finite number of independent random variables, which allows us to use
	the Berry--Esseen bound.
	In the case of superprocesses, no such particle representation exists. We turn instead to a L\'{e}vy--Khintchine type representation of the log-Laplace functional, and use Lemma \ref{lemA.3} to measure how close the distribution of the normalized martingale difference is to a standard normal distribution.
	The second step is to pass the limit to continuous time. This is
	achieved via Asmussen and Hering's martingale inequality, see \eqref{them3.0} below.

	Section \ref{S3} contains proofs of Theorems \ref{them0} and \ref{them1}. As will been seen, not only the normalisation and the limits in Theorem \ref{them1} are different when $\re\lambda<\lambda_{1}/2$, $\re\lambda=\lambda_{1}/2$ and $\re\lambda>\lambda_{1}/2$, the proofs are also different.
	The proofs for the cases $\re\lambda<\lambda_{1}/2$ and $\re\lambda>\lambda_{1}/2$, both of which rely
	heavily on Proposition \ref{lem3.5}, are given in Sections \ref{S3.1} and \ref{S3.3}, respectively.
	Essentially, Theorem \ref{them1}(i) and (iii) can be obtained by taking $s\to +\infty$ in Proposition \ref{lem3.5}.
	This approach fails in the critical case ($\re\lambda=\lambda_{1}/2$), since the limit superior of the normalized sequence $W_{t+s}(\lambda,g)-W_{t}(\lambda,g)$ becomes unbounded as $s$ increases. In this case, we use Lemma \ref{lemA.3.2} in the Appendix to replace Proposition \ref{lem3.5}. The former lemma follows directly from Stout's martingale version of Kolmogorov's LIL, which says that, under certain conditions, the limit superior of a martingale $(M_{n},\mathcal{H}_{n})_{n\ge 1}$ behaves like $\sqrt{2s^{2}_{n}\log\log s^{2}_{n}}$ as $n\to+\infty$, where $s^{2}_{n}=\sum_{k=1}^{n}\mathrm{E}\big[(M_{n}-M_{n-1})^{2}|\mathcal{H}_{n-1}\big]$. More specifically, by computing the conditional variance of the critical martingale $\re W_{n\sigma}(\lambda,g)$, we obtain that
	$$\sum_{k=1}^{n}\pp_{\delta_{x}}\left[\Big(\re W_{n\sigma}(\lambda,g)-\re W_{(n-1)\sigma}(\lambda,g)\Big)^{2}|\mathcal{F}_{(n-1)\sigma}\right]\sim c n\sigma W^{\varphi}_{\infty}\mbox{ as }n\to+\infty,$$
	where $c>0$ is written in terms of a second moment functional of the spatial branching process. Combining this with Lemma \ref{lemA.3.2} implies that the limit superior of $\re W_{n\sigma}(\lambda,g)$ behaves like $\sqrt{2c W^{\varphi}_{\infty}n\sigma\log\log n}$, which is exactly the result given in Lemma \ref{lem3.9}.
	With this established, we only need to extend the limit from lattice times to continuous time to complete the proof of Theorem \ref{them1} for the critical case.

	Section \ref{S4} contains the proof of Theorem \ref{them4}. The main idea is as follows. For $h=\sum_{j=1}^{m}g_{j}\in\widehat{\mathbb{T}}$ and $s,t\ge 0$, write
	\begin{align*}
		&\re\big(\langle h,X_{t+s}\rangle-\sum_{\re(\gamma_{k})>\lambda_{1}/2}\e^{\gamma_{k}(t+s)}W_{\infty}(\gamma_{k},g_{k})\big)\\
		=&
		\re\big(\langle
		\sum_{\re(\gamma_{k})\not=\lambda_{1}/2}g_{k}
		,X_{t+s}\rangle-\langle\sum_{\re(\gamma_{k})<\lambda_{1}/2}T_{s}g_{k},X_{t}\rangle-\sum_{\re(\gamma_{k})>\lambda_{1}/2}\e^{\gamma_{k}(t+s)}W_{\infty}(\gamma_{k},g_{k})\big)\\
		&+\re\big(\sum_{\re(\gamma_{k})=\lambda_{1}/2}\langle g,X_{t+s}\rangle\big)+\re\big(\langle \sum_{\re(\gamma_{k})<\lambda_{1}/2}T_{s}g_{k},X_{t}\rangle\big)\\
		=&\mbox{NC}_{h}(s,t)+\mbox{CC}_{h}(s+t)+\mbox{D}_{h}(s,t).
	\end{align*}
	Here NC stands for ``non-critical combination", CC for ``critical combination" and D for ``difference".
	The first step is to establish the LIL for $\mbox{NC}_{h}(s,t)$ as $t$ goes to infinity. This is done via Proposition \ref{lem4.2}, which follows in a similar way to Proposition \ref{lem3.5}. The second step is to show that
	$$\limsup_{t\to+\infty}\e^{-\lambda_{1}t/2}\mbox{CC}_{h}(t)/\sqrt{t\log\log t}$$
	is equal to the desired limit in Theorem \ref{them4}(ii).
	To achieve this, firstly we establish the LIL for
	the critical sum $\sum_{\re(\gamma_{k})=\lambda_{1}/2}W_{t}(\gamma_{k},g_{k})$
	in Proposition \ref{prop4.5}.
	Secondly, based on the identity
	$$\e^{-\lambda_{1}t/2}\mbox{CC}_{h}(t)=\re\big(\sum_{\re(\gamma_{k})=\lambda_{1}/2}\e^{\mathrm{i}t\im(\gamma_{k})}W_{t}(\gamma_{k},g_{k})\big),$$
	we show in Proposition \ref{prop:critical case} that the limit superior of $\e^{-\lambda_{1}t/2}\mbox{CC}_{h}(t)$ behaves like that of $\re\big(\sum_{\re(\gamma_{k})=\lambda_{1}/2}W_{t}(\gamma_{k},g_{k})\big)$.
	The key is to select an increasing sequence of times $\{t_{n}\}$ satisfying that $\e^{\mathrm{i}t_{n}\im(\gamma_{k})}\approx 1$.
	This is where we introduce the definition of \textit{syndetic} sequence, and employ \cite[Theorem 1.2]{Furstenberg} in the theory of dynamical systems.
	The third and last step is to show that $\mbox{D}_{h}(s,t)$ becomes negligible for $t$ and $s$ large.

	Some auxiliary results needed in the proofs are provided in the Appendix.

	\section{Preliminaries}\label{S2}

	\subsection{Some lower-order convergence results}\label{S2.1}
	
	Throughout this paper,
	we assume $(X,\pp)$ is either a $(P_{t},G)$-BMP satisfying (H1) and (H2), or a $(P_{t},\psi)$-superprocess satisfying (H1)-(H3).
We write $\mathbb{P}_{X_s}[F(X_{t})]$ as shorthand for
			\[
			\mathbb{P}_{X_s}[F(\widetilde{X}_{t})]= \mathbb{P}_\mu [F(X_{t+s})| \mathcal{F}_s]
			\]
			where $\widetilde{X}$ is an independent copy of the process starting from the random initial
			state $X_s$. The same convention is adopted for $\mathrm{Var}_{X_s}[F(X_t)]$.
	We further assume that $(\lambda,g)$ is an eigenpair with respect to the mean semigroup $T_{t}$, and set $c:=\re\lambda$.

	A sequence of functions $\{f_{n}:n\ge 1\}$ on $E$ is said to converge boundedly and pointwise to a function $f$,
	if there is a constant $M\ge 0$ such that $\|f_{n}\|_{\infty}\le M$ for $n$ sufficiently large and $f_{n}(x)\to f(x)$ as $n\to+\infty$ for all $x\in E$. We use $\stackrel{b}{\to}$ to denote the bounded pointwise convergence.

	Assumption (H1) implies that for every $f\in\mathcal{B}_{b}(E)$,
		\begin{equation}\label{new6}
			\|\e^{-\lambda_{1}t}T_{t}f-\langle f,\widetilde{\varphi}\rangle\varphi\|_{\infty}\le 2			 \|f\|_{\infty}\triangle_{t}.
		\end{equation}
		Immediately, one has
		$$\e^{-\lambda_{1}t}T_{t}f(x)\le \|f\|_{\infty}\big(2\triangle_{t}+\langle 1,\widetilde{\varphi}\rangle\|\varphi\|_{\infty}\big)\quad\forall x\in E,\ t\ge 0.$$
		Note by (H1) that $\triangle_{t}$ is bounded from above for $t$ sufficiently large. This together with \eqref{new4} and the above inequality yields that, there exists a constant $c_{1}>0$ such that
		\begin{equation}\label{eq1.9}
			\e^{-\lambda_{1}t}T_{t}f(x)\le c_{1}\|f\|_{\infty}\quad\forall x\in E,\ t\ge 0,\ f\in\mathcal{B}_{b}(E).
		\end{equation}

	\begin{lemma}\label{lemn2}
		\begin{description}
			\item[(i)] If $ c<\lambda_{1}/2$, then
			\begin{align*}
				&\e^{-(\lambda_{1}-2 c)t}\mathrm{Var}_{\delta_{x}}[\re W_{t}(\lambda,g)]-\frac{1}{2}\varphi(x)\Big(\frac{\langle\vartheta[g,\bar{g}],\widetilde{\varphi}\rangle}{\lambda_{1}-2 c}
				+\langle\cc[|g|],\widetilde{\varphi}\rangle
				+\re\big(\e^{-2\mathrm{i}t\im\lambda}
				\big(\frac{\langle \vartheta[g],\widetilde{\varphi}\rangle}{\lambda_{1}-2\lambda}
				+\langle \cc[g],\widetilde{\varphi}\rangle
				\big)\big)\Big)\\
				&\stackrel{b}{\to}0\mbox{ as }t\to+\infty.
			\end{align*}
			\item[(ii)] If $ c=\lambda_{1}/2$, then
			$$\frac{1}{t}\mathrm{Var}_{\delta_{x}}[\re W_{t}(\lambda,g)]\stackrel{b}{\to}
			\frac{1+1_{\{\im\lambda=0\}}}{2}\varphi(x)\langle \vartheta[g, \bar{g}],\widetilde{\varphi}\rangle
			\mbox{ as }t\to+\infty.$$
			\item[(iii)] If $ c>\lambda_{1}/2$, then
			$$\mathrm{Var}_{\delta_{x}}[\re W_{t}(\lambda,g)]\stackrel{b}{\to}\frac{1}{2}\int_{0}^{+\infty}\e^{-2c s}T_{s}\big(\vartheta[g,\bar{g}]+\re\big(\e^{-2\mathrm{i} s\im \lambda}\vartheta[g]\big)\big)(x)\mathrm{d}s
			-\cc[\re g](x)
			\mbox{ as }t\to+\infty.$$
			Moreover, it holds that
			\begin{align}\label{Sup-of-W}
				\sup_{t\ge 0,y\in E}\mathrm{Var}_{\delta_{y}}[\re W_{t}(\lambda,g)]<+\infty.
			\end{align}
		\end{description}
	\end{lemma}
	\begin{proof}  Note that $\re W_{t}(\lambda,g)=\e^{-c t}\langle \re\big(\e^{-\mathrm{i}t\im\lambda}g\big),X_{t}\rangle$. By \eqref{eq:variance}, we have
		\begin{align}\label{lemn1.1}
			\mathrm{Var}_{\delta_{x}}[\re W_{t}(\lambda,g)]=&\e^{-2ct}\Big(
			T_{t}\big(\cc\big[\re\big(\e^{-\mathrm{i}t\im\lambda}g\big)\big]\big)(x)-\cc
			\big[T_{t}\big(\re\big(\e^{-\mathrm{i}t\im\lambda}g\big)\big)\big](x)\nonumber\\
			&+\int_{0}^{t}T_{s}\big(\vartheta\big[T_{t-s}\big(\re\big(\e^{-\mathrm{i}t\im \lambda}g\big)\big)\big]\big)(x)\mathrm{d}s\Big).
		\end{align}
		A direct calculation gives that
		$$T_{t-s}\big(\re\big(\e^{-\mathrm{i}t\im \lambda}g\big)\big)=\re\big(\e^{-\mathrm{i}t\im \lambda}T_{t-s}g\big)=\e^{c(t-s)}\re\big(\e^{-\mathrm{i}s\im \lambda}g\big)=\frac{1}{2}\e^{c(t-s)}\big(\e^{-\mathrm{i}s\im\lambda}g+\e^{\mathrm{i}s\im\lambda}\bar{g}\big).$$
		Inserting this to the right hand side of \eqref{lemn1.1} yields that
		\begin{align}\label{step1}
			\mathrm{Var}_{\delta_{x}}[\re W_{t}(\lambda,g)]&=\frac{1}{4}\int_{0}^{t}\e^{-2cs}T_{s}\big(\vartheta\big[\e^{-\mathrm{i}s\im\lambda}g+\e^{\mathrm{i}s\im\lambda}\bar{g}\big]\big)(x)\mathrm{d}s\nonumber\\
			&\quad+\frac{1}{4}\e^{-2 c t}T_{t}\big(\cc\big[\e^{-\mathrm{i}t\im\lambda}g+\e^{\mathrm{i}t\im\lambda}\bar{g}\big]\big)(x)-\cc[\re g](x)\nonumber\\
			&=\frac{1}{2}\int_{0}^{t}\e^{-2c s}T_{s}\big(\vartheta[g,\bar{g}]+\re\big(\e^{-2\mathrm{i} s\im\lambda}\vartheta[g]\big)\big)(x)\mathrm{d}s\nonumber\\
			&\quad+\frac{1}{2}\e^{-2 c t}\Big(T_{t}\big(\cc[|g|]\big)(x)+\re\big(\e^{-2\mathrm{i}t\im\lambda}T_{t}\big(\cc[g]\big)(x)\big)\Big)-\cc[\re g](x)	.
		\end{align}
		We note that by Assumption (H1),
		\begin{align}\label{neweq7}
			\e^{-\lambda_{1}t}\Big(T_{t}\big(\cc[|g|]\big)(x)+\re\big(\e^{-2\mathrm{i}t\im\lambda}T_{t}\big(\cc[g]\big)(x)\big)\Big)-\varphi(x)\Big(\langle\cc[|g|],\widetilde{\varphi}\rangle
			+\re\big(\e^{-2\mathrm{i}t\im\lambda}\langle \cc[g],\widetilde{\varphi}\rangle\big)\Big)\stackrel{b}{\to}0
		\end{align}
		as $t\to+\infty$.
		Hence, the convergence results presented in (i)-(iii) follow immediately from \eqref{neweq7} and Lemma \ref{lemn1} in the Appendix.

		Next we prove \eqref{Sup-of-W}.
		Combining \eqref{step1} and \eqref{eq1.9}, we have
		for $t>0$ and $x\in E$,
		\begin{align*}
			\mathrm{Var}_{\delta_{x}}[\re W_{t}(\lambda,g)]	\le&\e^{-2 c t}T_{t}\big(\cc[|g|]\big)(x)+
			\int_{0}^{t}\e^{-2 c s}T_{s}(\vartheta[g,\bar{g}])(x)\mathrm{d}s\\
			\le&
			c_{1}
			\Big(\e^{(\lambda_{1}-2c)t}\|\cc[|g|]\|_{\infty}+\int_{0}^{t}\e^{(\lambda_{1}-2c)s}\|\vartheta[g,\bar{g}]\|_{\infty}\mathrm{d}s\big)\\
			\le&
			c_{1}
			\Big(\|\cc[|g|]\|_{\infty}+\int_{0}^{+\infty}\e^{(\lambda_{1}-2c)s}
			\|\vartheta[g,\bar{g}]\|_{\infty}\mathrm{d}s\Big),
		\end{align*}
		which implies \eqref{Sup-of-W}.
	\end{proof}

Next we establish a law of large numbers for the spatial branching process.
Analogous
results have been obtained for BMPs in \cite[Corollary 12.1 and Theorem 12.4]{HK}, and for the superprocesses in \cite[Proposition 2.3]{Y25}, albeit under a slightly different convergence assumption on the renormalized mean semigroup.

	\begin{lemma}\label{lem2.2}
		If $(X,\pp)$ is a $(P_{t},G)$-BMP (resp., a $(P_{t},\psi)$-superprocess), then for any $\mathbb{C}$-valued bounded measurable function $f$ and $\mu\in\mathcal{N}(E)$ (resp., $\mu\in\mf$),
		$$\lim_{t\to+\infty}\e^{-\lambda_{1}t}\langle f,X_{t}\rangle=\langle f,\widetilde{\varphi}\rangle W^{\varphi}_{\infty}\mbox{ in }L^{2}(\pp_{\mu}).$$
		Moreover, for any $\sigma>0$,
		$$\lim_{\mathbb{N}\ni n\to+\infty}\e^{-\lambda_{1} n\sigma}
		\langle f,X_{n\sigma}\rangle
		=\langle f,\widetilde{\varphi}\rangle W^{\varphi}_{\infty}\quad\pp_{\mu}\mbox{-a.s.}$$
	\end{lemma}
	
{ \begin{proof}
		Set $\widetilde{f}:= f- \langle f, \widetilde{\varphi}\rangle \varphi$, then $\langle \widetilde{f}, \widetilde{\varphi} \rangle=0$ and $e^{-\lambda_1 t}\langle f, X_t\rangle = e^{-\lambda_1 t} \langle \widetilde{f}, X_t\rangle + \langle f, \widetilde{\varphi}\rangle W_t^\varphi$.
		Note that $\pp_{\mu}[(W^{\varphi}_{t})^{2}]=\langle\varphi,\mu\rangle^{2}+\mathrm{Var}_{\mu}[W^{\varphi}_{t}]$. By Lemma \ref{lemn2}(iii), $\sup_{t\ge 0}\pp_{\mu}[(W^{\varphi}_{t})^{2}]<+\infty$, and hence $W^{\varphi}_{t}$ converges $\pp_{\mu}$-a.s. and in $L^{2}(\pp_{\mu})$ to $W^{\varphi}_{\infty}$. It remains to investigate the asymptotic behaviour of $\e^{-\lambda_{1}t}\langle \widetilde{f},X_{t}\rangle$.

		By \eqref{new6} and (H1), we have $\|t^{\rho}\e^{-\lambda_{1}t}T_{t}\widetilde{f}\|_{\infty}\le 2\|\widetilde{f}\|_{\infty}t^{\rho}\triangle_{t}\to 0$ as $t\to+\infty$. Hence,
		there exists $t_0>0$ such that
		$\|T_t\widetilde{f} \|_\infty  \leq
		{t^{-\rho}}	e^{\lambda_1 t}$ for all $t\ge t_0$.
		Combining this with \eqref{Second-moment-Sup}, \eqref{Second-moment-BMP} and \eqref{eq1.9}, we have for $t>2t_{0}$,
			\begin{align}\label{Moment-super-BMP}
				\mathbb{P}_{\mu} \big[\langle \widetilde{f}, X_{t}\rangle^2   \big]=&\langle T_{t}\widetilde{f},\mu\rangle^{2}+\mathrm{Var}_{\mu}\big[\langle\widetilde{f},X_{t}\rangle\big]\nonumber\\
				\le&\langle T_{t}\widetilde{f},\mu\rangle^{2}+\langle T_{t}(\widetilde{f}^{2}),\mu\rangle+\langle\int_{0}^{t}T_{t-s}\big(\vartheta[T_{s}\widetilde{f}]\big)\mathrm{d}s,\mu\rangle\nonumber\\
				\le&t^{-2\rho}\e^{2\lambda_{1}t}\langle 1,\mu\rangle^{2}+c_{1}\|\widetilde{f}^{2}\|_{\infty}\e^{\lambda_{1}t}\langle 1,\mu\rangle\nonumber\\
				&+\langle c^{2}_{1}\|\widetilde{f}\|^{2}_{\infty}\int_{0}^{t/2}\e^{2\lambda_{1}s}T_{t-s}\big(\vartheta[1]\big)\mathrm{d}s,\mu\rangle+\langle\int_{t/2}^{t}s^{-2\rho}\e^{2\lambda_{1}s}T_{t-s}\big
				(\vartheta[1]\big)\mathrm{d}s,\mu\rangle\nonumber\\
				\le&t^{-2\rho}\e^{2\lambda_{1}t}\langle 1,\mu\rangle^{2}+c_{1}\|\widetilde{f}^{2}\|_{\infty}\e^{\lambda_{1}t}\langle 1,\mu\rangle\nonumber\\
				&+\langle c^{3}_{1}\|\widetilde{f}\|^{2}_{\infty}\|\vartheta[1]\|_{\infty}
				\e^{\lambda_{1}t}\int_{0}^{t/2}\e^{\lambda_{1}s}\mathrm{d}s
				,\mu\rangle+\langle c_{1}\|\vartheta[1]\|_{\infty}\e^{\lambda_{1}t}\int_{t/2}^{t}s^{-2\rho}\e^{\lambda_{1}s}\mathrm{d}s,\mu\rangle\nonumber\\
				\le&c_{2}t^{-2\rho}\e^{2\lambda_{1}t}+c_{3}\e^{\lambda_{1}t}+c_{4}\e^{3\lambda_{1}t/2}.
			\end{align}
			Here $c_{i}$, $i=1,\cdots,4$
are positive constants independent of $t$. Consequently, we have $\pp_{\mu}\big[\big(\e^{-\lambda_{1}t}\langle \widetilde{f}, X_{t}\rangle\big)^2 \big]=O(t^{-2\rho})$ as $t\to+\infty$. Since $\rho>1/2$, $\sum_{n=1}^{+\infty}\pp_{\mu}\big[\big(\e^{-\lambda_{1}n\sigma}\langle \widetilde{f}, X_{n\sigma}\rangle\big)^2 \big]<+\infty$ for any $\sigma>0$. It follows that $\e^{-\lambda_{1}t}\langle \widetilde{f},X_{t}\rangle\to 0$ in $L^{2}(\pp_{\mu})$ as $t\to+\infty$, and by the Borel-Cantelli lemma $\e^{-\lambda_{1}n\sigma}\langle\widetilde{f},X_{n\sigma}\rangle\to 0$ $\pp_{\mu}$-a.s. as $n\to+\infty$. Hence we complete the proof.
		
	\end{proof}}

	\subsection{Almost sure behaviour of the martingale difference in a small time area}\label{S2.2}
	
	Our next result is on the almost sure behaviour of $\re W_{t+s}(\lambda,g)-\re W_{t}(\lambda,g)$ as $t\to+\infty$.
	 To simplify the notation, we define
		\begin{align}
			\Sigma_1(\lambda, g ;s)&:= \langle \vartheta[g,\bar{g}]+(\lambda_{1}-2c)\cc[|g|],\widetilde{\varphi}\rangle F_{s}(\lambda_{1}-2c) \quad \mbox{and} \label{Def-Sigma-1}  \\
			  \Sigma_2(\lambda, g ;s) &:= \langle\vartheta[g]+(\lambda_{1}-2\lambda)\cc[g],\widetilde{\varphi}\rangle F_{s}(\lambda_{1}-2\lambda) \label{Def-Sigma-2},
		\end{align}
		where $F_{s}(z):=z^{-1}(\e^{s z}-1)1_{\{z\not=0\}}+s1_{\{z=0\}}$ for $z\in \mathbb{C}$.
	\begin{prop}\label{lem3.5}
		For $s>0$ and $x\in E$,
		the following holds $\pp_{\delta_{x}}$-a.s.
		\begin{align*}
			& \limsup_{t\to+\infty}\slash \liminf_{t\to+\infty}\frac{\e^{(c-\frac{\lambda_{1}}{2})t}\big(\re W_{t+s}(\lambda,g)-\re W_{t}(\lambda,g)\big)}
			{\sqrt{\log t}}
		=(+\slash -)
			 \sqrt{\big(\Sigma_1(\lambda, g;s)+\big| \Sigma_2(\lambda,g;s)\big| \big)W^{\varphi}_\infty}.
		\end{align*}
	\end{prop}
	
	Our proof of the above result is divided into two steps:
	First to obtain the convergence along lattice times, and then extend it to all times.
	For the first step, we shall
	use the following lemma, whose proof
	is deferred to
	Section \ref{S2.3}.
	\begin{lemma}\label{lem3.4}
		Suppose $\sigma, s>0$, $r\ge 0$ and $x\in E$.
		If $\langle \vartheta[g,\bar{g}]+(\lambda_{1}-2c)\cc[|g|],\widetilde{\varphi}\rangle>0$, then
		\begin{align*}
			&\limsup_{\mathbb{N}\ni n\to+\infty}\frac{\e^{(c-\frac{\lambda_{1}}{2})(n\sigma+r)}\big(\re W_{n\sigma+r+s}(\lambda,g)-\re W_{n\sigma+r}(\lambda,g)\big)}
			{\sqrt{\log n}\cdot\sqrt{
				 \Sigma_1(\lambda,g;s) + \re\big( \e^{-2\mathrm{i}(n\sigma+r)\im\lambda } \Sigma_2(\lambda,g;s)\big)  }
					}
				\\
			&=\sqrt{W^{\varphi}_{\infty}}\quad\pp_{\delta_{x}}\mbox{-a.s.},
		\end{align*}
 where $\Sigma_1$ and $\Sigma_2$ are defined in \eqref{Def-Sigma-1} and \eqref{Def-Sigma-2}, respectively.
		Otherwise, if $\langle \vartheta[g,\bar{g}]+(\lambda_{1}-2c)\cc[|g|],\widetilde{\varphi}\rangle=0$, then
		$$\lim_{\mathbb{N}\ni n\to+\infty}\frac{\e^{(c-\frac{\lambda_{1}}{2})(n\sigma+r)}}{\sqrt{\log n}}\big(\re W_{n\sigma+r+s}(\lambda,g)-\re W_{n\sigma+r}(\lambda,g)\big)=0\quad\pp_{\delta_{x}}\mbox{-a.s.}$$
	\end{lemma}
	Before the proof of Proposition \ref{lem3.5}, we need the following useful inequality:
	\begin{equation}\label{neweq10}
		\big|\langle \vartheta[g]+(\lambda_{1}-2\lambda)\cc[g],\widetilde{\varphi}\rangle\big|\le \langle \vartheta[g,\bar{g}]+(\lambda_{1}-2c)\cc[|g|],\widetilde{\varphi}\rangle.
	\end{equation}
	This inequality holds trivially for superprocesses, since in that case $\cc[g]=\cc[|g|]=0$. Here we provide a proof that applies to both superprocesses and BMPs.
	Similar to \eqref{lemn1.1},
	we have for $t\ge 0$ and $y\in E$,
	\begin{align*}
		&\mathrm{Var}_{\delta_{y}}\big[\langle \re\big(\e^{-\mathrm{i}(t+s)\im\lambda}g\big),X_{s}\rangle\big]\\
		=&\frac{1}{2}\e^{2 c s}\Big[
		\e^{-2 c s}T_{s}(\cc[|g|])(y)-\cc[|g|](y)
		+\re\big(\e^{-2\mathrm{i}t\im\lambda}\big(
		\e^{-2\lambda s}T_{s}(\cc[g])(y)-\cc[g](y)
		\big)\big)\\
		&+\int_{0}^{s}\e^{-2 c r}T_{r}\big(\vartheta[g,\bar{g}]+\re\big(\e^{-2\mathrm{i}(t+r)\im\lambda}\vartheta[g]\big)\big)(y)dr\Big].
	\end{align*}
	A straightforward calculation yields that
	\begin{align}
		&\langle \mathrm{Var}_{\delta_{\cdot}}\big[\langle \re\big(\e^{-\mathrm{i}(t+s)\im\lambda}g\big),X_{s}\rangle\big],\widetilde{\varphi}\rangle
		 =	\frac{1}{2}\e^{2 c s}\big(\Sigma_1(\lambda,g;s)+\re\big(\e^{-2\mathrm{i}t\im\lambda} \Sigma_2(\lambda,g;s)\big)\big)  	
		\nonumber\\	=:&\frac{1}{2}\e^{2 c s}\mathfrak{F}_{s}(t).	\label{neweq11}
	\end{align}
	Since the term on the left hand side is nonnegative,  it follows that $\mathfrak{F}_{s}(t)\ge 0$ for all $t\ge 0$. If $\im\lambda=0$, then $0\le\mathfrak{F}_{s}(t)=2F_{s}(\lambda_{1}-2c)\langle \vartheta[g,\bar{g}]+(\lambda_{1}-2c)
	\cc[|g|]
	,\widetilde{\varphi}\rangle$, which immediately implies \eqref{neweq10}. Otherwise if $\im\lambda\not=0$, by taking $t:=t(s):=\big(\pi+\arg(
 \Sigma_2(\lambda,g;s)
	)\big)/2\im\lambda$, we get
	\begin{align}\label{neweq9}
		0&\le \mathfrak{F}_{s}(t)=
			 \Sigma_1(\lambda,g;s) - |\Sigma_2(\lambda,g;s)|\nonumber\\ & =
		F_{s}(\lambda_{1}-2c)\langle\vartheta[g,\bar{g}]+(\lambda_{1}-2c)	\cc[|g|]	,\widetilde{\varphi}\rangle-\big|F_{s}(\lambda_{1}-2\lambda)	\langle\vartheta[g]+(\lambda_{1}-2\lambda)	\cc[g]	,\widetilde{\varphi}\rangle\big|.
	\end{align}
	Note that $|F_{s}(\lambda_{1}-2\lambda)|/F_{s}(\lambda_{1}-2c)\to 1$ as $s\downarrow 0$. Taking $s\downarrow 0$ in the above inequality yields \eqref{neweq10}.

	\begin{proof}[Proof of Proposition \ref{lem3.5}]
		We only treat the ``$\limsup$" part. The ``$\liminf$" part can be proved similarly by investigating $\re W_{t+s}(\lambda,-g)-\re W_{t}(\lambda,-g)$.
		Fix $s>0$ and $x\in E$.
		It follows from \eqref{neweq10} and Lemma \ref{lem3.4} that
		there is an increasing sequence $\{t_{n}:n\ge 1\}$ (if $\langle \vartheta[g,\bar{g}]+(\lambda_{1}-2c)\cc[|g|],\widetilde{\varphi}\rangle>0$, set $t_{n}=n\sigma_{0}+r_{0}$, where $\sigma_{0}>0$ and $r_{0}\ge  0$
		satisfy
		that
		 $\e^{-2\mathrm{i}(k\sigma_{0}+r_{0})\im\lambda}\Sigma_2(\lambda,g;s)  =|\Sigma_2(\lambda,g;s)|$
		 for all $k\in\mathbb{N}$, otherwise, set
		$t_{n}=n\sigma_0$
		for some $\sigma_{0}>0$), such that, $\pp_{\delta_{x}}$-a.s.
		\begin{align*}
			& \limsup_{\mathbb{N}\ni n\to+\infty}\frac{\e^{(c-\frac{\lambda_{1}}{2})t_{n}}\big(\re W_{t_{n}+s}(\lambda,g)-\re W_{t_{n}}(\lambda,g)\big)}{\sqrt{\log t_{n}}}
			 = \sqrt{\big(\Sigma_1(\lambda,g;s)+|\Sigma_2(\lambda,g;s)|\big)W^{\varphi}_{\infty}}.
		\end{align*}
		We remark here that, by \eqref{neweq10}, the limit on the right hand side vanishes if $\langle \vartheta[g,\bar{g}]+(\lambda_{1}-2c)\cc[|g|],\widetilde{\varphi}\rangle=0$.
		Since $\limsup_{t\to+\infty}\ge\limsup_{t=t_{n}\to+\infty}$,
		it suffices to prove that, $\pp_{\delta_{x}}$-a.s.
		\begin{align}\label{eq3.13}
			&\limsup_{t\to+\infty}\frac{\e^{(c-\frac{\lambda_{1}}{2})t}\big(\re W_{t+s}(\lambda,g)-\re W_{t}(\lambda,g)\big)}{\sqrt{\log t}}
			 \le \sqrt{\big(\Sigma_1(\lambda,g;s)+|\Sigma_2(\lambda,g;s)|\big)W^{\varphi}_{\infty}}.
		\end{align}
		For this purpose, we shall use
		the following inequality for martingales (see, e.g. \cite[Appendix 6.3]{AH} for a discrete-time version, and \cite[proofs of IV.3.4 and VIII.12.3]{AH} for a continuous-time version): Let $(M_{t})_{t\in [0,N]}$ ($N\in (0,+\infty)$) be a martingale with respect to the filtration $\{\mathcal{H}_{t}: t\in [0,N]\}$ and $s^{2}_{t}:=\mathrm{Var}[M_{N}|\mathcal{H}_{t}]$. Then for every $m\in \R$,
		\begin{equation}\label{them3.0}
			\mathrm{P}\left(\sup_{t\in [0,N]}\left(M_{t}-\sqrt{2}s_{t}\right)>m\right)\le 2\mathrm{P}(M_{N}>m).
		\end{equation}
		Recall the definition of $\mathfrak{F}_{s}(t)$ given in \eqref{neweq11}.
		Obviously,
		${\mathfrak{F}}_{\sigma}(t)\le\Sigma_1 (\lambda, g;\sigma) + | \Sigma_2 (\lambda, g;\sigma) |$
		 for all $t\ge 0$  and $\sigma\ge 0$.
	Now fix an arbitrary $\sigma>0$. For $r\ge 0$ and $n\in\mathbb{N}$,
		define
		$$\epsilon_{\sigma,r}(n):=\e^{-(c-\frac{\lambda_{1}}{2})(n\sigma+r)}\sqrt{\log n}.$$
		Then by the
		Markov property, we have for any $\epsilon>0$ and $r\ge 0$,
		\begin{align}\label{eq3.10}
			&\pp_{\delta_{x}}\Big(\sup_{t\in [n\sigma,(n+1)\sigma)}\re W_{t+r}(\lambda,g)-\re W_{n\sigma+r}(\lambda,g)- \sqrt{2\mathrm{Var}_{X_{t+r}}[\re\big(\e^{-\lambda(t+r)}W_{(n+1)\sigma-t}(\lambda,g)\big)]}\nonumber\\
			&\quad\quad>\epsilon_{\sigma,r}(n)\big(\sqrt{(W^{\varphi}_{n\sigma}+\epsilon)\mathfrak{F}_{\sigma}(n\sigma+r)}+\epsilon\big)\,|\,\mathcal{F}_{n\sigma}\Big)\nonumber\\
			=&\pp_{X_{n\sigma}}\Big(\sup_{t\in [0,\sigma)}\re\big(\e^{-\lambda n\sigma}W_{t+r}(\lambda,g)\big)-\re\big(\e^{-\lambda n\sigma}W_{r}(\lambda,g)\big)-\sqrt{2\mathrm{Var}_{X_{t+r}}[\re\big(\e^{-\lambda(n\sigma+t+r)}W_{\sigma-t}(\lambda,g)\big)]}\nonumber\\
			&\quad\quad>\epsilon_{\sigma,r}(n)\big(\sqrt{(\e^{-\lambda_{1}n\sigma}W^{\varphi}_{0}+\epsilon)\mathfrak{F}_{\sigma}(n\sigma+r)}+\epsilon\big)\Big).
		\end{align}
		Noticing that for each $n\in \mathbb{N}$, initial measure $\mu$	and $r\ge 0$,
		$\{\re\big(\e^{-\lambda n\sigma}W_{t+r}(\lambda,g)\big)-\re\big(\e^{-\lambda n\sigma}W_{r}(\lambda,g)\big):\ t\ge 0\}$ is a $\pp_{\mu}$-martingale with respect to the filtration $\{\mathcal{F}_{t+r}:\ t\ge 0\}$. Moreover, by the Markov property, we have for $t\in [0,\sigma)$,
		\begin{align*}
			\mathrm{Var}[\re\big(\e^{-\lambda n\sigma}W_{\sigma+r}(\lambda,g)\big)-\re\big(\e^{-\lambda n\sigma}W_{r}(\lambda,g)\big)\,|\,\mathcal{F}_{t+r}]
			&=\mathrm{Var}[\re\big(\e^{-\lambda n\sigma}W_{\sigma+r}(\lambda,g)\big)\,|\,\mathcal{F}_{t+r}]\\
			&=\mathrm{Var}_{X_{t+r}}[\re\big(\e^{-\lambda(n\sigma+t+r)}W_{\sigma-t}(\lambda,g)\big)].
		\end{align*}
		Applying \eqref{them3.0} to the right hand side of \eqref{eq3.10}, we get
		\begin{align}\label{eq3.14}
			&\mbox{LHS of }\eqref{eq3.10}\nonumber\\
			&\le 2\pp_{X_{n\sigma}}\Big(\re\big(\e^{-\lambda n\sigma}W_{\sigma+r}(\lambda,g)\big)-\re\big(\e^{-\lambda n\sigma}W_{r}(\lambda,g)\big)>\epsilon_{\sigma,r}(n)\big(\sqrt{(\e^{-\lambda_{1}n\sigma}W^{\varphi}_{0}+\epsilon)\mathfrak{F}_{\sigma}(n\sigma+r)}+\epsilon\big)\Big)\nonumber\\
			&=2\pp_{\delta_{x}}\big(\re\big(W_{(n+1)\sigma+r}(\lambda,g)\big)-\re\big(W_{n\sigma+r}(\lambda,g)\big)>\epsilon_{\sigma,r}(n)
			\sqrt{\big(W^{\varphi}_{n\sigma}+\epsilon)\mathfrak{F}_{\sigma}(n\sigma+r)}+\epsilon\big)\,|\,\mathcal{F}_{n\sigma}\big),
		\end{align}
		where in last equality we used the Markov property.
		Since $\lim_{\mathbb{N}\ni n\to+\infty}W^{\varphi}_{n\sigma}=W^{\varphi}_{\infty}$ $\pp_{\delta_{x}}$-a.s., it follows from Lemma \ref{lem3.4} that, for any $\epsilon>0$,
		$$\re W_{(n+1)\sigma+r}(\lambda,g)-\re W_{n\sigma+r}(\lambda,g)\le \epsilon_{\sigma,r}(n)\big(\sqrt{(W^{\varphi}_{n\sigma}+\epsilon)\mathfrak{F}_{\sigma}(n\sigma+r)}+\epsilon\big)\mbox{ eventually }\pp_{\delta_{x}}\mbox{-a.s.}$$
		Thus by the conditional Borel--Cantelli lemma, we have
		$$\sum_{n=1}^{+\infty}\pp_{\delta_{x}}\left(\re W_{(n+1)\sigma+r}(\lambda,g)-\re W_{n\sigma+r}(\lambda,g)>\epsilon_{\sigma,r}(n)\big(\sqrt{(W^{\varphi}_{n\sigma}+\epsilon)\mathfrak{F}_{\sigma}(n\sigma+r)}+\epsilon\big)\,|\,\mathcal{F}_{n\sigma}\right)<+\infty\quad \pp_{\delta_{x}}\mbox{-a.s.}$$
		This together with \eqref{eq3.14} implies that, the sum from $n=1$ to infinity
		of the left hand side of \eqref{eq3.10} is
		a.s. finite.
		Again, by the conditional Borel--Cantelli lemma, we have
		\begin{align}\label{eq3.11}
			&\sup_{t\in [n\sigma,(n+1)\sigma)}\re W_{t+r}(\lambda,g)-\re W_{n\sigma+r}(\lambda,g)-\sqrt{2\mathrm{Var}_{X_{t+r}}[\re\big(\e^{-\lambda(t+r)}W_{(n+1)\sigma-t}(\lambda,g)\big)]}\nonumber\\
			\le&\epsilon_{\sigma,r}(n)\big(\sqrt{(W^{\varphi}_{n\sigma}+\epsilon)\mathfrak{F}_{\sigma}(n\sigma+r)}+\epsilon\big)\mbox{ eventually }\pp_{\delta_{x}}\mbox{-a.s.}
		\end{align}
		We write $n(t)=\lfloor t/\sigma\rfloor$.
		Note that
		\begin{align}
			&\frac{\e^{(c-\frac{\lambda_{1}}{2})t}\big(\re W_{t+s}(\lambda,g)-\re W_{t}(\lambda,g)\big)}{\sqrt{\log t} }\nonumber\\
			=&\frac{\e^{(c-\frac{\lambda_{1}}{2})t}\big(\re W_{t+s}(\lambda,g)-\re W_{n(t)\sigma+s}(\lambda,g)\big)}{\sqrt{\log t} }-\frac{\e^{(c-\frac{\lambda_{1}}{2})t}\big(\re W_{t}(\lambda,g)-\re W_{n(t)\sigma}(\lambda,g)\big)}{\sqrt{\log t} }\nonumber\\
			&\quad+\e^{(c-\frac{\lambda_{1}}{2})(t-n(t)\sigma)}\sqrt{\frac{\log n(t)}{\log t}}\cdot\frac{\e^{(c-\frac{\lambda_{1}}{2})n(t)\sigma}\big(\re W_{n(t)\sigma+s}(\lambda,g)-\re W_{n(t)\sigma}(\lambda,g)\big)}{\sqrt{\log n(t)}}.\label{eq3.9}
		\end{align}
		For the third term in the right hand side of \eqref{eq3.9}, Lemma \ref{lem3.4} yields that, $\pp_{\delta_{x}}$-a.s.,
		\begin{align}\label{neweq3}
			&\limsup_{t\to+\infty}\frac{\e^{(c-\frac{\lambda_{1}}{2})n(t)\sigma}\big(\re W_{n(t)\sigma+s}(\lambda,g)-\re W_{n(t)\sigma}(\lambda,g)\big)}{\sqrt{\log n(t)}}
			 \le \sqrt{\big(\Sigma_1(\lambda,g;s)+|\Sigma_2(\lambda,g;s)|\big)W^{\varphi}_{\infty}}.
		\end{align}
		The fist term in the right hand side of \eqref{eq3.9} can be written as
		$I_{\sigma,s}(t)+II_{\sigma,s}(t)$, where
		$$I_{\sigma,s}(t):=\frac{\e^{(c-\frac{\lambda_{1}}{2})t}}{\sqrt{\log t} }\big(\re W_{t+s}(\lambda,g)-\re W_{n(t)\sigma+s}(\lambda,g)-\sqrt{2\mathrm{Var}_{X_{t+s}}[\re\big(\e^{-\lambda(t+s)}W_{(n(t)+1)\sigma-t}(\lambda,g)\big)]}\big),$$
		and
		$$II_{\sigma,s}(t):=\frac{\e^{(c-\frac{\lambda_{1}}{2})t}}{\sqrt{\log t} }\sqrt{2\mathrm{Var}_{X_{t+s}}[\re\big(\e^{-\lambda(t+s)}W_{(n(t)+1)\sigma-t}(\lambda,g)\big)]}.$$
		Moreover,  $I_{\sigma,s}(t)$ can be written as
		\begin{align*}
			&\frac{\re W_{t+s}(\lambda,g)-\re W_{n(t)\sigma+s}(\lambda,g)-\sqrt{2\mathrm{Var}_{X_{t+s}}[\re\big(\e^{-\lambda(t+s)}W_{(n(t)+1)\sigma-t}(\lambda,g)\big)]}}{\epsilon_{\sigma,s}(n(t))}\\
			&\times\e^{(c-\frac{\lambda_{1}}{2})(t-n(t)\sigma-s)}\sqrt{\frac{\log n(t)}{\log t}}.
		\end{align*}
		Combining \eqref{eq3.11}, the inequality
		$$
		 \big|\re\big(\e^{-2\mathrm{i}t\im\lambda}\Sigma_2(\lambda,g;s)\big)\big|\le \Sigma_1(\lambda, g;s)
		$$
		(which is a consequence of \eqref{neweq9})
		and the fact that $W^{\varphi}_{n(t)\sigma}\to W^{\varphi}_{\infty}$ $\pp_{\delta_{x}}$-a.s., we get that
		\begin{equation}\label{lem3.5.1}
			\limsup_{t\to+\infty}I_{\sigma,s}(t)\le \e^{|c-\frac{\lambda_{1}}{2}|(\sigma+s)}
			 \sqrt{2 \Sigma_1(\lambda, g;\sigma)  W^{\varphi}_{\infty}}
			\quad\pp_{\delta_{x}}\mbox{-a.s.}
		\end{equation}
		On the other hand,
			 we have
			\begin{align}\label{lem3.5.2}
				& \e^{(2c-\lambda_{1})(t+s)}\mathrm{Var}_{X_{t+s}}[\re\big(\e^{-\lambda(t+s)}W_{(n(t)+1)\sigma-t}(\lambda,g)\big)]\nonumber\\
				& =\e^{-\lambda_{1}(t+s)}\langle\mathrm{Var}_{\delta_{\cdot}}
				[\re\big(\e^{-\mathrm{i}(t+s)\im\lambda }W_{(n(t)+1)\sigma-t}(\lambda,g)\big)],X_{t+s}\rangle.
			\end{align}
		{ 	We note that for $y\in E$ and $t\in [n \sigma,(n +1)\sigma)$,
			\begin{align*}
				&\mathrm{Var}_{\delta_{y}}[\re\big(\e^{-\mathrm{i}(t+s)\im\lambda }W_{(
					n
					+1)\sigma-t}(\lambda,g)\big)]\\
				=&
				\int_{0}^{(n +1)\sigma-t} \e^{-2c q}T_{q}\big(\vartheta[\re\big(\e^{-\mathrm{i}(t+s+q)\im\lambda}g\big)]\big)(y)\mathrm{d}q\\
				&+\e^{-2 c((n +1)\sigma-t)} \Big(
				T_{(n+1)
					\sigma-t}\big(\cc\big[\re\big(\e^{-\mathrm{i}((n
					+
					1)\sigma+s)\im\lambda}g\big)\big]\big)(y)-\cc\big[\re\big(\e^{-\mathrm{i}(t+s)\im\lambda}g\big)\big](y)\Big) \\
				\le&
				\int_{0}^{(n+1)\sigma-t}
				\e^{-2 c q}T_{q}(\vartheta[|g|])(y)\mathrm{d}q
				+\e^{-2 c(
					(n+1)\sigma-t)}T_{(n+1)\sigma-t}
				\big(\cc[|g|]\big)(y).
			\end{align*}
			Noticing that there is some constant	$c_{1}$
			such that $  \vartheta[|g|] \vee  \cc[|g|] \le c_{1}$,
			we have
			\begin{align}
				& \e^{-\lambda_{1}(t+s)}\langle\mathrm{Var}_{\delta_{\cdot}}
				[\re\big(\e^{-\mathrm{i}(t+s)\im\lambda }W_{(n+1)\sigma-t}(\lambda,g)\big)],X_{t+s}\rangle \nonumber\\
				\le &   \e^{2|c|\sigma} c_{1} \int_{0}^{(n+1)\sigma-t} e^{-\lambda_1(t+s)} \langle T_{q}1, X_{t+s}\rangle \mathrm{d}q
				+\e^{2|c|\sigma} c_{1} e^{-\lambda_1(t+s)} \langle T_{(n+1)\sigma-t}1,X_{t+s}\rangle.
			\end{align}
			Set $Q(x):= 1- \langle 1, \widetilde{\varphi} \rangle \varphi$.
			Then $T_{q}1=T_{q}Q+\langle 1,\widetilde{\varphi}\rangle \e^{\lambda_{1}q}\varphi$ for $q\ge 0$. Inserting to the right hand side yields that, for $t\in [n\sigma,(n+1)\sigma)$,
			\begin{align}\label{e2}
				& \e^{-\lambda_{1}(t+s)}\langle\mathrm{Var}_{\delta_{\cdot}}
				[\re\big(\e^{-\mathrm{i}(t+s)\im\lambda }W_{(n+1)\sigma-t}(\lambda,g)\big)],X_{t+s}\rangle \nonumber\\
				\le &\e^{2|c|\sigma}c_{1}\Big(\int_{t}^{(n+1)\sigma}\e^{-\lambda_{1}(t+s)}\langle T_{q-t}Q,X_{t+s}\rangle\mathrm{d}q+\langle 1,\widetilde{\varphi}\rangle W^{\varphi}_{t+s}\int_{0}^{(n+1)\sigma-t}\e^{\lambda_{1}q}\mathrm{d}q\Big)\nonumber\\
				&+\e^{2|c|\sigma}c_{1}\Big(\e^{-\lambda_{1}(t+s)}\langle T_{(n+1)\sigma-t}Q,X_{t+s}\rangle+\langle 1,\widetilde{\varphi}\rangle \e^{\lambda_{1}((n+1)\sigma-t)}W^{\varphi}_{t+s}\Big)\nonumber\\
				\le&\e^{(2|c|+\lambda_{1})\sigma}c_{1}\langle 1,\widetilde{\varphi}\rangle(\sigma+1)
W^{\varphi}_{t+s}
+\e^{2|c|\sigma}c_{1}\int_{t}^{(n+1)\sigma}\e^{-\lambda_{1}(t+s)}\big|\pp_{\delta_{x}}\big[\langle Q,X_{q+s}\rangle|\mathcal{F}_{t+s}\big]\big|\mathrm{d}q\nonumber\\
				&+\e^{2|c|\sigma}c_{1}\e^{-\lambda_{1}(t+s)}\big|\pp_{\delta_{x}}\big[\langle Q,X_{(n+1)\sigma+s}\rangle|\mathcal{F}_{t+s}\big]\big|\nonumber\\
				\leq &   \e^{(2|c|+\lambda_1)\sigma} c_{1}
				\langle 1,\widetilde{\varphi}\rangle(\sigma+1) W_{t+s}^\varphi
				+ \e^{2|c|\sigma} c_{1} e^{-\lambda_1(n\sigma+s)} \int_{n\sigma}^{(n+1)\sigma} \sup_{n \sigma \leq t\leq q} \big| \mathbb{P}_{\delta_x} \big[\langle Q, X_{q+s}\rangle \big| \mathcal{F}_{t+s} \big]\big| \mathrm{d}q
				\nonumber\\
				&+\e^{2|c| \sigma} c_{1}e^{-\lambda_1(n\sigma+s)} \sup_{n \sigma \leq t\leq (n+1)\sigma} \big| \mathbb{P}_{\delta_x} \big[\langle Q, X_{s+(n+1)\sigma}\rangle \big| \mathcal{F}_{t+s} \big]\big| \nonumber\\
				=:&  \e^{(2|c|+\lambda_1)\sigma} c_{1}
				\langle 1,\widetilde{\varphi}\rangle(\sigma+1) W_{t+s}^\varphi
				+ \e^{2|c|\sigma} c_{1} e^{-\lambda_1(n\sigma+s)} (J_{\sigma, s}(n) + K_{\sigma, s}(n)).
			\end{align}
			Applying Jensen's inequality and
			Doob's maximal inequality
			yields that
			\begin{align}\label{e1}
				& \mathbb{P}_{\delta_x}[J_{\sigma, s}(n)^2] +  \mathbb{P}_{\delta_x}[K_{\sigma, s}(n)^2] \nonumber\\
				\leq & \sigma \int_{n\sigma}^{(n+1)\sigma} \mathbb{P}_{\delta_x} \Big[\sup_{n \sigma \leq t\leq q} \big| \mathbb{P}_{\delta_x} \big[\langle Q, X_{s+q}\rangle \big| \mathcal{F}_{t+s} \big]\big|^2  \Big]\mathrm{d}q + \mathbb{P}_{\delta_x} \Big[
				\sup_{n \sigma \leq t\leq (n+1)\sigma}
				\big| \mathbb{P}_{\delta_x} \big[\langle Q, X_{s+(n+1)\sigma}\rangle \big| \mathcal{F}_{t+s} \big]\big|^2  \Big] \nonumber\\
				\leq &4\sigma \int_{n\sigma}^{(n+1)\sigma} \mathbb{P}_{\delta_x}\big[\langle Q, X_{s+q}\rangle^2   \big] \mathrm{d}q + 4\mathbb{P}_{\delta_x} \big[\langle Q, X_{s+(n+1)\sigma}\rangle^2   \big].
			\end{align}
		Applying a computation similar to that in \eqref{Moment-super-BMP} (with $\widetilde{f}$ and $\mu$ replaced by $Q$ and $\delta_{x}$, respectively), we can show that
			$$e^{-2\lambda_1(n\sigma+s)} \Big(\mathbb{P}_{\delta_x} [J_{\sigma, s}(n)^2] +  \mathbb{P}_{\delta_x} [K_{\sigma, s}(n)^2]\Big)=O(n^{-2\rho})\mbox{ as }n\to+\infty.$$
			Thus, by the Borel-Cantelli lemma, $e^{-\lambda_1(n\sigma+s)} (J_{\sigma, s}(n) + K_{\sigma, s}(n)) \to 0$ $\pp_{\delta_{x}}$-a.s. as $n\to\infty$. Combing this with \eqref{e2}, we get that
		\begin{align}\label{eq3.26}
			\limsup_{t\to\infty}  \e^{-\lambda_{1}(t+s)}\langle\mathrm{Var}_{\delta_{\cdot}}
			[\re\big(\e^{-\mathrm{i}(t+s)\im\lambda }
			W_{(n(t)+1)\sigma-t}(\lambda,g)
			\big)],X_{t+s}\rangle \leq  \e^{(2|c|+\lambda_1)\sigma} c_{1}
			\langle 1,\widetilde{\varphi}\rangle (\sigma+1)W^{\varphi}_{\infty},
		\end{align}
	and consequently,
	$$\limsup_{t\to+\infty}II_{\sigma,s}(t)=0
	\quad\pp_{\delta_{x}}\mbox{-a.s.}$$}
	This together with \eqref{lem3.5.1} yields that, 	$\pp_{\delta_{x}}$-a.s.,
	\begin{align*}
		&\limsup_{t\to+\infty}\frac{\e^{(c-\frac{\lambda_{1}}{2})t}\big(\re W_{t+s}(\lambda,g)-\re W_{n(t)\sigma+s}(\lambda,g)\big)}{\sqrt{\log t} }
		 \le \e^{|c-\frac{\lambda_{1}}{2}|(\sigma+s)}\sqrt{2\Sigma_1(\lambda,g;\sigma) W^{\varphi}_\infty}.
	\end{align*}
	Similarly by investigating $\re W_{t+s}(\lambda,-g)-\re W_{n(t)\sigma+s}(\lambda,-g)$, we can prove that, $\pp_{\delta_{x}}$-a.s.,
	\begin{align*}
		&\liminf_{t\to+\infty}\frac{\e^{(c-\frac{\lambda_{1}}{2})t}\big(\re W_{t+s}(\lambda,g)-\re W_{n(t)\sigma+s}(\lambda,g)\big)}{\sqrt{\log t} }
		 \ge - \e^{|c-\frac{\lambda_{1}}{2}|(\sigma+s)}\sqrt{2\Sigma_1(\lambda,g;\sigma) W^{\varphi}_\infty}.
	\end{align*}
	Note that $\Sigma_{1}(\lambda,g;\sigma)
	\to 0$ as $\sigma\to 0$. We have
	\begin{align}\label{lem3.5.3}
		&\lim_{\sigma\to 0}\liminf_{t\to+\infty}\frac{\e^{(c-\frac{\lambda_{1}}{2})t}\big(\re W_{t+s}(\lambda,g)-\re W_{n(t)\sigma+s}(\lambda,g)\big)}{\sqrt{\log t} }\nonumber\\
		=&\lim_{\sigma\to 0}\limsup_{t\to+\infty}\frac{\e^{(c-\frac{\lambda_{1}}{2})t}\big(\re W_{t+s}(\lambda,g)-\re W_{n(t)\sigma+s}(\lambda,g)\big)}{\sqrt{\log t} }=0\quad\pp_{\delta_{x}}\mbox{-a.s.}
	\end{align}
	Applying similar argument with $s=0$,
	we can prove that
	\begin{equation}\label{lem3.5.4}
		\lim_{\sigma\to 0}\liminf_{t\to+\infty}\frac{\e^{(c-\frac{\lambda_{1}}{2})t}\big(\re W_{t}(\lambda,g)-\re W_{n(t)\sigma}(\lambda,g)\big)}{\sqrt{\log t} }=0\quad\pp_{\delta_{x}}\mbox{-a.s.}
	\end{equation}
	In view of \eqref{neweq3}, \eqref{lem3.5.3} and \eqref{lem3.5.4}, we get \eqref{eq3.13} by letting first $t\to +\infty$ and then $\sigma\to 0$ in the right hand side of \eqref{eq3.9}.
\end{proof}

\subsection{Proof of Lemma \ref{lem3.4}}\label{S2.3}
The main tool in the proof of Lemma \ref{lem3.4} is the following lemma, see \cite[Appendix 7.2]{AH}.
\begin{lemma}\label{lem6.1}
	Let $\{\mathcal{G}_{n}:n\ge 1\}$ be an increasing sequence of $\sigma$-fields and $\{Z_{n}:n\ge 1\}$ be a sequence of random variables.
	Let $\Phi$ be the standard normal distribution.
	If
	$$\sum_{n=1}^{+\infty}\sup_{x\in\R}\left|\mathrm{P}\left(Z_{n}\le x|\mathcal{G}_{n}\right)-\Phi(x)\right|<+\infty
	\quad\mathrm{P}\mbox{-a.s.},
	$$
	then
	$$\limsup_{\mathbb{N}\ni n\to+\infty}\frac{Z_{n}}{\sqrt{2\log n}}\le 1
	\quad\mathrm{P}\mbox{-a.s.}
	$$
	If, furthermore, there exists
	a $k\in \mathbb{N}$ such that $Z_{n}$ is $\mathcal{G}_{n+k}$-measurable for each $n\ge 1$, then
	$$\limsup_{\mathbb{N}\ni n\to+\infty}\frac{Z_{n}}{\sqrt{2\log n}}= 1
	\quad\mathrm{P}\mbox{-a.s.}
	$$
\end{lemma}

We will prove Lemma \ref{lem3.4} separately for superprocesses and BMPs.

\subsubsection{Proof for superprocesses}
In the superprocess setting, our approach to Lemma \ref{lem3.4} relies on
a L\'{e}vy--Khintchine type representation of the log-Laplace functional.
Let us first introduce this type of representation.

Let $Q_{t}(\mu,\mathrm{d}\nu)$ denote the transition semigroup of the $(P_{t},\psi)$-superprocess.
The semigroup $Q_{t}$ is said to have the regular branching property in the sense that,
$$\int_{\mf}\e^{-\nu(f)}Q_{t}(\mu,\mathrm{d}\nu)=\e^{-\langle V_{t}f,\mu\rangle}\quad\forall \mu\in\mf,\ f\in\mathcal{B}^{+}_{b}(E),$$
where $V_{t}f(x)=
-\log\int_{\mf}\e^{-\nu(f)}Q_{t}(\delta_{x},\mathrm{d}\nu)$.
Consequently, by
\cite[Theorem 2.4 and Theorem 1.36]{Li}, $(V_{t})_{t\ge 0}$ is a cumulant semigroup on $\mathcal{B}^{+}_{b}(E)$ with the following
canonical representation
\begin{equation}\label{eq:canonical representation for V_{t}}
	V_{t}f(x)=\Lambda_{t}(x,f)+\int_{\mfo}\left(1-\e^{-\nu(f)}\right)L_{t}(x,\mathrm{d}\nu)\quad\forall f\in\mathcal{B}^{+}_{b}(E),\ t\ge 0,\ x\in E.
\end{equation}
Here $\Lambda_{t}(x,\mathrm{d}y)$ is a bounded kernel on $E$ and $(1\wedge \nu(1))L_{t}(x,\mathrm{d}\nu)$ is a bounded kernel from $E$ to $\mfo$.
We recall that, under (H2) the second moment of $\langle f,X_{t}\rangle$ is finite for every $f\in\mathcal{B}_{b}(E)$. Using the fact that $\pp_{\delta_{x}}\left[\langle f,X_{t}\rangle\right]=\left.-\frac{\partial}{\partial \theta}\e^{-V_{t}(\theta f)(x)}\right|_{\theta=0}$ and $\pp_{\delta_{x}}\left[\langle f,X_{t}\rangle^{2}\right]=\left.\frac{\partial^{2}}{\partial \theta^{2}}\e^{-V_{t}(\theta f)(x)}\right|_{\theta=0}$, we get from \eqref{eq:canonical representation for V_{t}} and \cite[Proposition 1.39]{Li}
that for every $x\in E$ and $t\ge 0$,
$\int_{\mfo}\nu(1)^{2}L_{t}(x,\mathrm{d}\nu)<+\infty$,
\begin{equation}\label{eq:second moment2}
	T_{t}f(x)=\pp_{\delta_{x}}\left[\langle f,X_{t}\rangle\right]=\Lambda_{t}(x,f)+\int_{\mfo}\nu(f)L_{t}(x,\mathrm{d}\nu)\mbox{ and }
	\mathrm{Var}_{\delta_{x}}[\langle f,X_{t}\rangle]=
	\int_{\mfo}|\nu(f)|^{2}L_{t}(x,\mathrm{d}\nu)
\end{equation}
for every $f\in\mathcal{B}^{+}_{b}(E)$.

It follows from \eqref{eq1} and \eqref{eq:canonical representation for V_{t}} that for every $\mu\in\mf$, $t\ge 0$ and $f\in\mathcal{B}^{+}_{b}(E)$,
\begin{equation}\label{eq2}
	\pp_{\mu}\left[\e^{-\langle f,X_{t}\rangle}\right]=\exp\Big\{-\Big(\int_{E}\Lambda_{t}(x,f)\mu(\mathrm{d}x)+\int_{E}\mu(\mathrm{d}x)\int_{\mfo}\left(1-\e^{-\nu(f)}\right)L_{t}(x,\mathrm{d}\nu)\Big)\Big\}.
\end{equation}
This implies that, for every $t\ge 0$ and $\mu\in\mf$, the distribution of the random measure $X_{t}(\cdot)$ under $\pp_{\mu}$ is equal to
$$\int_{E}\Lambda_{t}(x,\cdot)\mu(\mathrm{d}x)+\int_{\mf}\nu(\cdot)N^{\mu}_{t}(\mathrm{d}\nu),$$
where $N^{\mu}_{t}(\mathrm{d}\nu)$ denotes a Poisson random measure on $\mf$ with intensity $\int_{E}L_{t}(x,\mathrm{d}\nu)\mu(\mathrm{d}x)$.

Let $\mathbb{C}^{+}:=\{z\in\mathbb{C}:\ \re z\ge 0\}$.
The result of
\cite[Lemma A.1]{Y25} implies that, for every $\mathbb{C}^{+}$-valued bounded measurable function $f$ on $E$, the expectation of $\exp\{-\int_{\mf}\nu(f)N^{\mu}_{t}(\mathrm{d}\nu)\}$ exists and is equal to $$\exp\{\int_{E}\mu(\mathrm{d}x)\int_{\mfo}\left(\e^{-\nu(f)}-1\right)L_{t}(x,\mathrm{d}\nu)\}.$$ From this, one can see that the restriction that $f$ is real-valued in \eqref{eq2} is unnecessary, and \eqref{eq2} is indeed valid for all $\mathbb{C}^{+}$-valued bounded measurable functions $f$ on $E$.

For complex-valued random variables $X$ and $Y$, we use the standard definitions $\mathrm{Cov}(X,Y):=\mathrm{E}[(X-\mathrm{E}X)\overline{(Y-\mathrm{E}Y)}]$ and $\mathrm{Var}[X]:=\mathrm{Cov}(X,X)=\mathrm{E}[|X-EX|^{2}]$.
	Note that the identities in \eqref{eq:second moment2} remain valid
	for any $\mathbb{C}^{+}$-valued bounded measurable function $f$ on $E$. Moreover, we have
	$$\mathrm{Cov}_{\delta_{x}}\big(\langle f,X_{t}\rangle,\langle \bar{f},X_{t}\rangle\big)=\pp_{\delta_{x}}\big[\big(\langle f,X_{t}\rangle-T_{t}f(x)\big)^{2}\big]=\int_{\mfo}\nu(f)^{2}L_{t}(x,\mathrm{d}\nu).$$
	Combining this with \eqref{eq2} yields, for every $\mathbb{C}^{+}$-valued bounded measurable function $f$,
\begin{align}\label{eq3}
	\pp_{\mu}\left[\e^{-\langle f,X_{t}\rangle}\right]
	=\exp\Big\{-\langle T_{t}f-\frac{1}{2}
	\mathrm{Cov}_{\delta_{\cdot}}\big(\langle f,X_{t}\rangle,\langle\bar{f},X_{t}\rangle\big)
	+{\mathfrak{L}_{t}}[f],\mu\rangle\Big\},
\end{align}
where
$$\mathfrak{L}_{t}[f](x):=\int_{\mfo}\left(1-\e^{-\nu(f)}-\nu(f)+\frac{1}{2}\nu(f)^{2}\right)L_{t}(x,\mathrm{d}\nu).$$

The following two lemmas are needed in the proof of Lemma \ref{lem3.4}.
\begin{lemma}\label{lemA.3}
	Let $\Phi$ be the standard normal distribution and $F$ be an arbitrary distribution on $\R$. Let $f(\theta):=\int_{-\infty}^{+\infty}\e^{\mathrm{i}\theta y}\mathrm{d} F(y)$ be the characteristic function of $F$. Then
	$$\sup_{x\in \R}\left|F(x)-\Phi(x)\right|\le \frac{1}{\pi}\int_{-\infty}^{+\infty}\left|\frac{f(\theta)-\e^{-\frac{1}{2}\theta^{2}}}{\theta}\right|\mathrm{d}\theta.$$
\end{lemma}
\begin{proof}  The Berry--Esseen basic inequality (cf. \cite[4.1.a]{LB}) implies that for every $T>0$ and $b>\frac{1}{2\pi}$,
	$$\sup_{x\in\R}\left|F(x)-\Phi(x)\right|\le b\int_{-T}^{T}\left|\frac{f(\theta)-\e^{-\frac{1}{2}\theta^{2}}}{\theta}\right|\mathrm{d}\theta+2b T\sup_{x\in\R}\int_{-\kappa/T}^{\kappa/T}|\Phi(x+y)-\Phi(x)|\mathrm{d}y.$$
	Here $\kappa=\kappa(b)$ is a positive constant depending only on $b$.
	
	Putting $b=\frac{1}{\pi}$ and using the fact that $\sup_{x\in\R}|\Phi'(x)|\le \frac{1}{\sqrt{2\pi}}$, we get
	\begin{eqnarray*}\sup_{x\in\R}\left|F(x)-\Phi(x)\right|&\le& \frac{1}{\pi}\int_{-T}^{T}\left|\frac{f(\theta)-\e^{-\frac{1}{2}\theta^{2}}}{\theta}\right|\mathrm{d}\theta+\frac{2T}{\pi}\sup_{x\in\R}\int_{-\kappa/T }^{\kappa/T}\frac{1}{\sqrt{2\pi}}|y|\mathrm{d}y\\
		&\le&\frac{1}{\pi}\int_{-T}^{T}\left|\frac{f(\theta)-\e^{-\frac{1}{2}\theta^{2}}}{\theta}\right|\mathrm{d}\theta+\frac{\sqrt{2}\kappa^{2}}{\pi^{3/2}T}.
	\end{eqnarray*}
	Hence we get the desired result by letting $T\to+\infty$.
\end{proof}

\begin{lemma}\label{lem3.2}
	For every $t>0$ and $f\in\mathcal{B}^{+}_{b}(E)$, $x\mapsto \int_{\mfo}\nu(f)^{4}L_{t}(x,\mathrm{d}\nu)$ is a bounded function on $E$.
\end{lemma}
\begin{proof}  Fix arbitrary $t>0$ and $f\in\mathcal{B}^{+}_{b}(E)$. The characteristic function of $\langle f,X_{t}\rangle$ can be represented by
	\begin{equation}\label{lem3.2.1}
		\pp_{\delta_{x}}\left[\e^{\mathrm{i}\theta\langle f,X_{t}\rangle}\right]=\e^{\Phi^{f}_{t,x}(\theta)},
	\end{equation}
	where
	$$\Phi^{f}_{t,x}(\theta)=\mathrm{i}\theta\Lambda_{t}(x,f)-\int_{\mfo}\left(1-\e^{\mathrm{i}\theta\nu(f)}\right)L_{t}(x,\mathrm{d}\nu).$$
	Taking derivatives with respect to $\theta$ in both sides of \eqref{lem3.2.1}, we get that
	\begin{align*}
		\int_{\mfo}\nu(f)^{4}L_{t}(x,\mathrm{d}\nu)=&\left.\frac{\partial^{4}}{\partial \theta^{4}}\Phi^{f}_{t,x}(\theta)\right|_{\theta=0}\\
		=&\pp_{\delta_{x}}\left[\langle f,X_{t}\rangle^{4}\right]
		+3\pp_{\delta_{x}}\left[\langle f,X_{t}\rangle\right]^{4}+6\pp_{\delta_{x}}\left[\langle f,X_{t}\rangle\right]^{2}\cdot\mathrm{Var}_{\delta_{x}}\left[\langle f,X_{t}\rangle\right]\\
		&\quad-3 \mathrm{Var}_{\delta_{x}}\left[\langle f,X_{t}\rangle\right]^{2}
	+4\pp_{\delta_{x}}\left[\langle f,X_{t}\rangle\right]\cdot\pp_{\delta_{x}}\left[\langle f,X_{t}\rangle^{3}\right].
	\end{align*}
	By H\"{o}lder's inequality we have
	$$\int_{\mfo}\nu(f)^{4}L_{t}(x,\mathrm{d}\nu)\lesssim \pp_{\delta_{x}}\left[\langle f,X_{t}\rangle^{4}\right]\quad\forall x\in E.$$
	So it suffices to prove that
	$\sup_{x\in E}\pp_{\delta_{x}}\left[\langle f,X_{t}\rangle^{4}\right]<+\infty.$
	\cite[Lemma A.3]{Y25} implies that there is some
	$c_{0}>0$
	such that
	\begin{equation}\label{eq3.6}
		\pp_{\delta_{x}}\left[\langle f,X_{t}\rangle^{4}\right]\le
		c_{0}
		\left(T_{t}f(x)^{4}+
		\mathrm{Var}_{\delta_{x}}[\langle f,X_{t}\rangle]^{2}
		+\int_{0}^{t}T_{s}\left(\vartheta[\mathrm{Var}_{\delta_{\cdot}}[\langle f,X_{t-s}\rangle]]\right)(x)\mathrm{d}s\right)\quad\forall x\in E.
	\end{equation}
	By \eqref{eq1.9},  for all $x\in E$ and $r\ge 0$, we have
	\begin{eqnarray*}
		\mathrm{Var}_{\delta_{x}}[\langle f,X_{r}\rangle]&=&\int_{0}^{r}T_{r-u}\left(\vartheta[T_{u}f]\right)(x)\mathrm{d}u
		\le
		c^{2}_{1}
		\int_{0}^{r}T_{r-u}\left(\vartheta[\e^{\lambda_1u}\|f\|_{\infty}]\right)(x)\mathrm{d}u\\
		&\le&
		c^{2}_{1}
		\|f\|^{2}_{\infty}\|\vartheta[1]\|_{\infty}\int_{0}^{r}\e^{2 \lambda_1 u}T_{r-u}1(x)\mathrm{d}u\le
		c^{3}_1
		\|f\|^{2}_{\infty}\|\vartheta[1]\|_{\infty}\int_{0}^{r}\e^{\lambda_1 (u+r)}\mathrm{d}u.
	\end{eqnarray*}
	From this, it is easy to show that $(r,x)\mapsto \mathrm{Var}_{\delta_{x}}[\langle f,X_{r}\rangle]$ is bounded on $[0,t]\times E$. Using this, we can easily verify that $x\mapsto \int_{0}^{t}T_{s}\left(\vartheta[\mathrm{Var}_{\delta_{\cdot}}[\langle f,X_{t-s}\rangle]]\right)(x)\mathrm{d}s$ is bounded on $E$.
	Hence by \eqref{eq3.6}, $x\mapsto \pp_{\delta_{x}}\left[\langle f,X_{t}\rangle^{4}\right]$ is bounded on $E$.
\end{proof}

\begin{remark}\label{remark2}
	\rm	
	Lemma \ref{lem3.2} implies that for any $t>0$, $f\in\mathcal{B}^{+}_{b}(E)$ and $p\in [2,4]$, the function $x\mapsto \int_{\mfo}\nu(f)^{p}L_{t}(x,\mathrm{d}\nu)$ is bounded on $E$. In fact, we have
	\begin{align*}
		\int_{\mfo}\nu(f)^{p}L_{t}(x,\mathrm{d}\nu)&\le \|f\|^{p}_{\infty}\int_{\mfo}\nu(1)^{p}L_{t}(x,\mathrm{d}\nu)\\
		&\le\|f\|^{p}_{\infty}\Big(\int_{\nu(1)\le 1}\nu(1)^{2}L_{t}(x,\mathrm{d}\nu)+\int_{\nu(1)>1}\nu(1)^{4}L_{t}(x,\mathrm{d}\nu)\Big).
	\end{align*}
	We also note that $x\mapsto \int_{\mfo}\nu(1)^{2}L_{t}(x,\mathrm{d}\nu)=\mathrm{Var}_{\delta_{x}}[\langle 1,X_{t}\rangle]=\int_{0}^{t}T_{t-s}(\vartheta
	[T_{s}1)](x)\mathrm{d}s$ is bounded on $E$.
	Hence we get the conclusion.
\end{remark}

\begin{proof}[Proof of Lemma \ref{lem3.4} for a $(P_{t},\psi)$-superprocess]
	Fix $x\in E$, $\sigma,s>0$ and $r\ge 0$.
	Let $B:=(B_{t})_{t\ge 0}$
	be a standard Brownian motion independent of $(X_{t})_{t\ge 0}$.
	For $n\in\mathbb{N}$ and $\epsilon>0$,
	define
	$\mathcal{G}_{n}:=\mathcal{F}_{ n\sigma+r}\vee \sigma(B_{s}:\ s\le n\sigma+r)$
	and
	\begin{align*}
		V_{n,\epsilon}:=&\sqrt{\e^{-\lambda_{1}(n\sigma+r+s)}\langle \mathrm{Var}_{\delta_{\cdot}}[\langle \re\big(\e^{-\mathrm{i}(n\sigma+r+s)\im\lambda }g\big),X_{s}\rangle],X_{ n\sigma+r }\rangle+\epsilon^{2}s},\\
		Z_{n,\epsilon}:=&\frac{\e^{(c-\frac{\lambda_{1}}{2})(n\sigma+r+s)}\big(\re W_{n\sigma+r+s}(\lambda,g)-\re W_{n\sigma+r}(\lambda,g)\big)+\epsilon\big(B_{n\sigma+r+s}-B_{n\sigma+r}\big)}{V_{n,\epsilon}}.
	\end{align*}
	Similar to \eqref{step1}, for $t\ge 0$ and $y\in E$, we have
	\begin{align*}
		\mathrm{Var}_{\delta_{y}}[\langle \re\big(\e^{-\mathrm{i}(t+s)\im\lambda }g\big),X_{s}\rangle]=&\int_{0}^{s}T_{r}\big(\vartheta\big[T_{s-r}\big(\re\big(\e^{-\mathrm{i}(t+s)\im\lambda}g\big)\big)\big]\big)(y)\mathrm{d}r\\
		=&\int_{0}^{s}\e^{2c(s-r)}T_{r}\big(\vartheta[\re\big(\e^{-\mathrm{i}(t+r)\im\lambda}g\big)]\big)(y)\mathrm{d}r\\
		=&\frac{1}{2}\e^{2cs}\int_{0}^{s}\e^{-2 c r}T_{r}\big(\vartheta[g,\bar{g}]+\re\big(
		\e^{-2\mathrm{i}(t+r)\im\lambda}
		\vartheta[g]\big)\big)(y)\mathrm{d}r.
	\end{align*}
	Hence,
	\begin{align}\label{lem3.4.1}
		&\e^{-\lambda_{1}(n\sigma+r+s)}\langle \mathrm{Var}_{\delta_{\cdot}}[\langle \re\big(\e^{-\mathrm{i}(n\sigma+r+s)\im\lambda}g\big),X_{s}\rangle],X_{n\sigma+r}\rangle\nonumber\\
		=&\frac{1}{2}\e^{(2 c-\lambda_{1}) s}\big[\e^{-\lambda_{1}(n\sigma+r)}\langle\int_{0}^{s}\e^{-2cu}T_{u}(\vartheta[g,\bar{g}])\mathrm{d}u,X_{n\sigma+r}\rangle\nonumber\\
		&+\re\big(\e^{-2\mathrm{i}(n\sigma+r)\im\lambda }\cdot\e^{-\lambda_{1}(n\sigma+r)}\langle\int_{0}^{s}\e^{-2\lambda u}T_{u}(\vartheta[g])\mathrm{d}u,X_{n\sigma+r}\rangle\big)\big].
	\end{align}
	By the Markov property and Lemma \ref{lem2.2},
	we have
	\begin{eqnarray*}
		&&\pp_{\delta_{x}}\left(\lim_{\mathbb{N}\ni n\to+\infty}\e^{-\lambda_{1}(n\sigma+r)}\langle\int_{0}^{s}\e^{-2cu}T_{u}(\vartheta[g,\bar{g}])\mathrm{d}u,X_{n\sigma+r}\rangle
		=\langle \int_{0}^{s}\e^{-2cu}T_{u}(\vartheta[g,\bar{g}])\mathrm{d}u,\widetilde{\varphi}\rangle W^{\varphi}_{\infty}\right)\\
		&=&\pp_{\delta_{x}}\left(\pp_{X_{r}}\left(\lim_{\mathbb{N}\ni n\to+\infty}\e^{-\lambda_{1}n\sigma}\langle\int_{0}^{s}\e^{-2cu}T_{u}(\vartheta[g,\bar{g}])\mathrm{d}u,
		X_{n\sigma}
		\rangle
		=\langle \int_{0}^{s}\e^{-2cu}T_{u}(\vartheta[g,\bar{g}])\mathrm{d}u,\widetilde{\varphi}\rangle W^{\varphi}_{\infty}\right)\right)\\
		&=&1.
	\end{eqnarray*}
	Note that $\langle \int_{0}^{s}\e^{-2cu}T_{u}(\vartheta[g,\bar{g}])\mathrm{d}u,\widetilde{\varphi}\rangle=\int_{0}^{s}\e^{(\lambda_{1}-2c)u}\langle\vartheta[g,\bar{g}],\widetilde{\varphi}\rangle \mathrm{d}u=
	 \Sigma_{1}(\lambda, g;s)
	$.
	This with the above equation yields that
	$$\lim_{\mathbb{N}\ni n\to+\infty}\e^{-\lambda_{1}(n\sigma+r)}\langle\int_{0}^{s}\e^{-2cu}T_{u}(\vartheta[g,\bar{g}])\mathrm{d}u,X_{n\sigma+r}\rangle=
	 \Sigma_{1}(\lambda, g;s)
	W^{\varphi}_{\infty}\quad\pp_{\delta_{x}}\mbox{-a.s.}$$
	Similarly, we can prove that
	$$\lim_{\mathbb{N}\ni n\to+\infty}\e^{-\lambda_{1}(n\sigma+r)}\langle\int_{0}^{s}\e^{-2\lambda u}T_{u}(\vartheta[g])\mathrm{d}u,X_{n\sigma+r}\rangle=
	 \Sigma_{2}(\lambda, g;s)
	 W^{\varphi}_{\infty}\quad\pp_{\delta_{x}}\mbox{-a.s.}$$
	It then follows from \eqref{lem3.4.1} that, $\pp_{\delta_{x}}$-a.s.,
	\begin{equation}\label{lem3.2.0}
		\lim_{\mathbb{N}\ni n\to+\infty}V_{n,\epsilon}^{2}-\Big[\frac{1}{2}\e^{(2c-\lambda_{1})s}
		\mathfrak{F}_{s}(n\sigma+r)
		W^{\varphi}_{\infty}+\epsilon^{2}s\Big]=0.
	\end{equation}	
	Here $\mathfrak{F}_{s}$ is defined in \eqref{neweq11}.
	Consequently, $\pp_{\delta_{x}}$-a.s.,
	\begin{equation}\label{lem3.2.2}
		\lim_{\mathbb{N}\ni n\to+\infty}\frac{V_{n,\epsilon}^{2}}{\frac{1}{2}\e^{(2c-\lambda_{1})s}
		\mathfrak{F}_{ s	}(n\sigma+r)
W^{\varphi}_{\infty}+\epsilon^{2}s}=1.
	\end{equation}
	Next we shall prove that
	\begin{equation}\label{neweq1}
		\limsup_{\mathbb{N}\ni n\to+\infty}\frac{Z_{n,\epsilon}}{\sqrt{2\log n}}=1 \quad \pp_{\delta_{x}}\mbox{-a.s.}
	\end{equation}
	By Lemma \ref{lem6.1}, it suffices to prove that
	\begin{equation}\label{lem6.3.1}
		\sum_{n=1}^{+\infty}\sup_{y\in\R}\left|\pp_{\delta_{x}}\left(Z_{n,\epsilon}\le y \,|\, \mathcal{G}_{n}\right)-\Phi(y)\right|<+\infty\quad\pp_{\delta_{x}}\mbox{-a.s.}
	\end{equation}
	Let $f_{n,\epsilon}(\theta):=\pp_{\delta_{x}}\left[\e^{\mathrm{i}\theta Z_{n,\epsilon}}|\mathcal{G}_{n}\right]$ for $\theta\in\R$.
	Since $B$ and $X$ are independent, and $V_{n,\epsilon}$ is $\mathcal{G}_{n}$-measurable, we have
	\begin{align}\label{lem6.3.2}
		&f_{n,\epsilon}(\theta)=\e^{-\frac{\epsilon^{2}s}{2}\hat{\theta}^{2}}\pp_{X_{n\sigma+r}}\Big[\exp\Big\{\mathrm{i}\hat{\theta}\e^{-\frac{\lambda_{1}}{2}(n\sigma+r+s)}\big(\langle \re\big(\e^{-\mathrm{i}(n\sigma+r+s)\im\lambda }g\big),X_{s}\rangle\nonumber\\
		&\qquad\qquad-\langle T_{s}\big(\re\big(\e^{-\mathrm{i}(n\sigma+r+s)\im\lambda }g\big)\big),X_{0}\rangle\big)\Big\}\Big]\bigg|_{\hat{\theta}=\frac{\theta}{V_{n,\epsilon}}}.
	\end{align}
	Recall from \eqref{eq3} that for every $\mu\in\mf$, $s\ge 0$ and $\theta \in \R$,
	$$\pp_{\mu}\left[\exp\{\mathrm{i}\theta \left(\langle f,X_{s}\rangle-\langle T_{s}f,X_{0}\rangle\right)\}\right]=\exp\Big\{\langle-\frac{1}{2}\theta ^{2}\mathrm{Var}_{\delta_{\cdot}}[\langle f,X_{s}\rangle]
	-\mathfrak{L}_{s}[-\mathrm{i}\theta  f],\mu\rangle\Big\}.$$
	Inserting to \eqref{lem6.3.2} yields that
	$$f_{n,\epsilon}(\theta)=\exp\Big\{-\frac{1}{2}\theta^{2}-\langle\mathfrak{L}_{s}\left[-\mathrm{i}\frac{\theta}{V_{n,\epsilon}}\e^{-\frac{\lambda_{1}}{2}(n\sigma+r+s)}\re\big(\e^{-\mathrm{i}(n\sigma+r+s)\im\lambda }g\big)\right],X_{ n\sigma+r }\rangle\Big\}.$$
	By Lemma \ref{lemA.3}, we have
	\begin{align}
		&\sup_{y\in\R}\left|\pp_{\delta_{x}}\left(Z_{n,\epsilon}\le y\,|\,\mathcal{G}_{n}\right)-\Phi(y)\right|\nonumber\\
		\le&\frac{1}{\pi}\int_{-\infty}^{+\infty}\left|\frac{f_{n,\epsilon}(\theta)-\e^{-\frac{1}{2}\theta^{2}}}{\theta}\right|\mathrm{d}\theta\nonumber\\
		=&\frac{1}{\pi}\int_{-\infty}^{+\infty}\frac{1}{|\theta|}\e^{-\frac{1}{2}\theta^{2}}\left|\exp\big\{-\langle\mathfrak{L}_{s}\left[-\mathrm{i}\frac{\theta}{V_{n,\epsilon}}\e^{-\frac{\lambda_{1}}{2}
			(n\sigma+r+s)}\re\big(\e^{-\mathrm{i}(n\sigma+r+s)\im\lambda }g\big)\right],X_{ n\sigma+r }\rangle\big\}-1\right|\mathrm{d}\theta.\label{eq3.12}
	\end{align}
	According to
	the inequality
	\begin{equation}\label{eq:inequality2}
		\left|\e^{-z}-1+z-\frac{z^{2}}{2}\right|\le|z|^{2}\left(\frac{|z|}{6}\wedge 1\right)\quad\forall z\in\mathbb{C}^{+},
	\end{equation}
	we have for every $\hat{\theta}\in\R$, $s>0$, $f\in\mathcal{B}_{b}(E)$ and $y\in E$,
	\begin{eqnarray}\label{lem6.3.4}
		|\mathfrak{L}_{s}[-\mathrm{i}\hat{\theta} f](y)|&\le&\int_{\mfo}\left|1-\e^{\mathrm{i}\hat{\theta}\nu(f)}+\mathrm{i}\hat{\theta}\nu(f)-\frac{1}{2}\hat{\theta}^{2}\nu(f)^{2}\right| L_{s}(y,\mathrm{d}\nu)\nonumber\\
		&\le&\frac{|\hat{\theta}|^{3}}{6}\int_{\mfo}\nu(|f|)^{3} L_{s}(y,\mathrm{d}\nu).
	\end{eqnarray}
	Consequently, we have
	\begin{align}\label{lem6.3.5}
		&\left|\langle\mathfrak{L}_{s}\left[-\mathrm{i}\frac{\theta}{V_{n,\epsilon}}\e^{-\frac{\lambda_{1}}{2}(n\sigma+r+s)}\re\big(\e^{-\mathrm{i}(n\sigma+r+s)\im\lambda }g\big)\right],X_{ n\sigma+r }\rangle\right|\nonumber\\
		\le& \frac{|\theta|^{3}}{6V_{n,\epsilon}^{3}}\e^{-\frac{3}{2}\lambda_{1}(n\sigma+r+s)}
		\langle \int_{\mfo}\nu(|g|)^{3} L_{s}(\cdot,\mathrm{d}\nu),X_{ n\sigma+r }\rangle.
	\end{align}
	Recall from Remark \ref{remark2}
	that $y\mapsto \int_{\mfo}\nu(|g|)^{3}L_{s}(y,\mathrm{d}\nu)$ is a bounded 	function on $E$. Thus by
	Lemma \ref{lem2.2},
	$\langle \int_{\mfo}\nu(|g|)^{3} L_{s}(\cdot,\mathrm{d}\nu),X_{ n\sigma+r }\rangle=O(\e^{\lambda_{1} (n\sigma+r) })$ $\pp_{\delta_{x}}$-a.s.
	Note that $V_{n,\epsilon}\ge \epsilon\sqrt{s}>0$ for all $n$.
	It then follows from \eqref{lem6.3.5} that $\left|\langle\mathfrak{L}_{s}\left[-\mathrm{i}\frac{\theta}{V_{n,\epsilon}}\e^{-\frac{\lambda_{1}}{2}(n\sigma+r+s)}\re\big(\e^{-\mathrm{i}(n\sigma+r+s)\im\lambda }g\big)\right],X_{ n\sigma+r }\rangle\right|=|\theta|^{3}O(\e^{-\frac{\lambda_{1}}{2}(n\sigma+r)})$
	$\pp_{\delta_{x}}$-a.s. Using the fact that for
	$z\in\mathbb{C}$, $|\e^{ z}-1|=O(|z|)$ as $|z|\to 0$,
	we get for $|\theta|\leq n$,
	$$\left|\exp\big\{-\langle\mathfrak{L}_{s}\left[-\mathrm{i}\frac{\theta}{V_{n,\epsilon}}
	\e^{-\frac{\lambda_{1}}{2}(n\sigma+r+s)}
	\re\big(\e^{-\mathrm{i}(n\sigma+r+s)\im\lambda }g\big)\right],X_{ n\sigma+r }\rangle\big\}-1\right|=|\theta|^{3}O(\e^{-\frac{\lambda_{1}}{2}(n\sigma+r)})
	\mbox{ as }n\to+\infty.$$
	On the other hand, for $|\theta|>n$,
	since there exists some random variable $\mathcal{C}>0$ such that $V_{n,\epsilon}^2\leq \mathcal{C}$ by \eqref{lem3.2.2},  it holds that
	\[
	\left|f_{n,\epsilon}(\theta)\right|\leq \e^{-\frac{\epsilon^2 s}{2 V_{n,\epsilon}^2}\theta^2}\leq  \e^{-\frac{\epsilon^2 s}{2\mathcal{C}}\theta^2}.
	\]
	Here the first inequality is from \eqref{lem6.3.2}.
	Inserting to \eqref{eq3.12} yields that,
	\begin{align*}
		\sup_{y\in\R}\left|\pp_{\delta_{x}}\big(Z_{n,\epsilon}\le y|\mathcal{G}_{n}\big)-\Phi(y)\right| & \le
		\frac{1}{\pi}\int_{-n}^{n}
		|\theta|^{2}\e^{-\frac{1}{2}\theta^{2}}
		O(\e^{-\frac{\lambda_{1}}{2}(n\sigma+r)}) \mathrm{d}\theta
		+ \frac{1}{\pi}\int_{\{|\theta|>n\}}\frac{\e^{-\frac{1}{2}\theta^{2}}+\e^{-\frac{\epsilon^2 s}{2\mathcal{C}}\theta^2} }{|\theta|} \mathrm{d}\theta  \nonumber\\
		&=
		O(\e^{-\frac{\lambda_{1}}{2}(n\sigma+r)})
		\mbox{ as }n\to+\infty,
	\end{align*}
	which implies \eqref{lem6.3.1}.
	Now it follows by \eqref{lem3.2.2} and \eqref{neweq1} that
	\begin{align}\label{neweq2}
		&\limsup_{\mathbb{N}\ni n\to+\infty}\frac{\e^{(c-\frac{\lambda_{1}}{2})(n\sigma+r+s)}\big(\re W_{n\sigma+r+s}(\lambda,g)-\re W_{n\sigma+r}(\lambda,g)\big)+\epsilon\big(B_{n\sigma+r+s}-B_{n\sigma+r}\big)}{\sqrt{\log n}\sqrt{W^{\varphi}_{\infty}\e^{(2c-\lambda_{1})s}
		\mathfrak{F}_{s}(n\sigma+r)
		+2\epsilon^{2}s}}\nonumber\\
		=&1\quad\pp_{\delta_{x}}\mbox{-a.s.}
	\end{align}
	If $\langle \vartheta[g,\bar{g}],\widetilde{\varphi}\rangle=0$, then $\langle \vartheta[g],\widetilde{\varphi}\rangle=0$, and the above equality yields that
	$$\limsup_{\mathbb{N}\ni n\to+\infty}\frac{\e^{(c-\frac{\lambda_{1}}{2})(n\sigma+r+s)}\big(\re W_{n\sigma+r+s}(\lambda,g)-\re W_{n\sigma+r}(\lambda,g)\big)}{\epsilon\sqrt{2s\log n}}+\frac{B_{n\sigma+r+s}-B_{n\sigma+r}}{\sqrt{2s\log n}}=1\quad\pp_{\delta_{x}}\mbox{-a.s.}$$
	Hence we have, $\pp_{\delta_{x}}$-a.s.,
	$$\limsup_{\mathbb{N}\ni n\to+\infty}\frac{\e^{(c-\frac{\lambda_{1}}{2})(n\sigma+r+s)}\big(\re W_{n\sigma+r+s}(\lambda,g)-\re W_{n\sigma+r}(\lambda,g)\big)}{\sqrt{2s\log n}}\le \epsilon\big(1-\liminf_{\mathbb{N}\ni n\to+\infty}\frac{B_{n\sigma+r+s}-B_{n\sigma+r}}{\sqrt{2s\log n}}\big) =2\epsilon.
	$$
	We conclude the second assertion by letting $\epsilon\to 0$.
	
	Now suppose $\langle \vartheta[g,\bar{g}],\widetilde{\varphi}\rangle>0$. We note that, if $\im\lambda=0$, then
	 $\re\big(\e^{-2\mathrm{i}t\im\lambda }
	 \Sigma_{2}(\lambda,g;s)	\big)= \Sigma_{1}(\lambda,g;s)$. Otherwise if $\im\lambda\not=0$, one can verify by an elementary calculation that $|F_{s}(\lambda_{1}-2\lambda)|<F_{s}(\lambda_{1}-2c)$. Thus,
	 $\big|\re\big(\e^{-2\mathrm{i}t\im\lambda } \Sigma_{2}(\lambda,g;s) \big)\big|	< \Sigma_{1}(\lambda,g;s).$
	In either case, there exist positive constants $c_{1}(s)$ and $c_{2}(s)$, such that
	\begin{align}\label{Def-F}
		c_{1}(s)\le \mathfrak{F}_{s}(t)	\le c_{2}(s)\quad\forall t\ge 0.
	\end{align}
	Note that,
	\begin{align*}
		&\frac{\e^{(c-\frac{\lambda_{1}}{2})(n\sigma+r)}\big(\re W_{n\sigma+r+s}(\lambda,g)-\re W_{n\sigma+r}(\lambda,g)\big)}{\sqrt{\log n}\sqrt{\mathfrak{F}_{s}(n\sigma+r)}}\\
		=&\frac{\e^{(c-\frac{\lambda_{1}}{2})(n\sigma+r+s)}\big(\re W_{n\sigma+r+s}(\lambda,g)-\re W_{n\sigma+r}(\lambda,g)\big)+\epsilon\big(B_{n\sigma+r+s}-B_{n\sigma+r}\big)}{\sqrt{\log n}\sqrt{W^{\varphi}_{\infty}\e^{(2c-\lambda_{1})s}\mathfrak{F}_{s}(n\sigma+r)+2\epsilon^{2}s}}\cdot
		\sqrt{W^{\varphi}_{\infty}+\frac{2\epsilon^{2}s}{\e^{(2c-\lambda_{1})s}\mathfrak{F}_{s}(n\sigma+r)}}\\
		&-\frac{\epsilon\sqrt{2s}}{\sqrt{\e^{(2c-\lambda_{1})s}\mathfrak{F}_{s}(n\sigma+r)}}\cdot \frac{B_{n\sigma+r+s}-B_{n\sigma+r}}{\sqrt{2s\log n}}.
	\end{align*}
	It then follows by \eqref{neweq2} that, $\pp_{\delta_{x}}$-a.s.,
	\begin{align*}
		\sqrt{W^{\varphi}_{\infty}+\frac{2\epsilon^{2}s}{\e^{(2c-\lambda_{1})s}c_{2}(s)}}-\frac{\epsilon\sqrt{2s}}{\sqrt{c_{1}(s)\e^{(2c-\lambda_{1})s}}}\le& \limsup_{\mathbb{N}\ni n\to+\infty}\frac{\e^{(c-\frac{\lambda_{1}}{2})(n\sigma+r)}\big(\re W_{n\sigma+r+s}(\lambda,g)-\re W_{n\sigma+r}(\lambda,g)\big)}{\sqrt{\log n}\sqrt{\mathfrak{F}_{s}(n\sigma+r)}}\\
		\le &\sqrt{W^{\varphi}_{\infty}+\frac{2\epsilon^{2}s}{\e^{(2c-\lambda_{1})s}c_{1}(s)}}+\frac{\epsilon\sqrt{2s}}{\sqrt{c_{1}(s)\e^{(2c-\lambda_{1})s}}}.
	\end{align*}
	The	first conclusion thus follows by letting $\epsilon\to 0$.
	Therefore, we complete the proof.
\end{proof}

\subsubsection{Proof for BMPs}
We will use the following Berry--Esseen inequality (see, e.g. \cite[4.1.b]{LB}).
\begin{lemma}\label{lem:BE inequality}
	Let $\eta_{1},\cdots,\eta_{n}$ be independent random variables with $\mathrm{E}[\eta_{j}]=0$ and $\mathrm{E}[|\eta_{j}|^{3}]<+\infty$, $j=1,\cdots,n$. Then there exists a positive constant $C_{1}$ (independent of the laws of $\eta_{1},\cdots,\eta_{n}$), such that,
	$$\sup_{\R}\Big|\mathrm{P}\left(\frac{\sum_{j=1}^{n}\eta_{j}}{\sqrt{\mathrm{Var}[\sum_{j=1}^{n}\eta_{j}]}}\le y\right)-\Phi(y)\Big|\le C_{1}\frac{\sum_{j=1}^{n}\mathrm{E}[|\eta_{j}|^{3}]}{\mathrm{Var}[\sum_{j=1}^{n}\eta_{j}]^{3/2}}.$$
\end{lemma}
\begin{proof}[Proof of Lemma \ref{lem3.4} for a $(P_{t},G)$-BMP]
	Fix $x\in E$, $\sigma,s>0$ and $r\ge 0$.
	Let $Z:=(Z_{t})_{t\ge 0}$
	be an independent continuous-time binary branching process with branching rate $1$ and $Z_0=1$.
	It is well-known that $W_{t}(Z):=\e^{-t}Z_{t}$ is a nonnegative martingale which converges a.s. and in $L^{2}$ to a limit $W_{\infty}(Z)$, and $W_{\infty}(Z)>0$ a.s.
	
	For $n\in\mathbb{N}$ and $\varepsilon>0$, let $\mathcal{G}_{n}:=\mathcal{F}_{ n\sigma+r}\vee \sigma(Z_{s}:\ s\le n\sigma)$
	and
	\begin{align*}
		Y^{\pm}_{n,\varepsilon}:=&\e^{(c-\frac{\lambda_{1}}{2})(n\sigma+r+s)}\big(\re W_{n\sigma+r+s}(\lambda,g)-\re W_{n\sigma+r}(\lambda,g)\big)\pm \varepsilon\e^{n\sigma/2}\big(W_{(n+1)\sigma}(Z)-W_{n\sigma}(Z)\big),\\
		V_{n,\varepsilon}:=&\sqrt{\mathrm{Var}[Y^{\pm}_{n,\varepsilon}|\mathcal{G}_{n}]}=\sqrt{\e^{-\lambda_{1}(n\sigma+r+s)}\langle \mathrm{Var}_{\delta_{\cdot}}[\langle \re\big(\e^{-\mathrm{i}(n\sigma+r+s)\im\lambda }g\big),X_{s}\rangle],X_{ n\sigma+r }\rangle+\varepsilon^{2}W_{n\sigma}(Z)\mathrm{Var}[W_{\sigma}(Z)]}.
	\end{align*}
	The last equality is from the Markov property and the independence of $X$ and $Z$.
	We adopt the notation $\mathfrak{F}_{s}$ as defined in \eqref{neweq11}.	Applying a computation similar to the one leading to \eqref{lem3.2.2},
	we can show that
	\begin{align}\label{lem3.2.2'}
		&\lim_{\mathbb{N}\ni n\to+\infty}\frac{V_{n,\varepsilon}^{2}}{\frac{1}{2}\e^{(2c-\lambda_{1})s}
			\mathfrak{F}_s(n\sigma+r)
			W^{\varphi}_{\infty}+\varepsilon^{2}W_{\infty}(Z)\mathrm{Var}[W_{\sigma}(Z)]}
		=1\quad\pp_{\delta_{x}}\mbox{-a.s.}
	\end{align}
	Consequently, $\liminf_{\mathbb{N}\ni n\to+\infty}V^{2}_{n,\varepsilon}\ge\varepsilon^{2}W_{\infty}(Z)\mathrm{Var}[W_{\sigma}(Z)]>0$ $\pp_{\delta_{x}}$-a.s.
	
	Let $Z^{\pm}_{n,\varepsilon}:=Y^{\pm}_{n,\varepsilon}/V_{n,\varepsilon}$.
	We shall prove that
	\begin{equation}\label{neweq1'}
		\limsup_{\mathbb{N}\ni n\to+\infty}\frac{Z^{\pm}_{n,\varepsilon}}{\sqrt{2\log n}}=1 \quad \pp_{\delta_{x}}\mbox{-a.s.}
	\end{equation}
	By Lemma \ref{lem6.1}, it suffices to prove that
	\begin{equation}\label{lem6.3.1'}
		\sum_{n=1}^{+\infty}\sup_{y\in\R}\left|\pp_{\delta_{x}}\left(Z^{\pm}_{n,\varepsilon}\le y \,|\, \mathcal{G}_{n}\right)-\Phi(y)\right|<+\infty\quad\pp_{\delta_{x}}\mbox{-a.s.}
	\end{equation}
	We observe that, by the branching property, $W_{n\sigma+r+s}(\lambda,g)-W_{n\sigma+r}(\lambda,g)$ can be written as $$\e^{-\lambda(n\sigma+r)}\sum_{i=1}^{N_{n\sigma+r}}W^{(i)}_{s}(\lambda,g),$$
	where, given $\mathcal{G}_{n}$, $W^{(i)}_{s}(\lambda,g)$ $i=1,\cdots,N_{n\sigma+r}$ are independent random variables, and $W^{(i)}_{s}(\lambda,g)$ is equal in law to $(W_{s}(\lambda,g)-W_{0}(\lambda,g),\pp_{\delta_{x_{i}(n\sigma+r)}})$. Hence we have
	$$Y^{\pm}_{n,\varepsilon}=\e^{(c-\frac{\lambda_{1}}{2})s-\frac{\lambda_{1}}{2}(n\sigma+r)}\sum_{i=1}^{N_{n\sigma+r}}\re\big(\e^{-\mathrm{i}(n\sigma+r)\im\lambda}W^{(i)}_{s}(\lambda,g)\big)
	\pm\varepsilon\e^{n\sigma/2}\big(W_{(n+1)\sigma}(Z)-W_{n\sigma}(Z)\big).$$
	In particular, given $\mathcal{G}_{n}$, $\re\big(\e^{-\mathrm{i}(n\sigma+r)\im\lambda}W^{(i)}_{s}(\lambda,g)\big)$, $i=1,\cdots,N_{n\sigma+r}$ and $W_{(n+1)\sigma}(Z)-W_{n\sigma}(Z)$ are independent random variables, and each of them has mean $0$. We also note that
	$$\pp_{\delta_{x}}\big[\big|\re\big(\e^{-\mathrm{i}(n\sigma+r)\im\lambda}W^{(i)}_{s}(\lambda,g)\big)\big|^{3}|\mathcal{G}_{n}\big]\le\pp_{\delta_{x}}\big[|W^{(i)}_{s}(\lambda,g)|^{3}|\mathcal{G}_{n}
	\big]=\pp_{\delta_{x_{i}(n\sigma+r)}}\big[\big|W_{s}(\lambda,g)-W_{0}(\lambda,g)\big|^{3}\big],$$
	and
	\begin{align*}
		& \e^{n\sigma/2}\big(W_{(n+1)\sigma}(Z)-W_{n\sigma}(Z)\big) = \e^{-n\sigma/2} \sum_{i=1}^{Z_{n\sigma}} \big(W_\sigma^{(i)} (Z)-1 \big).
	\end{align*}
	Here, given $\mathcal{G}_{n}$, $W_\sigma^{(i)}(Z)$ are i.i.d. copies of $W_\sigma (Z)$.
	Let $\mathrm{E}$ denote the expectation operator associated with the binary branching process.
	Thus, by Lemma \ref{lem:BE inequality}, we have
	\begin{align}\label{eq3.12'}
		\sup_{y\in \R}\big|\pp_{\delta_{x}}\big(Z^{\pm}_{n,\varepsilon}\le y |\mathcal{G}_{n}\big)-\Phi(y)\big|
		=&\sup_{y\in \R}\left|\pp_{\delta_{x}}\left(\frac{Y^{\pm}_{n,\varepsilon}}{V_{n,\varepsilon}}\le y |\mathcal{G}_{n}\right)-\Phi(y)\right|\nonumber\\
		\le &\frac{C_{1}}{V^{3}_{n,\varepsilon}}\Big(\e^{3(c-\frac{\lambda_{1}}{2})s-\frac{3\lambda_{1}}{2}(n\sigma+r)}\sum_{i=1}^{N_{n\sigma+r}}\pp_{\delta_{x_{i}(n\sigma+r)}}\big[\big
		|W_{s}(\lambda,g)-W_{0}(\lambda,g)\big|^{3}\big]\nonumber\\
		&+\varepsilon^{3}\e^{-n\sigma/2}W_{n\sigma}(Z)
		\mathrm{E}[|W_{\sigma}(Z)-1	|^3]\Big).
	\end{align}
We note that, under (H2), the evolution equation for the $k$-th moment ($k\in \{1,\cdots,4\}$) of the BMP is given in \cite[Lemma 4.1]{DH1} (see, also \cite[Proposition 9.1]{HK}). Using this equation and the Gronwall's inequality, one can apply an inductive argument, with the base case given by \eqref{new4}, to show that
for every $f\in\mathcal{B}^{+}_{b}(E)$, $t\ge 0$ and $k\in \{1,\cdots,4\}$, the function $x\mapsto\pp_{\delta_{x}}\big[\langle f,X_{t}\rangle^{k}\big]$ is bounded on $E$.
	Hence by Lemma \ref{lem2.2}, $\sum_{i=1}^{N_{n\sigma+r}}\pp_{\delta_{x_{i}(n\sigma+r)}}\big[\big
	|W_{s}(\lambda,g)-W_{0}(\lambda,g)\big|^{3}\big]=\langle \pp_{\delta_{\cdot}}\big[\big
	|W_{s}(\lambda,g)-W_{0}(\lambda,g)\big|^{3}\big],X_{n\sigma+r}\rangle=O(\e^{\lambda_{1}(n\sigma+r)})$ as $n\to+\infty$.
	Recall that $\liminf_{\mathbb{N}\ni n\to+\infty}V_{n,\varepsilon}>0$ a.s.	
	Inserting to \eqref{eq3.12'} yields that,
	$$\sup_{y\in\R}\left|\pp_{\delta_{x}}\big(Z^{\pm}_{n,\varepsilon}\le y|\mathcal{G}_{n}\big)-\Phi(y)\right|=O(\e^{-\frac{\lambda_{1}\wedge 1}{2}n\sigma})\mbox{ as }n\to+\infty\quad\pp_{\delta_{x}}\mbox{-a.s.}$$
	Thus we prove \eqref{lem6.3.1'}.
	
	It follows by \eqref{lem3.2.2'} and \eqref{neweq1'} that
	\begin{align}\label{neweq2'}
		&\limsup_{\mathbb{N}\ni n\to+\infty}\frac{\e^{(c-\frac{\lambda_{1}}{2})(n\sigma+r+s)}\big(\re W_{n\sigma+r+s}(\lambda,g)-\re W_{n\sigma+r}(\lambda,g)\big)\pm
			\varepsilon\e^{n\sigma/2}\big(W_{(n+1)\sigma}(Z)-W_{n\sigma}(Z)\big)}{\sqrt{\log n}\sqrt{W^{\varphi}_{\infty}\e^{(2c-\lambda_{1})s}
				\mathfrak{F}_s(n\sigma+r)
				+2\varepsilon^{2}W_{\infty}(Z)\mathrm{Var}[W_{\sigma}(Z)]}}\nonumber\\
		=&1\quad\pp_{\delta_{x}}\mbox{-a.s.}
	\end{align}
	In particular, by taking $g=0$, we get, almost surely,
	\begin{align}\label{LIL-GW}
		\limsup_{\mathbb{N}\ni n\to+\infty}\slash \liminf_{\mathbb{N}\ni n\to+\infty}\frac{\e^{n\sigma/2}\big(W_{(n+1)\sigma}(Z)-W_{n\sigma}(Z)\big)}{\sqrt{2(\log n) W_{\infty}(Z)\mathrm{Var}[W_{\sigma}(Z)]}}=(+\slash -)1.
	\end{align}
	By combining \eqref{neweq2'}-\eqref{LIL-GW} and taking the limit $\varepsilon\downarrow 0$ similarly as in the end of the superprocess proof, we obtain the desired result.
\end{proof}

\section{Proofs of Theorems \ref{them0} and \ref{them1}}\label{S3}

Recall that $(\lambda, g)$ is an eigenpair and that $c=\re\lambda$.
In this section, we investigate the almost sure
behaviour
of the martingale $W_{t}(\lambda,g)$. We show that there exist three regimes, $c<\lambda_{1}/2$, $c=\lambda_{1}/2$ and $c>\lambda_{1}/2$, such that each regime results in different limiting behaviour.
\subsection{Case I: $c<\lambda_{1}/2$}\label{S3.1}

\begin{lemma}
	Suppose $c<\lambda_{1}/2$. Let
	$$D_{+}\slash D_{-}:=\limsup_{t\to+\infty}\slash \liminf_{t\to+\infty}\frac{e^{(c-\frac{\lambda_{1}}{2})t}\re W_{t}(\lambda,g)}{\sqrt{\log t}}.$$
	Then both $D_{+}$ and $D_{-}$ are finite $\pp_{\delta_{x}}$-a.s.
\end{lemma}

\begin{proof}
	Fix arbitrary $s,\sigma>0$.
	Let
	$$l(s):=\sqrt{\big( \Sigma_1(\lambda,g;s)+|\Sigma_{2}(\lambda,g;s)|
		\big)W^{\varphi}_{\infty}},$$
	and $\hat{\epsilon}_{\sigma}(n):=\e^{(\frac{\lambda_{1}}{2}-c)n\sigma}\sqrt{\log n}$ for $n\in\mathbb{N}$.
	Proposition \ref{lem3.5} implies that, for any $\varepsilon>0$,
	$$\re W_{n\sigma+s}(\lambda,g)-\re W_{n\sigma}(\lambda,g)\le\hat{\epsilon}_{\sigma}(n)(l(s)+\varepsilon)\ \mbox{ eventually }\pp_{\delta_{x}}\mbox{-a.s.}$$
	In particular, by setting $\sigma=s=\delta/2$ for some $\delta>0$, we get
	$$\re W_{(n+1)\delta/2}(\lambda,g)-\re W_{n\delta/2}(\lambda,g)\le
	\hat{\epsilon}_{\delta/2}(n)\Big(
	l\big(\delta/2\big)
	+\varepsilon\Big)\mbox{ eventually }\pp_{\delta_{x}}\mbox{-a.s.}
	$$
	Consequently, we have
	\begin{align}\label{eq3.15}
		\re W_{n\delta}(\lambda,g)-\re W_{n\delta/2}(\lambda,g)&=
		\sum_{q=n}^{2n-1}\left(\re W_{(q+1)\delta/2}(\lambda,g)-\re W_{q\delta/2}(\lambda,g)\right)\nonumber\\	&\le \Big(l\big(\delta/2\big)
		+\varepsilon\Big)\sum_{q=n}^{2n-1}\hat{\epsilon}_{\delta/2}(q), \ \mbox{ eventually }\pp_{\delta_{x}}\mbox{-a.s.}
	\end{align}
	Moreover, it holds by standard calculation that,
	for some constant $c_1(\delta)>0$,
	\begin{align*}
		\sum_{q=n}^{2n-1}\hat{\epsilon}_{\delta/2}(q)&=\sum_{q=n}^{2n-1}\sqrt{\log q}\e^{(\frac{\lambda_{1}}{2}-c)\frac{q\delta}{2}}\\	&=\e^{(\frac{\lambda_{1}}{2}-c)n\delta}\sqrt{\log n}\Big(\sqrt{\frac{\log(2n-1)}{\log n}}\sum_{q=n}^{2n-1}\sqrt{\frac{\log q}{\log (2n-1)}}\e^{-(\frac{\lambda_{1}}{2}-c)(n-\frac{q}{2})\delta}\Big)\\
		&\le c_{1}(\delta)\e^{(\frac{\lambda_{1}}{2}-c)n\delta}\sqrt{\log n}.
	\end{align*}
	The final inequality is from the fact that $c<\lambda_{1}/2$. Inserting this to the right hand side of \eqref{eq3.15} and taking $n\to+\infty$, we get
	\begin{equation}\label{eq3.16}
		\limsup_{\mathbb{N}\ni n\to+\infty}\frac{\e^{-(\frac{\lambda_{1}}{2}-c)n\delta}}{\sqrt{\log n}}\left(\re W_{n\delta}(\lambda,g)-\re W_{n\delta/2}(\lambda,g)\right)<+\infty\quad\pp_{\delta_{x}}\mbox{-a.s.}
	\end{equation}
	We note that for any $\varepsilon>0$,
	\begin{align*}
		& \pp_{\delta_{x}}\left(\left|\e^{-(\frac{\lambda_{1}}{2}-c)n\delta}\re W_{n\delta/2}(\lambda,g)\right|>\varepsilon\right)
		\nonumber\\
		=&\pp_{\delta_{x}}\left(\e^{-(\frac{\lambda_{1}}{2}-c)n\delta}\left|\re W_{n\delta/2}(\lambda,g)-\re g(x)+\re g(x)\right|>\varepsilon\right)\\
		\le& \pp_{\delta_{x}}\left(\e^{-(\frac{\lambda_{1}}{2}-c)n\delta}\left|\re W_{n\delta/2}(\lambda,g)-\re g(x)\right|>\varepsilon-\e^{-(\frac{\lambda_{1}}{2}-c)n\delta}|\re g(x)|\right).
	\end{align*}
	Since $g$ is bounded and $c<\lambda_{1}/2$, we have $\e^{-(\frac{\lambda_{1}}{2}-c)n\delta}\re g(x)\to 0$ as $n\to+\infty$. Thus, for $n$ sufficiently large,
	\begin{align}\label{lem3.9.1}
		\pp_{\delta_{x}}\left(\left|\e^{-(\frac{\lambda_{1}}{2}-c)n\delta}\re W_{n\delta/2}(\lambda,g)\right|>\varepsilon\right)
		&\le\pp_{\delta_{x}}\left(\e^{-(\frac{\lambda_{1}}{2}-c)n\delta}\left|\re W_{n\delta/2}(\lambda,g)-\re g(x)\right|>\frac{1}{2}\varepsilon\right)\nonumber\\
		&\le\frac{4}{\varepsilon^{2}}\e^{-(\lambda_{1}-2c)n\delta}\mathrm{Var}_{\delta_{x}}[\re W_{n\delta/2}(\lambda,g)].
	\end{align}
	Recall from Lemma \ref{lemn2}(i) that
	$$\e^{-(\lambda_{1}-2c)t}\mathrm{Var}_{\delta_{x}}[\re W_{t}(\lambda,g)]-\frac{1}{2}\varphi(x)\Big(\frac{\langle\vartheta[g,\bar{g}],\widetilde{\varphi}\rangle}{\lambda_{1}-2c}+\langle \cc[|g|],\widetilde{\varphi}\rangle+\re\big(\e^{-2\mathrm{i}t\im\lambda}\big(\frac{\langle \vartheta[g],\widetilde{\varphi}\rangle}{\lambda_{1}-2\lambda}
	+\langle \cc[g],\widetilde{\varphi}\rangle
	\big)\big)\Big)\to 0$$
	as $t\to+\infty$.
	Hence we have $\sup_{t>0}\e^{-(\lambda_{1}-2c)t}\mathrm{Var}_{\delta_{x}}[\re W_{t}(\lambda,g)]<+\infty$.
	This together with \eqref{lem3.9.1} yields that
	$$\pp_{\delta_{x}}\left(\left|\e^{-(\frac{\lambda_{1}}{2}-c)n\delta}\re W_{n\delta/2}(\lambda,g)\right|>\varepsilon\right)\lesssim \e^{-(\frac{\lambda_{1}}{2}-c)n\delta}\mbox{ for $n$ sufficiently large,}$$
	which implies that
	$$\sum_{n=1}^{+\infty}\pp_{\delta_{x}}\left(\left|\e^{-(\frac{\lambda_{1}}{2}-c)n\delta}\re W_{n\delta/2}(\lambda,g)\right|>\varepsilon\right)<+\infty.$$
	Therefore,
	by the Borel--Cantelli lemma, we have
	$$\lim_{\mathbb{N}\ni n\to+\infty}\e^{-(\frac{\lambda_{1}}{2}-c)n\delta}\re W_{n\delta/2}(\lambda,g)=0 \quad \pp_{\delta_{x}}\mbox{-a.s.}$$
	Inserting to \eqref{eq3.16} yields that
	\begin{equation}\nonumber
		\limsup_{\mathbb{N}\ni n\to+\infty}\frac{\e^{-(\frac{\lambda_{1}}{2}-c)n\delta}}{\sqrt{\log n}}\re W_{n\delta}(\lambda,g)<+\infty\quad\pp_{\delta_{x}}\mbox{-a.s.}
	\end{equation}
	Now we treat the continuous-time case.
	Recall from \eqref{eq3.11}
	that for any $\varepsilon>0$,
	\begin{align*}
		&\sup_{t\in [n\delta,(n+1)\delta)}\re W_{t}(\lambda,g)-\re W_{n\delta}(\lambda,g)-
		\sqrt{2\mathrm{Var}_{X_{t}}[\re\big(\e^{-\lambda t}W_{(n+1)\delta-t}(\lambda,g)\big)]}\\
		\le&
		\epsilon_{\delta,0}(n)\big(\sqrt{(W^{\varphi}_{n\delta}+\varepsilon)\mathfrak{F}_{\delta}(n\delta)}+\varepsilon\big)
		\mbox{ eventually }\pp_{\delta_{x}}\mbox{-a.s.},
	\end{align*}
	where $\epsilon_{\delta,0}(n)=\e^{-(c-\frac{\lambda_{1}}{2})n\delta}\sqrt{\log n}$
	and $\mathfrak{F}_{\delta}(n\delta)$ is defined in
{\eqref{neweq11}}.
	Let $n(t):=\lfloor t/\delta \rfloor$.
	Therefore, $\pp_{\delta_{x}}$-a.s. when $t$ is large enough,
	\begin{align}\label{eq3.19}
		\frac{\e^{-(\frac{\lambda_{1}}{2}-c)t}\re W_{t}(\lambda,g)}{\sqrt{\log t}}
		\le& \frac{\e^{-(\frac{\lambda_{1}}{2}-c)t}}{\sqrt{\log t}}\Big(\re W_{n(t)\delta}(\lambda,g)+
		\epsilon_{\delta,0}(n(t))\big(\sqrt{(W^{\varphi}_{n(t)\delta}+\varepsilon)\mathfrak{F}_{\delta}(n(t)\delta)}+\varepsilon\big)\nonumber\\
		&+\sqrt{2\mathrm{Var}_{X_{t}}[\re\big(\e^{-\lambda t}W_{(n(t)+1)\delta-t}(\lambda,g)\big)]}\,\Big)\nonumber\\
		=&\e^{-(\frac{\lambda_{1}}{2}-c)(t-n(t)\delta)}\sqrt{\frac{\log n(t)}{\log t}}\frac{\e^{-(\frac{\lambda_{1}}{2}-c)n(t)\delta}\re W_{n(t)\delta}(\lambda,g)}{\sqrt{\log n(t)}}\nonumber\\
		&+
		\frac{\e^{-(\frac{\lambda_{1}}{2}-c)t} \epsilon_{\delta,0}(n(t))}{\sqrt{\log t}}
		\big(\sqrt{(W^{\varphi}_{n(t)\delta}+\varepsilon)\mathfrak{F}_{\delta}(n(t)\delta)}+\varepsilon\big)\nonumber\\
		&+\frac{\e^{-(\frac{\lambda_{1}}{2}-c)t}}{\sqrt{\log t}}\sqrt{2\mathrm{Var}_{X_{t}}[\re\big(\e^{-\lambda t}W_{(n(t)+1)\delta-t}(\lambda,g)\big)]}.
	\end{align}
	Note that $n(t)\sim t$ as $t\to+\infty$. We have
	$$\limsup_{t\to+\infty}
	\frac{\e^{-(\frac{\lambda_{1}}{2}-c)t}\epsilon_{\delta,0}(n(t))}{\sqrt{\log t}}
	\big(\sqrt{(W^{\varphi}_{n(t)\delta}+\varepsilon)\mathfrak{F}_{\delta}(n(t)\delta)}+\varepsilon\big)<+\infty\quad\pp_{\delta_{x}}\mbox{-a.s.}$$
	By \eqref{eq3.26},
	$$\frac{\e^{-(\frac{\lambda_{1}}{2}-c)t}}{\sqrt{\log t}}\sqrt{\mathrm{Var}_{X_{t}}[\re\big(\e^{-\lambda t}W_{(n(t)+1)\delta-t}(\lambda,g)\big)]}
	\to 0
	\mbox{ as }t\to+\infty \quad\pp_{\delta_{x}}\mbox{-a.s.}$$
	Hence we get the desired result by letting $t\to +\infty$ in the right hand side of \eqref{eq3.19}.
	Similarly, by investigating $\e^{(c-\frac{\lambda_{1}}{2})t}\re W_{t}(\lambda,-g)/\sqrt{\log t}$, we can prove $D_{-}>-\infty$ $\pp_{\delta_{x}}$-a.s.
\end{proof}

\begin{lemma}\label{lem:modula}
	Suppose $\{\eta_{t}: t\ge 0\}$ is a family of $\mathbb{C}$-valued random variables in a probability space $(\Omega,\mathcal{H},\mathrm{P})$. If for any $z\in\mathbb{C}$ with $|z|=1$,
	\begin{align}\label{condition}
		\limsup_{t\to +\infty} \re ( z  \eta_t) =1\quad\mathrm{P}\mbox{-a.s.,}
	\end{align}
	then
	\[
	\limsup_{t\to+\infty} |\eta_t| =1\quad\mathrm{P}\mbox{-a.s.}
	\]
\end{lemma}
\begin{proof}
	It is easy to see that $\limsup_{t\to+\infty}|\eta_{t}|\ge 1$ $\mathrm{P}$-a.s.,
	since $|\eta_t| \geq \re(\eta_t)$.
	By taking
	$ z = \pm 1, \pm \mathrm{i}$ in \eqref{condition}, 	we get
	$\limsup_{t\to+\infty} |\eta_t|\leq \limsup_{t\to+\infty} |\re(\eta_t)| + \limsup_{t\to+\infty} |\im(\eta_t)|\leq 2$.
	Let $\mathbb{S}^{1}:=\{z\in\mathbb{C}:\ |z|=1\}$. For each $\varepsilon>0$,
	one can	find a finite set $\{z_{1},\cdots,z_{N}\}\subset \mathbb{S}^{1}$ satisfying that, for each $z\in \mathbb{S}^{1}$, there is some $1\le j\le N$ such that $|z-z_{j}|\le \varepsilon$.
	For $t\ge 0$, let $ z (t)$ be an $\mathbb{S}^{1}$-valued random variable such that $|\eta_t|= \re( z (t) \eta_t)$.
	Then there exists $i(t)\in \{1,\cdots,N\}$,
	such that $| z (t)-  z _{i(t)}| \leq \varepsilon$. Hence, $\mathrm{P}$-a.s.,
	\begin{align*}
		\limsup_{t\to+\infty} |\eta_t| & =	\limsup_{t\to+\infty}  \re( z (t) \eta_t)\leq  \limsup_{t\to+\infty}  \re( z _{i(t)} \eta_t) + \varepsilon \limsup_{t\to+\infty} |\eta_t|\nonumber\\
		& \leq 	\limsup_{t\to\infty}  \max_{1\leq i\leq N}  \re ( z _i \eta_t)  + 2\varepsilon \leq 1+2\varepsilon.
	\end{align*}
	Letting $\varepsilon\downarrow 0$, we get $\limsup_{t\to+\infty}|\eta_{t}|\le 1$ $\mathrm{P}$-a.s.
\end{proof}

\begin{proof}[Proof of Theorem \ref{them1}(i)]
	Suppose $(\lambda,g)$ is an eigenpair with $c=\re\lambda<\lambda_{1}/2$.
	For the first equation of Theorem \ref{them1}(i), we only need to show that
	\begin{equation}\label{them1.1}
		D_{+}\slash D_{-}=(+\slash -)\sqrt{
		\Sigma(\lambda,g)W^{\varphi}_{\infty}
		}\quad\pp_{\delta_{x}}\mbox{-a.s.}
	\end{equation}
	Let $s>0$. We have the following decomposition
	\begin{align}\label{them1.6}
		& \frac{\e^{(c-\frac{\lambda_{1}}{2})t}\big(\re W_{t+s}(\lambda,g)-\re W_{t}(\lambda,g)\big)}
		{\sqrt{\log t}} \nonumber \\&=\frac{\e^{(c-\frac{\lambda_{1}}{2})(t+s)}\re W_{t+s}(\lambda,g)}{\sqrt{\log(t+s)}}\cdot
		\sqrt{\frac{\log(t+s)}{\log t}}\cdot \e^{(\frac{\lambda_{1}}{2}-c)s}  - \frac{\e^{(c-\frac{\lambda_{1}}{2})t}\re W_{t}(\lambda,g)}{\sqrt{\log t}}.
	\end{align}
	In view of  Proposition \ref{lem3.5}, by
	letting $t\to+\infty$ in \eqref{them1.6}, we get that, $\pp_{\delta_{x}}$-a.s.,
	\begin{align*}
		D_{+}\e^{(\frac{\lambda_{1}}{2}-c)s}   -D_{+}\le& \sqrt{
			  (\Sigma_1(\lambda,g;s) +|\Sigma_{2}(\lambda,g;s)|)W^{\varphi}_\infty }
			\le
			 D_{+}\e^{(\frac{\lambda_{1}}{2}-c)s}  -D_{-},
	\end{align*}
	and
	\begin{align*}
		D_{-}\e^{(\frac{\lambda_{1}}{2}-c)s} -D_{+}\le& -\sqrt{
			(\Sigma_1(\lambda,g;s) +|\Sigma_{2}(\lambda,g;s)|)W^{\varphi}_\infty }
			\le
			 D_{-}\e^{(\frac{\lambda_{1}}{2}-c)s} -D_{-}.
	\end{align*}
	It follows from the above inequalities that
	\begin{equation}\label{them1.7}
		\sqrt{W^{\varphi}_{\infty}}\frac{a(s)}{b(s)-1}\Big(1-\frac{2}{b(s)}\Big)\le D_{+}\le \frac{a(s)}{b(s)-1}\sqrt{W^{\varphi}_{\infty}}\quad\pp_{\delta_{x}}\mbox{-a.s.},
	\end{equation}
	where $a(s):=\sqrt{
		(\Sigma_1(\lambda,g;s) +|\Sigma_{2}(\lambda,g;s)|)  }
		$
	and $b(s):=\exp\{(\frac{\lambda_{1}}{2}-c)s\}$. Using the fact that $c<\lambda_{1}/2$, one can easily show that $b(s)\to +\infty$, and
	$$\frac{a(s)}{b(s)-1}\to
	\sqrt{
		\frac{\langle\vartheta[g,\bar{g}],\widetilde{\varphi}\rangle}{\lambda_{1}-2c}+\langle\cc[|g|],\widetilde{\varphi}\rangle+\Big|\frac{\langle\vartheta[g],\widetilde{\varphi}\rangle}{\lambda_{1}-2\lambda}+\langle\cc[g],\widetilde{\varphi}\rangle\Big|}
		 = \sqrt{\Sigma(\lambda,g)}
		\mbox{ as }s\to+\infty.$$
	Hence, by taking $s\to+\infty$ in \eqref{them1.7}, we get
	$$D_{+}=\sqrt{
		 \Sigma(\lambda,g) W^{\varphi}_\infty
		}  \quad \pp_{\delta_{x}}\mbox{-a.s.}$$
	Similarly, one can prove that $D_{-}=-\sqrt{
		  \Sigma(\lambda,g)W^{\varphi}_{\infty}
		}$ $\pp_{\delta_{x}}$-a.s.
	Hence we prove \eqref{them1.1}.
	The second equation of Theorem \ref{them1}(i)
	follows from the first equation and Lemma \ref{lem:modula}.
\end{proof}

\subsection{Case II: $ c =\lambda_{1}/2$}\label{S3.2}

\begin{lemma}\label{lem3.9}
	Suppose $c=\lambda_{1}/2$.
	Then the
	following holds $\pp_{\delta_{x}}$-a.s. for every $x\in E$ and $\sigma>0$.
	$$\limsup_{\mathbb{N}\ni n\to+\infty}\slash \liminf_{\mathbb{N}\ni n\to+\infty}\frac{\re W_{n\sigma}(\lambda,g)}{\sqrt{n\sigma\log\log n}}=(+\slash -)\sqrt{
		\Sigma(\lambda,g)
		W^{\varphi}_{\infty}},$$
where $\Sigma(\lambda,g)=\big(1+1_{\{\im\lambda=0\}}\big)\langle\vartheta[g,\bar{g}],\widetilde{\varphi}\rangle$.
\end{lemma}
\begin{proof}
	Fix $x\in E$ and $\sigma>0$. We only need to prove the ``$\limsup$" part. For $n\ge 1$,
	let
	$Y_{n}:=\re W_{n\sigma}(\lambda,g)-\re W_{(n-1)\sigma}(\lambda,g)$ and $\mathcal{G}_{n}:=\mathcal{F}_{n\sigma}$.
	It is easy to see that $\{Y_{n}:n\ge 1\}$ is a martingale difference sequence with filtration $\{\mathcal{G}_{n}:n\ge 1\}$.
	Let $s^{2}_{n}:=\sum_{k=1}^{n}\pp_{\delta_{x}}\left[Y^{2}_{k}|\mathcal{G}_{k-1}\right]$. Since $c=\lambda_{1}/2$, we have by the Markov property,
	\begin{align}\label{lem3.3.1}
		\pp_{\delta_{x}}[Y^{2}_{k}|\mathcal{G}_{k-1}]&=\mathrm{Var}_{\delta_{x}}[\re W_{k\sigma}(\lambda,g)|\mathcal{F}_{(k-1)\sigma}]\nonumber\\
		&=\e^{-\lambda_{1}(k-1)\sigma}\mathrm{Var}_{X_{(k-1)\sigma}}[\re\big(\e^{-\mathrm{i}(k-1)\sigma\im\lambda}W_{\sigma}(\lambda,g)\big)]\nonumber\\
		&=\frac{1}{2}\e^{-\lambda_{1}(k-1)\sigma}\langle\int_{0}^{\sigma}\e^{-\lambda_{1}s}T_{s}
		\big(\vartheta[g,\bar{g}]+\re\big(\e^{-2\mathrm{i}((k-1)\sigma+s)\im\lambda}\vartheta[g]\big)\big)\mathrm{d}s\nonumber\\
		&\quad+\e^{-\lambda_{1}\sigma}T_{\sigma}(\cc[|g|])-
		\cc[|g|]
		+\re\big(\e^{-2\mathrm{i}(k-1)\sigma \im\lambda}\big(\e^{-2\lambda\sigma}T_{\sigma}(\cc[g])-\cc[g]\big)\big)
		,X_{(k-1)\sigma}\rangle.\quad
	\end{align}
	On one hand,
	if $\im\lambda=0$, then $\pp_{\delta_{x}}[Y^{2}_{k}|\mathcal{G}_{k-1}]=\e^{-\lambda_{1}(k-1)\sigma}\langle\int_{0}^{\sigma}\e^{-\lambda_{1}s}T_{s}
	\big(\vartheta[g]\big)\mathrm{d}s+\e^{-\lambda_{1}\sigma}T_{\sigma}(\cc[|g|])-\cc[|g|],X_{(k-1)\sigma}\rangle$.
	Thus by
	Lemma \ref{lem2.2},
	$$\lim_{\mathbb{N}\ni k\to+\infty}\pp_{\delta_{x}}[Y^{2}_{k}|\mathcal{G}_{k-1}]=\langle\int_{0}^{\sigma}\e^{-\lambda_{1}s}T_{s}
	\big(\vartheta[g]\big)\mathrm{d}s,\widetilde{\varphi}\rangle W^{\varphi}_{\infty}=\sigma\langle\vartheta[g],\widetilde{\varphi}\rangle W^{\varphi}_{\infty}\quad\pp_{\delta_{x}}\mbox{-a.s.}$$
	Consequently, we have
	\begin{equation}\label{lem3.10.2}
		\lim_{\mathbb{N}\ni n\to+\infty}\frac{s^{2}_{n}}{n}=\lim_{\mathbb{N}\ni n\to+\infty}\frac{\sum_{k=1}^{n}\pp_{\delta_{x}}[Y^{2}_{k}|\mathcal{G}_{k-1}]}{n}
		=\sigma\langle\vartheta[g],\widetilde{\varphi}\rangle W^{\varphi}_{\infty}\quad\pp_{\delta_{x}}\mbox{-a.s.}
	\end{equation}
	On the other hand, if $\im\lambda\not=0$,
	then \eqref{lem3.3.1} yields that
	\begin{align}
		&\pp_{\delta_{x}}[Y^{2}_{k}|\mathcal{G}_{k-1}]\nonumber\\
		=&\frac{1}{2}\e^{-\lambda_{1}(k-1)\sigma}\langle\int_{0}^{\sigma}\e^{-\lambda_{1}s}
		T_{s}(\vartheta[g,\bar{g}])\mathrm{d}s+\e^{-\lambda_{1}\sigma}T_{\sigma}(\cc[|g|])-\cc[|g|],X_{(k-1)\sigma}\rangle\nonumber\\
		&+\frac{1}{2}\re\big(\e^{-2\mathrm{i}(k-1)\sigma\im\lambda}\cdot\e^{-\lambda_{1}(k-1)\sigma} \langle\int_{0}^{\sigma}\e^{-(\lambda_{1}+2\mathrm{i}\im\lambda)s}T_{s}(\vartheta[g])\mathrm{d}s+\e^{-(\lambda_{1}+2\mathrm{i}\im\lambda)\sigma}T_{\sigma}(\cc[g])-\cc[g],X_{(k-1)\sigma}\rangle\big)\nonumber\\
		=:&\frac{1}{2}I(k,\sigma)+\frac{1}{2}\re\big(\e^{-2\mathrm{i}(k-1)\sigma\im\lambda}II(k,\sigma)\big).\label{lem3.10.8}
	\end{align}
	Since by Lemma \ref{lem2.2},
	\begin{align*}
		\lim_{\mathbb{N}\ni k\to+\infty}II(k,\sigma)=&\langle\int_{0}^{\sigma}\e^{-(\lambda_{1}+2\mathrm{i}\im\lambda)s}T_{s}(\vartheta[g])\mathrm{d}s+\e^{-(\lambda_{1}+2\mathrm{i}\im\lambda)\sigma}T_{\sigma}(\cc[g])-\cc[g],\widetilde{\varphi}\rangle W^{\varphi}_{\infty}\\
		=&\Big(\frac{1-\e^{-2\mathrm{i}\sigma\im\lambda}}{2\mathrm{i}\im\lambda}\langle\vartheta[g],\widetilde{\varphi}\rangle
		+\big(\e^{-2\mathrm{i}\sigma\im\lambda}-1\big)\langle\cc[g],\widetilde{\varphi}\rangle \Big)W^{\varphi}_{\infty}\quad\pp_{\delta_{x}}\mbox{-a.s.},
	\end{align*}
	we have, $\pp_{\delta_{x}}$-a.s.,
	\begin{align}\label{lem3.10.1}
		&\frac{1}{n}
		\sum_{k=1}^{n}
		\re\Big(\e^{-2\mathrm{i}(k-1)\sigma\im\lambda}\big(II(k,\sigma)-\big(\frac{1-\e^{-2\mathrm{i}\sigma\im\lambda}}{2\mathrm{i}\im\lambda}\langle\vartheta[g],\widetilde{\varphi}\rangle
		+\big(\e^{-2\mathrm{i}\sigma\im\lambda}-1\big)\langle\cc[g],\widetilde{\varphi}\rangle\big)W^{\varphi}_{\infty}\big)\Big)\nonumber\\
		\le&\frac{1}{n}\sum_{k=1}^{n}\big|II(k,\sigma)-\big(\frac{1-\e^{-2\mathrm{i}\sigma\im\lambda}}{2\mathrm{i}\im\lambda}\langle\vartheta[g],\widetilde{\varphi}\rangle
		+\big(\e^{-2\mathrm{i}\sigma\im\lambda}-1\big)\langle\cc[g],\widetilde{\varphi}\rangle\big)W^{\varphi}_{\infty}\big|
		\to
		0\mbox{ as }n\to+\infty.
	\end{align}
	A simple calculation gives that
	\begin{align*}
		&\frac{1}{n}\sum_{k=1}^{n}\re\Big(\e^{-2\mathrm{i}(k-1)\sigma\im\lambda}\big(\frac{1-\e^{-2\mathrm{i}\sigma\im\lambda}}{2\mathrm{i}\im\lambda}\langle\vartheta[g],\widetilde{\varphi}\rangle
		+\big(\e^{-2\mathrm{i}\sigma\im\lambda}-1\big)\langle\cc[g],\widetilde{\varphi}\rangle
		\big)\Big)\\
		&=\frac{1}{n}\re\Big(\sum_{k=1}^{n}\frac{\e^{-2\mathrm{i}(k-1)\sigma\im\lambda}-\e^{-2\mathrm{i}k\sigma\im\lambda}}{2\mathrm{i}\im\lambda}\langle\vartheta[g],\widetilde{\varphi}\rangle
		+\big(\e^{-2\mathrm{i}k\sigma\im\lambda}-\e^{-2\mathrm{i}(k-1)\sigma\im\lambda}\big)\langle\cc[g],\widetilde{\varphi}\rangle
		\Big)\\
		&=\frac{1}{n}\re\Big(\frac{1-\e^{-2\mathrm{i}n\sigma\im\lambda}}{2\mathrm{i}\im\lambda}\langle\vartheta[g],\widetilde{\varphi}\rangle
		+\big(\e^{-2\mathrm{i}n\sigma\im\lambda}-1\big)\langle\cc[g],\widetilde{\varphi}\rangle
		\Big)\to 0\mbox{ as }n\to+\infty.
	\end{align*}
	This together with \eqref{lem3.10.1} yields that
	\begin{equation}\label{neweq12}
		\lim_{\mathbb{N}\ni n\to+\infty}\frac{1}{n}\sum_{k=1}^{n}\re\big(\e^{-2\mathrm{i}(k-1)\sigma\im\lambda}
		II(k,\sigma)
		\big)=0 \quad\pp_{\delta_{x}}\mbox{-a.s.}
	\end{equation}
	Hence if $\im\lambda\not=0$, we have by \eqref{lem3.10.8},
	\begin{align}\label{lem3.10.3}
		&\lim_{\mathbb{N}\ni n\to+\infty}\frac{s^{2}_{n}}{n}=\lim_{\mathbb{N}\ni n\to+\infty}\frac{1}{2}\frac{\sum_{k=1}^{n}\e^{-\lambda_{1}(k-1)\sigma}
			\langle\int_{0}^{\sigma}\e^{-\lambda_{1}s}T_{s}(\vartheta[g,\bar{g}])\mathrm{d}s+\e^{-\lambda_{1}\sigma}T_{\sigma}(\cc[|g|])-\cc[|g|],X_{(k-1)\sigma}\rangle}{n}\nonumber\\
		&=\frac{1}{2}\sigma
		\langle\vartheta[g,\bar{g}],\widetilde{\varphi}\rangle W^{\varphi}_{\infty}\quad\pp_{\delta_{x}}\mbox{-a.s.}
	\end{align}
	Note that $\re W_{n\sigma}(\lambda,g)=\re W_{0}(\lambda,g)+\sum_{k=1}^{n}Y_{n}$. In view of Lemma \ref{lemA.3.2} in the Appendix,
	the desired conclusion follows
	from \eqref{lem3.10.2} and \eqref{lem3.10.3},
	provided that
	\begin{equation}\label{eq3.20}
		\sup_{n\ge 1}\pp_{\delta_{x}}\left[Y^{4}_{n}|\mathcal{G}_{n-1}\right]<+\infty\quad\pp_{\delta_{x}}\mbox{-a.s.}
	\end{equation}
We state this as a lemma below. Once established, it will complete the proof of the current lemma.
\end{proof}

	\begin{lemma}\label{lem-eq3.20}
		For $\sigma>0$ and $n\geq 1$, set
$Y_{n}:=\re W_{n\sigma}(\lambda,g)-\re W_{(n-1)\sigma}(\lambda,g)$ and $\mathcal{G}_{n}:=\mathcal{F}_{n\sigma}$.
Then, for every $x\in E$,
		\[
		\sup_{n\ge 1}\pp_{\delta_{x}}\left[Y^{4}_{n}|\mathcal{G}_{n-1}\right]<+\infty\quad\pp_{\delta_{x}}\mbox{-a.s.}
		\]
	\end{lemma}

We will prove this lemma separately for BMPs and superprocesses.

\noindent\textit{Proof of Lemma \ref{lem-eq3.20} for a $(P_{t},G)$-BMP:}
	We shall use the following simple result: If $\eta_{1},\cdots,\eta_{n}$ are independent random variables with $\mathrm{E}[\eta_{j}]=0$, $j=1,\cdots,n$, then
	\begin{equation}\label{eq:i.r.v.}
		\mathrm{E}\Big[\big(\sum_{j=1}^{n}\eta_{j}\big)^{4}\Big]=\sum_{j=1}^{n}\mathrm{E}[\eta^{4}_{j}]+3\sum_{i,j=1;i\not=j}^{n}\mathrm{E}[\eta^{2}_{i}]\mathrm{E}[\eta^{2}_{j}].
	\end{equation}
	We note that, by the branching property,
	$$W_{n\sigma}(\lambda,g)-W_{(n-1)\sigma}(\lambda,g)=\e^{-\lambda(n-1)\sigma}\sum_{k=1}^{N_{(n-1)\sigma}}W^{(k)}_{\sigma}(\lambda,g),$$
	where, given $\mathcal{F}_{(n-1)\sigma}$, $W^{(k)}_{\sigma}(\lambda,g)$, $k=1,\cdots,N_{(n-1)\sigma}$ are independent random variables, and $W^{(k)}_{\sigma}(\lambda,g)$ is equal in law to $\big(W_{\sigma}(\lambda,g)-W_{0}(\lambda,g),\pp_{\delta_{x_{k}((n-1)\sigma)}}\big)$. Thus, by \eqref{eq:i.r.v.}, we have
	\begin{align*}
		&\pp_{\delta_{x}}\big[Y^{4}_{n}\,|\,\mathcal{G}_{n-1}\big]\\
		=&\pp_{\delta_{x}}\big[\big(\re W_{n\sigma}(\lambda,g)-\re W_{(n-1)\sigma}(\lambda,g)\big)^{4}\,|\,\mathcal{F}_{(n-1)\sigma}\big]\\
		=&\sum_{k=1}^{N_{(n-1)\sigma}}\pp_{\delta_{x}}\big[
		\re\big(\e^{-\lambda(n-1)\sigma}W^{(k)}_{\sigma}(\lambda,g)\big)^{4}
		\,|\,\mathcal{F}_{(n-1)\sigma}\big]\\
		&+3\sum_{k,j=1;k\not=j}^{N_{(n-1)\sigma}}\pp_{\delta_{x}}\big[
		\re\big(\e^{-\lambda(n-1)\sigma}W^{(k)}_{\sigma}(\lambda,g)\big)^{2}
		\,|\,\mathcal{F}_{(n-1)\sigma}\big]\cdot\pp_{\delta_{x}}
		\big[
		\re\big(\e^{-\lambda(n-1)\sigma}W^{(j)}_{\sigma}(\lambda,g)\big)^{2}
		\,|\,\mathcal{F}_{(n-1)\sigma}\big].
	\end{align*}
	By the inequality $\sum_{i,j=1;i\not=j}^{n}\mathrm{E}[\eta^{2}_{i}]\mathrm{E}[\eta^{2}_{j}]\leq \big(
	\sum_{i=1}^{n}
	\mathrm{E}[\eta^{2}_{i}] \big)^2$,
	we can continue the computation to get
	\begin{align}
		&\pp_{\delta_{x}}\big[Y^{4}_{n}\,|\,\mathcal{G}_{n-1}\big]\nonumber\\
		=&\e^{-2\lambda_{1}(n-1)\sigma}\sum_{k=1}^{N_{(n-1)\sigma}}\pp_{\delta_{x_{k}(n-1)\sigma}}\big[
		\re\big(\e^{-\mathrm{i}(n-1)\sigma\im\lambda}(W_{\sigma}(\lambda,g)-W_{0}(\lambda,g))\big)^{4}
		\big]\nonumber\\
		&+3\e^{-2\lambda_{1}(n-1)\sigma}\sum_{k,j=1;k\not=j}^{N_{(n-1)\sigma}}\pp_{\delta_{x_{k}(n-1)\sigma}}\big[
		\re\big(\e^{-\mathrm{i}(n-1)\sigma\im\lambda}(W_{\sigma}(\lambda,g)-W_{0}(\lambda,g))\big)^{2}
		\big]\nonumber\\
		&\qquad\qquad\qquad\cdot\pp_{\delta_{x_{k}(n-1)\sigma}}\big[
		\re\big(\e^{-\mathrm{i}(n-1)\sigma\im\lambda}(W_{\sigma}(\lambda,g)-W_{0}(\lambda,g))\big)^{2}
		\big]\nonumber\\
		\le& \e^{-2\lambda_{1}(n-1)\sigma}\langle \pp_{\delta_{\cdot}}\big[\big|W_{\sigma}(\lambda,g)-W_{0}(\lambda,g)\big|^{4}\big],X_{(n-1)\sigma}\rangle\nonumber\\
		&+3\Big(\e^{-\lambda_{1}(n-1)\sigma}\langle\pp_{\delta_{\cdot}}\big[|W_{\sigma}(\lambda,g)-W_{0}(\lambda,g)|^{2}\big],X_{(n-1)\sigma}\rangle\Big)^{2}.\label{eq3.20.2}
	\end{align}
Recall from the argument below \eqref{eq3.12'} that
for every $f\in\mathcal{B}^{+}_{b}(E)$, $t>0$ and $k\in \{1,\cdots,4\}$, $x\mapsto \pp_{\delta_{x}}[\langle f,X_{t}\rangle^{k}]$ is bounded on $E$. This further implies that, $x\mapsto \pp_{\delta_{x}}\big[|W_{\sigma}(\lambda,g)-W_{0}(\lambda,g)|^{k}\big]$, $k=2,4$,
	are bounded on $E$. Thus, by Lemma \ref{lem2.2}, both terms in the right hand side of \eqref{eq3.20.2} converge a.s. to some finite limits as $n\to+\infty$.
Hence we complete the proof. \qed

\medskip

\noindent\textit{Proof of Lemma \ref{lem-eq3.20} for a $(P_{t},\psi)$-superprocess:}
	By the Markov property, we have
	\begin{align}\label{lem3.10.7}
		&\pp_{\delta_{x}}[Y^{4}_{n}|\mathcal{G}_{n-1}]\nonumber
		\\=&\pp_{\delta_{x}}[\big(\re W_{n\sigma}(\lambda,g)-\re W_{(n-1)\sigma}(\lambda,g)\big)^{4}|\mathcal{F}_{(n-1)\sigma}]\nonumber\\
		=&\e^{-2\lambda_{1}(n-1)\sigma}\pp_{X_{(n-1)\sigma}}\big[\big(\re\big(\e^{-\mathrm{i}(n-1)\sigma\im\lambda}(W_{\sigma}(\lambda,g)
		-W_{0}(\lambda,g))\big)\big)^{4}\big]\nonumber\\
		=&\e^{-2\lambda_{1}(n-1)\sigma}\pp_{X_{(n-1)\sigma}}\big[\big(\cos((n-1)\sigma\im\lambda)(\re W_{\sigma}(\lambda,g)-\re W_{0}(\lambda,g))\nonumber\\
		&+\sin((n-1)\sigma\im\lambda)(\im W_{\sigma}(\lambda,g)-\im W_{0}(\lambda,g))\big)^{4}\big]\nonumber\\
		\le&8\e^{-2\lambda_{1}(n-1)\sigma}\big[\pp_{X_{(n-1)\sigma}}\big[\big(\re W_{\sigma}(\lambda,g)-\re W_{0}(\lambda,g)\big)^{4}\big]+\pp_{X_{(n-1)\sigma}}\big[\big(\im W_{\sigma}(\lambda,g)-\im W_{0}(\lambda,g)\big)^{4}\big]\big]. \quad
	\end{align}
	We first deal with $\e^{-2\lambda_{1}(n-1)\sigma}\pp_{X_{(n-1)\sigma}}\big[\big(\re W_{\sigma}(\lambda,g)-\re W_{0}(\lambda,g)\big)^{4}\big]$.
	It is proved in \cite{PY} that
	$\re W_{t}(\lambda,g)-\re W_{0}(\lambda,g)$ is an $\mathcal{F}_{t}$-martingale with quadratic variation
	$$\int_{0}^{t}\e^{-\lambda_{1}s}\langle 2b
	\re\big(\e^{-\mathrm{i}s\im\lambda}g\big)^{2}
	,X_{s}\rangle \mathrm{d}s+\sum_{s\le t}\e^{-\lambda_{1}s}\langle \re\big(\e^{-\mathrm{i}s\im\lambda}g\big),\triangle X_{s}\rangle^{2}1_{\{\triangle X_{s}\not=0\}}.$$
	Here we use the notation $\triangle X_{s}=X_{s}-X_{s-}$ for the jump of $X$ at time $s$.
	Hence by the Burkholder--Davis--Gundy inequality, we have
	\begin{align}
		&\pp_{X_{(n-1)\sigma}}\left[\big(\re W_{\sigma}(\lambda,g)-\re W_{0}(\lambda,g)\big)^{4}\right]\nonumber\\
		\lesssim&\pp_{X_{(n-1)\sigma}}\left[\left(\int_{0}^{\sigma}\e^{-\lambda_{1}s}\langle 2b
		\re\big(\e^{-\mathrm{i}s\im\lambda}g\big)^{2}
		,X_{s}\rangle \mathrm{d}s+\sum_{s\le \sigma}\e^{-\lambda_{1}s}\langle \re\big(\e^{-\mathrm{i}s\im\lambda}g\big),\triangle X_{s}\rangle^{2}1_{\{\triangle X_{s}\not=0\}}\right)^{2}\right]\nonumber\\
		\lesssim&\pp_{X_{(n-1)\sigma}}\left[\left(\int_{0}^{\sigma}\e^{-\lambda_{1}s}\langle 2b
		\re\big(\e^{-\mathrm{i}s\im\lambda}g\big)^{2}
		,X_{s}\rangle \mathrm{d}s\right)^{2}\right]
		\nonumber\\ &\qquad \qquad +
		\pp_{X_{(n-1)\sigma}}\left[\left(\sum_{s\le \sigma}\e^{-\lambda_{1}s}\langle \re\big(\e^{-\mathrm{i}s\im\lambda}g\big),\triangle X_{s}\rangle^{2}1_{\{\triangle X_{s}\not=0\}}\right)^{2}\right]\nonumber\\
		\le&\pp_{X_{(n-1)\sigma}}\left[\left(\int_{0}^{\sigma}\e^{-\lambda_{1}s}\langle
		2b|g|^2 ,X_{s}\rangle \mathrm{d}s\right)^{2}\right]+\pp_{X_{(n-1)\sigma}}\left[\left(\sum_{s\le \sigma}\e^{-\lambda_{1}s}\langle |g|,\triangle X_{s}\rangle^{2}1_{\{\triangle X_{s}\not=0\}}\right)^{2}\right].
		\label{eq3.24}
	\end{align}
	By repeating the computation in the proof of \cite[Lemma A.3]{Y25}, one can show that, for any $\mu\in\mf$, $t>0$ and nonnegative measurable function $(s,y)\mapsto f_{s}(y)$ on $\R^{+}\times E$,
	\begin{equation}\label{lem3.10.4}
		\pp_{\mu}\left[\Big(\int_{0}^{t}\langle f_{s},X_{s}\rangle \mathrm{d}s\Big)^{2}\right]
		=\langle \int_{0}^{t}T_{s}f_{s}\mathrm{d}s,\mu\rangle^{2}+\langle \int_{0}^{t}T_{s}\big(\vartheta[G(t-s,s,f)]\big)\mathrm{d}s,\mu\rangle,
	\end{equation}
	where $G(u,v,f)(y):=\int_{0}^{u}T_{r}f_{v+r}(y)\mathrm{d}r$ for $u,v\ge 0$ and $y\in E$, and
	\begin{align}\label{lem3.10.5}
		\pp_{\mu}\left[\Big(\sum_{s\le t}\langle f_{s},\triangle X_{s}\rangle^{2}1_{\{\triangle X_{s}\not=0\}}\Big)^{2}\right]
		\le&2\int_{0}^{t}\langle T_{s}
		\tilde{f}_{s}
		\mu\rangle\cdot\langle T_{s}(G(t-s,s,
		\tilde{f}
		)),\mu\rangle \mathrm{d}s\nonumber\\
		&+2\int_{0}^{t}\mathrm{d}s\int_{0}^{s}\langle T_{s-r}\big(\vartheta[T_{r}
		\tilde{f}_{s}
		,T_{r}(G(t-s,s,
		\tilde{f}
		))]\big),\mu\rangle \mathrm{d}r\nonumber\\
		&+\langle \int_{0}^{t}T_{s}\big(\vartheta[G(t-s,s,
		\tilde{f}
		)]\big)\mathrm{d}s,\mu\rangle,
	\end{align}
	where $
	\tilde{f}_{s}(y)
	:=\int_{\mfo}\nu(f_{s})^{2}H(y,\mathrm{d}\nu)$ for $(s,y)\in\R^{+}\times E$.
	From \eqref{lem3.10.4}, we have
	\begin{align}\label{lem3.10.6}
		&\e^{-2\lambda_{1}(n-1)\sigma}\pp_{X_{(n-1)\sigma}}\left[\left(\int_{0}^{\sigma}\e^{-\lambda_{1}s}
		\langle 2b|g|^2,X_{s}\rangle \mathrm{d}s\right)^{2}\right]\nonumber\\
		=&\big(\e^{-\lambda_{1}(n-1)\sigma}\langle \int_{0}^{\sigma}2\e^{-\lambda_{1}s}T_{s}(b|g|^{2})\mathrm{d}s,X_{(n-1)\sigma}\rangle\big)^{2}\nonumber\\
		&+\e^{-2\lambda_{1}(n-1)\sigma}\langle\int_{0}^{\sigma}T_{s}\big(\vartheta[G(\sigma-s,s,2\e^{-\lambda_{1}s}b|g|^{2})]\big)\mathrm{d}s,X_{(n-1)\sigma}\rangle.
	\end{align}
	Using the boundedness of $|g|$ and
	the boundedness of $x\mapsto \vartheta[1](x)=2b(x)+\int_{\mfo}\nu(1)^{2}H(x,\mathrm{d}\nu)$ on $E$,
	one can show that both $y\mapsto \int_{0}^{\sigma}2\e^{-\lambda_{1}s}T_{s}(b|g|^{2})\mathrm{d}s$ and $y\mapsto \int_{0}^{\sigma}T_{s}\big(\vartheta[G(\sigma-s,s,2\e^{-\lambda_{1}s}b|g|^{2})]\big)(y)\mathrm{d}s$ are bounded functions on $E$. It then follows by
	Lemma \ref{lem2.2}
	that, the right hand side of \eqref{lem3.10.6} converges $\pp_{\delta_{x}}$-a.s. to some finite random variable as $n\to+\infty$. Hence we get
	$$\sup_{n\ge 1}\e^{-2\lambda_{1}(n-1)\sigma}\pp_{X_{(n-1)\sigma}}\left[\left(\int_{0}^{\sigma}\e^{-\lambda_{1}s}\langle 2b|g|^{2},X_{s}\rangle \mathrm{d}s\right)^{2}\right]<+\infty \quad\pp_{\delta_{x}}\mbox{-a.s.}$$
	Similarly, using \eqref{lem3.10.5} and
	Lemma \ref{lem2.2},
	one can show that
	$$\sup_{n\ge 1}\e^{-2\lambda_{1}(n-1)\sigma}\pp_{X_{(n-1)\sigma}}\left[\big(\sum_{s\le \sigma}\e^{-\lambda_{1}s}\langle |g|^{2},\triangle X_{s}\rangle^{2}1_{\{\triangle X_{s}\not=0\}}\big)^{2}\right]<+\infty\quad\pp_{\delta_{x}}\mbox{-a.s.}$$
	Thus by \eqref{eq3.24}, we have proved that
	$$\sup_{n\ge 1}\e^{-2\lambda_{1}(n-1)\sigma}\pp_{X_{(n-1)\sigma}}\left[\big(\re W_{\sigma}(\lambda,g)-\re W_{0}(\lambda,g)\big)^{4}\right]<+\infty\quad\pp_{\delta_{x}}\mbox{-a.s.}$$
	Similarly, by investigating $\pp_{X_{(n-1)\sigma}}\left[\big(\re W_{\sigma}(\lambda,-\mathrm{i}g)-\re W_{0}(\lambda,-\mathrm{i}g)\big)^{4}\right]$, we can prove that
	$$\sup_{n\ge 1}\e^{-2\lambda_{1}(n-1)\sigma}\pp_{X_{(n-1)\sigma}}\left[\big(\im W_{\sigma}(\lambda,g)-\im W_{0}(\lambda,g)\big)^{4}\right]<+\infty\quad\pp_{\delta_{x}}\mbox{-a.s.}$$
Therefore, combining \eqref{lem3.10.7} and the above inequalities, we complete the proof of this lemma in the superprocess setting.
\qed

\begin{proof}[Proof of Theorem \ref{them1}(ii)]
	Suppose $(\lambda,g)$ is an eigenpair with
	$c=\re\lambda=\lambda_{1}/2$.
	Let $\sigma>0$. For $t>0$, let $n(t):=\lfloor t/\sigma\rfloor$. Then we have
	\begin{equation}\label{eq3.29}
		\frac{\re W_{t}(\lambda,g)}{\sqrt{t\log\log t}}=\frac{\re W_{t}(\lambda,g)-\re W_{n(t)\sigma}(\lambda,g)}{\sqrt{t\log\log t}}+\frac{\sqrt{n(t)\sigma\log \log n(t)}}{\sqrt{t\log\log t}}\cdot\frac{\re W_{n(t)\sigma}(\lambda,g)}{\sqrt{n(t)\sigma\log\log n(t)}}.
	\end{equation}
	It follows from \eqref{lem3.5.3} that
	$$\lim_{\sigma\to 0}\liminf_{t\to+\infty}\frac{\re W_{t}(\lambda,g)-\re W_{n(t)\sigma}(\lambda,g)}{\sqrt{t\log\log t}}=\lim_{\sigma\to 0}\limsup_{t\to+\infty}\frac{\re W_{t}(\lambda,g)-\re W_{n(t)\sigma}(\lambda,g)}{\sqrt{t\log\log t}}=0\quad\pp_{\delta_{x}}\mbox{-a.s.}$$
	Combining the above limit with Lemma \ref{lem3.9} and taking first $t\to+\infty$ and then $\sigma \to 0$ in \eqref{eq3.29}, we conclude that
	$$\limsup_{t\to+\infty}\frac{\re W_{t}(\lambda,g)}{\sqrt{t\log\log t}}
	 = \sqrt{\Sigma(\lambda,g)W^{\varphi}_\infty}
	\quad\pp_{\delta_{x}}\mbox{-a.s.}$$
	The second equation of Theorem \ref{them1}(ii) follows from the first equation and Lemma \ref{lem:modula}. This completes the proof of Theorem \ref{them1} (ii).
\end{proof}

\subsection{Case III: $c>\lambda_{1}/2$}\label{S3.3}

\noindent\textit{Proof of Theorem \ref{them0}:}
We note that
$$\pp_{\delta_{x}}[|W_{t}(\lambda,g)|^{2}]=\pp_{\delta_{x}}[
\big|\re W_{t}(\lambda,g)\big|^{2}
]+\pp_{\delta_{x}}[
\big|\im W_{t}(\lambda,g)\big|^{2}
]=\mathrm{Var}_{\delta_{x}}[\re W_{t}(\lambda,g)]+\mathrm{Var}_{\delta_{x}}[\im W_{t}(\lambda,g)]+|g(x)|^{2}.$$
Since $\im W_{t}(\lambda,g)=\re W_{t}(\lambda,-\mathrm{i}g)$, it follows by Lemma \ref{lemn2}(iii) that $$\sup_{t>0}\mathrm{Var}_{\delta_{x}}[\re W_{t}(\lambda,g)]<+\infty\quad \mbox{and}\quad  \sup_{t>0}\mathrm{Var}_{\delta_{x}}[\im W_{t}(\lambda,g)]<+\infty.$$
Therefore,
$W_{t}(\lambda,g)$
is an
$L^2(\mathbb{P}_{\delta_x})$-bounded
martingale, which implies the desired result.
\qed

\begin{lemma}\label{lem3.6}
	Suppose $c>\lambda_{1}/2$. Let
	$$E_{+}\slash E_{-}:=\limsup_{t\to+\infty}\slash \liminf_{t\to+\infty}\frac{\e^{(c-\frac{\lambda_{1}}{2})t}\big(\re W_{\infty}(\lambda,g)-\re W_{t}(\lambda,g)\big)}{\sqrt{\log t}}.$$
	Then both $E_{+}$ and $E_{-}$ are finite $\pp_{\delta_{x}}$-a.s.
\end{lemma}
\begin{proof} Fix an arbitrary $\sigma>0$. Let
	$c_{1}(\sigma):=\sqrt{2
		 \Sigma_1(\lambda,g;\sigma)
		}$.
	Lemma \ref{lem3.4} implies that, for any $\epsilon>0$,
	$$\re W_{(n+1)\sigma}(\lambda,g)-\re W_{n\sigma}(\lambda,g)\le
	c_{1}(\sigma)\e^{-(c-\frac{\lambda_{1}}{2})n\sigma}\sqrt{\log n}(\sqrt{W^{\varphi}_{\infty}}+\epsilon)
	\mbox{ eventually }\pp_{\delta_{x}}\mbox{-a.s.}$$
	It follows that, $\pp_{\delta_{x}}$-a.s., for $n$ sufficiently large,
	\begin{align}\label{low2}
		& \frac{\e^{(c-\frac{\lambda_{1}}{2})n\sigma}}{\sqrt{\log n}}\sum_{k=n}^{+\infty}\big(\re W_{(k+1)\sigma}(\lambda,g)-\re W_{k\sigma}(\lambda,g)\big)\nonumber\\
		& \le	c_{1}(\sigma)(\sqrt{W^{\varphi}_{\infty}}+\epsilon)
		\left(\sum_{j=0}^{+\infty}\e^{-(c-\frac{\lambda_{1}}{2})j\sigma}\sqrt{\frac{\log(n+j)}{\log n}}\right).
	\end{align}
	Recall that $W_{t}(\lambda,g)\to W_{\infty}(\lambda,g)$ $\pp_{\delta_{x}}$-a.s. as $t\to+\infty$. The above inequality yields that
	\begin{equation}\label{lem3.4.2}
		\limsup_{\mathbb{N}\ni\to +\infty}\frac{\e^{(c-\frac{\lambda_{1}}{2})n\sigma}\big(\re W_{\infty}(\lambda,g)-\re W_{n\sigma}(\lambda,g)\big)}{\sqrt{\log n}}<+\infty\quad\pp_{\delta_{x}}\mbox{-a.s.}
	\end{equation}
	For $t\ge 0$, let $n(t):=\lfloor t/\sigma \rfloor$ and $\hat{t}:=t-n(t)\sigma$. We have
	\begin{align*}
		\frac{\e^{(c-\frac{\lambda_{1}}{2})t}\big(\re W_{\infty}(\lambda,g)-\re W_{t}(\lambda,g)\big)}{\sqrt{\log t}}=&\e^{(c-\frac{\lambda_{1}}{2})\hat{t}}\sqrt{\frac{\log n(t)}{\log t}}\frac{\e^{(c-\frac{\lambda_{1}}{2})n(t)\sigma}\big(\re W_{\infty}(\lambda,g)-\re W_{n(t)\sigma}(\lambda,g)\big)}{\sqrt{\log n(t)}}\\
		&+\frac{\e^{(c-\frac{\lambda_{1}}{2})t}\big(\re W_{t}(\lambda,-g)-\re W_{n(t)\sigma}(\lambda,-g)\big)}{\sqrt{\log t}}.
	\end{align*}	
	Hence,
	\begin{align}\label{lem3.6.1}
		E_{+}\le&\e^{(c-\frac{\lambda_{1}}{2})\sigma}\limsup_{\mathbb{N}\ni n\to+\infty}\frac{\e^{(c-\frac{\lambda_{1}}{2})n\sigma}\big(\re W_{\infty}(\lambda,g)-\re W_{n\sigma}(\lambda,g)\big)}{\sqrt{\log n}}\nonumber\\
		&+\limsup_{t\to+\infty}\frac{\e^{(c-\frac{\lambda_{1}}{2})t}\big(\re W_{t}(\lambda,-g)-\re W_{n(t)\sigma}(\lambda,-g)\big)}{\sqrt{\log t}}.
	\end{align}
	By \eqref{lem3.5.3}, the second term in the right hand side of \eqref{lem3.6.1} is a.s.
	finite,
	which together with \eqref{lem3.4.2}
	implies that $E_{+}<+\infty$, $\pp_{\delta_{x}}$-a.s.
	Similarly, by investigating $\re W_{\infty}(\lambda,-g)-\re W_{t}(\lambda,-g)$, one can prove that $E_{-}$ is a.s. finite.
\end{proof}

\begin{proof}[Proof of Theorem \ref{them1}(iii)]
	Suppose $(\lambda,g)$ is an eigenpair with $c=\re\lambda>\lambda_{1}/2$.
	We note that, for any $s>0$,
	\begin{align*}
		\frac{\e^{(c-\frac{\lambda_{1}}{2})t}\big(\re W_{t+s}(\lambda,g)-\re W_{t}(\lambda,g)\big)}{\sqrt{\log t}}
		=&\frac{\e^{(c-\frac{\lambda_{1}}{2})t}\big(\re W_{\infty}(\lambda,g)-\re W_{t}(\lambda,g)\big)}{\sqrt{\log t}}\\
		&-\e^{-(c-\frac{\lambda_{1}}{2})s}\sqrt{\frac{\log(t+s)}{\log t}}\frac{\e^{(c-\frac{\lambda_{1}}{2})(t+s)}\big(\re W_{\infty}(\lambda,g)-\re W_{t+s}(\lambda,g)\big)}{\sqrt{\log (t+s)}}.
	\end{align*}
	Combining Proposition \ref{lem3.5} with Lemma \ref{lem3.6} and taking $t\to+\infty$ in the above equation, we get, $\pp_{\delta_{x}}$-a.s.,
	\begin{align*}
		E_{+}-\e^{-(c-\frac{\lambda_{1}}{2})s}E_{+}\le&\sqrt{
		 \big(\Sigma_1(\lambda,g;s) + \big|\Sigma_2(\lambda,g;s) \big)\big|W^{\varphi}_\infty
		}
		\le
		E_{+}-\e^{-(c-\frac{\lambda_{1}}{2})s}E_{-},
	\end{align*}
	and
	\begin{align*}
		E_{-}-\e^{-(c-\frac{\lambda_{1}}{2})s}E_{+}\le&-\sqrt{
			 \big(\Sigma_1(\lambda,g;s) + \big|\Sigma_2(\lambda,g;s) \big)\big|W^{\varphi}_\infty }
			\le
			 E_{-}-\e^{-(c-\frac{\lambda_{1}}{2})s}E_{-}.
	\end{align*}	
	It follows that
	$$\sqrt{W^{\varphi}_{\infty}}a(s)\Big(2-\frac{1}{b(s)}\Big)\le E_{+}\le \sqrt{W^{\varphi}_{\infty}}\frac{a(s)}{b(s)}\quad\pp_{\delta_{x}}\mbox{-a.s.},$$
	where $a(s):=\sqrt{
		 \Sigma_1(\lambda,g;s)+ | \Sigma_2(\lambda,g;s)|
		}$ and
	$b(s):=1-\e^{-(c-\frac{\lambda_{1}}{2})s}$. By letting $s\to +\infty$, we get that $E_{+}
		 = \sqrt{\Sigma(\lambda,g) W^{\varphi}_{\infty}}
		$, $\pp_{\delta_{x}}$-a.s.
	Similarly, by investigating $\re W_{\infty}(\lambda,-g)-\re W_{t}(\lambda,-g)$, we can prove that
	$$E_{-}=-\sqrt{
		 \Sigma(\lambda, g) W^{\varphi}_\infty
		}\quad \pp_{\delta_{x}}\mbox{-a.s.}$$	The second equation of Theorem \ref{them1}(iii) follows from the first equation and Lemma \ref{lem:modula}.
	Hence we complete the proof.	
\end{proof}

\section{Proof of Theorem \ref{them4}}\label{S4}

Throughout this section, we fix an arbitrary $h=\sum_{k=1}^{m}g_{k}\in
\widehat{\mathbb{T}}$.
Recall from the end of Section \ref{sec1.4} that
	$$\re\big(\langle h,X_{t+s}\rangle -\sum_{\re(\gamma_{k})>\lambda_{1}/2}\e^{\gamma_{k}(t+s)}W_{\infty}(\gamma_{k},g_{k})\big)$$
	can be decomposed into the sum of three terms: the non-critical combination $\mathrm{NC}_{h}(s,t)$, the critical combination $\mathrm{CC}_{h}(s+t)$ and the remainder $\mathrm{D}_{h}(s,t)$.
	With this decomposition at hand, the proof of Theorem \ref{them4} proceeds in three main steps.
	
	First, we establish an LIL for $\mathrm{NC}_{h}(s,t)$ as $t\to+\infty$, with normalization factor $\e^{-\lambda_{1}t/2}/\sqrt{2\log t}$. This follows directly from Proposition \ref{lem4.2}, which can be proved similarly to Proposition \ref{lem3.5}, with the main step being to establish the convergence along lattice times (Lemma \ref{lem4.1}).
	
	Second, we obtain an LIL for the critical combination $\mathrm{CC}_{h}(t)$ as $t\to+\infty$, with the normalization factor $\e^{-\lambda_{1}t/2}/\sqrt{2t\log\log t}$. This is accomplished via Lemma \ref{lem5.1}, Proposition \ref{prop4.5} and Proposition \ref{prop:critical case}.
	
	Finally, we show that, $\pp_{\delta_{x}}$-a.s.,
	$\e^{-\lambda_{1}(s+t)/2}\mathrm{D}_{h}(s,t)/\sqrt{\log (t+s)}\to 0$ as first $t\to+\infty$ and then $s\to+\infty$.

We now turn to the detailed proof of these three steps. Let
$$h_{sm}:=\sum_{\re(\gamma_{k})\le \lambda_{1}/2}g_{k}\qquad \mbox{ and }\qquad h_{la}:=\sum_{\re(\gamma_{k})>\lambda_{1}/2}g_{k}.$$
For $0<r<s\le +\infty$, define
\begin{align*}
	\mathcal{Z}_t^{s,r}:= e^{-\lambda_1 t/2} \Big( \langle h, X_{t+r}\rangle  -  \langle T_r h_{sm}, X_t\rangle    - \sum_{\re(\gamma_{k})> \frac{\lambda_1}{2}}   e^{\gamma_k (t+r)} W_{s+t}(\gamma_k, g_k)   \Big)
\end{align*}
and
\begin{align}\label{Def-of-H}
	\mathcal{H}_{s,r}:=    \langle h, X_r\rangle -	\langle T_{r}h_{sm},X_{0}\rangle	-  \sum_{\re(\gamma_{k})> \frac{\lambda_1}{2}} e^{\gamma_k r}W_{s}(\gamma_k, g_k).
\end{align}
\begin{lemma}\label{lem4.1}
	For every $r,\sigma>0$, $\pp_{\delta_{x}}$-a.s.,
	\begin{align*}
		\limsup_{\mathbb{N}\ni n\to +\infty}\frac{\re\big(\mathcal{Z}^{\infty,r}_{n\sigma}\big)}{\sqrt{2\log n}}
		&=  \sqrt{\langle \mathrm{Var}_{\delta_\cdot}[\re(\mathcal{H}_{\infty,r})], \widetilde{\varphi} \rangle W_\infty^\varphi}.
	\end{align*}
\end{lemma}

This lemma will be proved separately for superprocesses and BMPs. We first consider
the superprocess case. Fix $r,\sigma>0$ and an arbitrary $s\in (r,+\infty]$.
	We write $\mathcal{Z}^{s,r}_{t}$ and $\mathcal{H}_{s,r}$, respectively, as $\mathcal{Z}_{t}$ and $\mathcal{H}$ when the values of $r$ and $s$ are unimportant.
	The regular branching property of the superprocess implies that  for every $x\in E$,
	the law of $\re(\mathcal{H})$ under $\mathbb{P}_{\delta_x}$ is infinitely divisible with mean $\pp_{\delta_{x}}[\re(\mathcal{H})]$=0.
	Let $\Phi_x^\mathcal{H}$ be the characteristic exponent of $\re(\mathcal{H})$, i.e.,
	\begin{equation}\label{eq:Phi-H}
		\e^{-\Phi^{\mathcal{H}}_{x}(\theta)}=\pp_{\delta_{x}}\left[\exp\{\mathrm{i} \theta \re(\mathcal{H})\}\right]\quad\forall \theta \in\R.
	\end{equation}
	Then $\Phi_{x}^{\mathcal{H}}$ can be represented as
	$$\Phi^{\mathcal{H}}_{x}(\theta)=\frac{1}{2}\mathrm{Var}_{\delta_{x}}[\re(\mathcal{H})]\theta^{2}+R_{\mathcal{H}}(x,\theta),$$
	where
	$$R_{\mathcal{H}}(x,\theta):=\int_{\R}\big(1-\e^{\mathrm{i}\theta u}+\mathrm{i} \theta  u-\frac{1}{2}\theta^{2} u^{2}\big)\Upsilon^{\mathcal{H}}(x,\mathrm{d}u),$$
	and for every $x\in E$, $\Upsilon^{\mathcal{H}}(x,\mathrm{d}u)$ is a measure on $\R$ satisfying that $\int_{\R}u^{2}\Upsilon^{\mathcal{H}}(x,\mathrm{d}u)<+\infty$.

	\begin{lemma}\label{lem3.1}
		$x\mapsto \int_{\R}u^{4}\Upsilon^{\mathcal{H}}(x,\mathrm{d}u)$ is a bounded function on $E$.
	\end{lemma}
	\begin{proof} The idea of this proof is similar to that of Lemma \ref{lem3.2}. By taking derivatives at $\theta=0$, we get
		$$\int_{\R}u^{4}\Upsilon^{\mathcal{H}}(x,\mathrm{d}u)=-\left.\frac{\mathrm{d}^{4}}{\mathrm{d}\theta^{4}}\Phi^{\mathcal{H}}_{x}(\theta)\right|_{\theta=0}
		=\pp_{\delta_{x}}[
		\re(\mathcal{H})^{4}
		]-
		3
		\mathrm{Var}_{\delta_{x}}[\re(\mathcal{H})]^{2}
.$$
		Thus, it suffices to prove that
		\begin{equation}\label{lem3.1.1}
			\sup_{x\in E}\pp_{\delta_{x}}[		
			\re(\mathcal{H})^{4}
			]<+\infty.
		\end{equation}
		We note that
		$$
		\re(\mathcal{H})^{4}
		\le 2^{3m}\big[\big(\langle \re h,X_{r}\rangle-\langle \re(T_{r}h_{sm}),X_{0}\rangle\big)^{4}+\sum_{\re(\gamma_{k})>\lambda_{1}/2}
		\re\big(\e^{\gamma_{k}r}W_{s}(\gamma_{k},g_{k})\big)^{4}
		\big].$$
		As mentioned
		in the proof of Lemma \ref{lem3.2}, $\sup_{x\in E}\pp_{\delta_{x}}\big[\langle f,X_{t}\rangle^{4}\big]<+\infty$ for any $f\in\mathcal{B}^{+}_{b}(E)$ and $t\in (0,+\infty)$.
		Hence,
		\eqref{lem3.1.1} trivially holds for $s<+\infty$. Now suppose $s=+\infty$. Since $\big|\re\big(\e^{\gamma_{k}r}W_{\infty}(\gamma_{k},g_{k})\big)\big|\le |\e^{\gamma_{k}r}||W_{\infty}(\gamma_{k},g_{k})|$,
		\eqref{lem3.1.1} follows once we prove that,
		for an eigenpair $(\lambda,g)$ with
		$c=\re\lambda>\lambda_{1}/2$,
		\begin{equation}\label{eqnew1}
			\sup_{x\in E}\pp_{\delta_{x}}\left[
			\re\big( W_{\infty}(\lambda,g)\big)^{4}
			\right]<+\infty.
		\end{equation}
		In fact, \eqref{Sup-of-W} yields that both $\sup_{t\ge 0,y\in E}\mathrm{Var}_{\delta_{y}}[\re W_{t}(\lambda,g)]$ and $\sup_{t\ge 0,y\in E}\mathrm{Var}_{\delta_{y}}[\im W_{t}(\lambda,g)]$ are finite. Let
		\[
		M_{1}:=\sup_{t\ge 0,y\in E}\mathrm{Var}_{\delta_{y}}[|W_{t}(\lambda,g)|]<\infty.
		\]
		By \cite[Lemma A.3]{Y25}, there exists
		some $c_{1}>0$ such that for every $x\in E$ and $t>0$,
		\begin{align*}
			&\pp_{\delta_{x}}\left[
			\re\big(W_{t}(\lambda,g)\big)^{4}
			\right]\nonumber\\
			=&\e^{-4c t}\pp_{\delta_{x}}\big[\langle \re\big(\e^{-\mathrm{i}t\im\lambda}g\big),X_{t}\rangle^{4}\big]\nonumber\\
			\le& c_{1}\left(
			\re(g(x))^{4}
			+
		\mathrm{Var}_{\delta_{x}}[\re W_{t}(\lambda,g)]^{2}
		+\int_{0}^{t}\e^{-4 c s}T_{s}\left(\vartheta[\mathrm{Var}_{\delta_{\cdot}}[\re\big(\e^{-\mathrm{i}s\im\lambda}W_{t-s}(\lambda,g)\big)]]\right)(x)\mathrm{d}s\right)\nonumber\\
			\le & c_{1}\left(	
			\re(g(x))^{4}
			+M_1^2+M_1^2\int_{0}^{t}\e^{-4 c s}T_{s}\left(\vartheta[1]\right)(x)\mathrm{d}s\right).
		\end{align*}
		Noticing that both $\vartheta[1](\cdot)$ and $g(\cdot)$ are bounded functions on $E$, we get by \eqref{eq1.9} that for any $t>0$ and $x\in E$,
		\begin{align}\label{eq3.36}
			&\pp_{\delta_{x}}\left[
			\re\big(W_{t}(\lambda,g)\big)^{4}
			\right]\leq
			c_{2}\left( \|g\|_\infty^{4}+M_1^2+\|\vartheta[1]\|_\infty M_1^2
			\int_{0}^{\infty}
			\e^{2(\lambda_1-2c) s}\mathrm{d}s\right).
		\end{align}
		It follows from
		the Fatou's lemma and \eqref{eq3.36} that
		\begin{align*}
			& \pp_{\delta_{x}}\left[
			\re\big(W_{\infty}(\lambda,g)\big)^{4}
			\right]\le \liminf_{t\to+\infty}\pp_{\delta_{x}}\left[
			\re\big(W_{t}(\lambda,g)\big)^{4}
			\right]\nonumber\\
			& \le
			c_{2}\left( \|g\|_\infty^{4}+M_1^2+\|\vartheta[1]\|_\infty M_1^2
			\int_0^\infty
			\e^{2(\lambda_1-2c) s}\mathrm{d}s\right),
		\end{align*}
		which implies \eqref{eqnew1}.
	\end{proof}

Now we are ready to give the proof for superprocesses.

\begin{proof}[Proof of Lemma \ref{lem4.1} for a $(P_{t},\psi)$-superprocess]
Fix $r,\sigma>0$ and an arbitrary $s\in (r,+\infty]$. Recall the definitions of $\mathcal{Z}_{t}$, $\mathcal{H}$, $\Phi^{\mathcal{H}}_{x}(\theta)$ and $R_{\mathcal{H}}(x,\theta)$ given below Lemma \ref{lem4.1}.
	Let $B:=(B_{t})_{t\ge 0}$
	be a standard Brownian motion	 independent of $X$. Set $\widehat{\mathcal{G}}_t:= \mathcal{F}_t \bigvee \sigma(B_s: s\leq t)$.
	For $\varepsilon>0$ and $t\ge 0$,
	define
	\begin{align}\label{variance}
		\mathcal{V}^{2}_{t,\varepsilon}=(\mathcal{V}^{s,r}_{t,\varepsilon})^{2}
		&:=  \mathrm{Var}_{\delta_x} \Big[ \re(\mathcal{Z}_t + \varepsilon (B_{t+r}-B_t) )  \Big| \widehat{\mathcal{G}}_t\Big] = e^{-\lambda_1 t}\langle \mathrm{Var}_{\delta_{\cdot}}[
		\re(\mathcal{H})
		], X_t\rangle + \varepsilon^2 r .
	\end{align}
	By Lemma \ref{lemA.3}, we have
	\begin{align*}
		&\sup_{y\in \R} \left|\mathbb{P}_{\delta_x} \left[ \frac{\re (\mathcal{Z}_t+\varepsilon (B_{t+r}-B_t) )}{\mathcal{V}_{t,\varepsilon} } \leq y\Big|\widehat{\mathcal{G}}_t\right] -\Phi(y)\right|   \nonumber\\
		& \leq \frac{1}{\pi} \int_{-\infty}^{+\infty} \frac{1}{|\theta|} \left| \mathbb{P}_{\delta_x}\Big[ \exp\bigg\{ \mathrm{i}\frac{\theta}{\mathcal{V}_{t,\varepsilon} } \re (\mathcal{Z}_t+\varepsilon (B_{t+r}-B_t) )  \bigg\} \Big| \widehat{\mathcal{G}}_t \Big]-\e^{-\frac{1}{2}\theta^2} \right| \mathrm{d}\theta \nonumber\\
		& =  \frac{1}{\pi} \int_{-\infty}^{+\infty} \frac{1}{|\theta|} \left| e^{-\frac{\theta^2 \varepsilon^2 r}{2\mathcal{V}_{t,\varepsilon}^2}}
		\left.\mathbb{P}_{X_t}\Big[ \exp\bigg\{ \mathrm{i}\hat{\theta}\e^{-\lambda_{1}t/2} \re(\mathcal{H}) \bigg\} \Big]\right|_{\hat{\theta}=\frac{\theta}{\mathcal{V}_{t,\varepsilon}}}
		-\e^{-\frac{1}{2}\theta^2} \right| \mathrm{d}\theta.
	\end{align*}
	According to the definitions of $\Phi_x^{\mathcal{H}}(\theta)$ and $\mathcal{V}_{t,\varepsilon}^2$, we conclude from the above inequality that
	\begin{align*}
		&  \sup_{y\in \R} \left|\mathbb{P}_{\delta_x} \left( \frac{\re (\mathcal{Z}_t + \varepsilon (B_{t+r}-B_t) )}{\mathcal{V}_{t,\varepsilon} } \leq y\Big|\widehat{\mathcal{G}}_t\right) -\Phi(y)\right|  \nonumber\\
		& \leq \frac{1}{\pi} \int_{-\infty}^{+\infty} \frac{1}{|\theta|} \left| \exp\left\{-\frac{\theta^2 \varepsilon^2 r}{2\mathcal{V}_{t,\varepsilon}^2} - \langle \Phi_{\cdot}^{\mathcal{H}}\Big(  \frac{e^{-\lambda_1 t/2}\theta}{\mathcal{V}_{t,\varepsilon}}\Big) , X_t\rangle \right\}-\e^{-\frac{1}{2}\theta^2} \right| \mathrm{d}\theta\nonumber\\
		& = \frac{1}{\pi} \int_{-\infty}^{+\infty} \frac{e^{-\frac{1}{2}\theta^2}}{|\theta|} \left| \exp\left\{-\langle R_{\mathcal{H}}\Big(\cdot,  \frac{e^{-\lambda_1 t/2}\theta}{\mathcal{V}_{t,\varepsilon}}\Big), X_t\rangle \right\}-1 \right| \mathrm{d}\theta.
	\end{align*}
	Based on this, we can apply a similar argument
	as in the proof of Lemma \ref{lem3.4} (with Lemma \ref{lem3.2} replace by Lemma \ref{lem3.1}) to show that,
	$$\sup_{y\in \R} \left|\mathbb{P}_{\delta_x} \left( \frac{\re (\mathcal{Z}_{n\sigma} + \varepsilon (B_{n\sigma+r}-B_{n\sigma}) )}{\mathcal{V}_{n\sigma,\varepsilon} } \leq y\Big|\widehat{\mathcal{G}}_{n\sigma}\right) -\Phi(y)\right|=O(\e^{-\frac{\lambda_{1}n\sigma}{2}})\mbox{ as $n\to+\infty$}\quad\pp_{\delta_{x}}\mbox{-a.s.}$$
	The above equation together with Lemma \ref{lem6.1} implies that,
	for $s\in (r,+\infty)$,
	\begin{align*}
		\limsup_{\mathbb{N}\ni n\to+\infty} \frac{\re (\mathcal{Z}^{s,r}_{n\sigma} + \varepsilon (B_{{n\sigma}+r}-B_{n\sigma}) )}{\sqrt{2\log n}\, \mathcal{V}^{s,r}_{n\sigma,\varepsilon}}=1\quad\pp_{\delta_{x}}\mbox{-a.s.}
	\end{align*}
	By Lemma \ref{lem2.2}, $(\mathcal{V}_{n\sigma,\varepsilon}^{s,r})^2\to \langle \mathrm{Var}_{\delta_\cdot}[\re(\mathcal{H}_{s,r})], \widetilde{\varphi}\rangle W_\infty^\varphi + \varepsilon^2 r$ as $n\to+\infty$, $\pp_{\delta_{x}}$-a.s. Hence, we conclude from the above equality that,
	for $s\in (r,+\infty)$,
	\begin{align*}
		\limsup_{\mathbb{N}\ni n\to+\infty} \frac{\re (\mathcal{Z}_{n\sigma}^{s,r} )+\varepsilon(B_{n\sigma +r}-B_{n\sigma})}{\sqrt{2\log n}} =\sqrt{\langle \mathrm{Var}_{\delta_\cdot}[\re(\mathcal{H}_{s,r})], \widetilde{\varphi}\rangle W_\infty^\varphi + \varepsilon^2 r}\quad\pp_{\delta_{x}}\mbox{-a.s.}
	\end{align*}
	Noticing that
	$\limsup_{\mathbb{N}\ni n\to+\infty} \frac{ B_{n\sigma +r}-B_{n\sigma}}{\sqrt{2 r\log n}}=1$
	almost surely, we obtain that
	\begin{align*}
		\sqrt{\langle \mathrm{Var}_{\delta_\cdot}[\re(\mathcal{H}_{s,r})], \widetilde{\varphi}\rangle W_\infty^\varphi + \varepsilon^2 r} -
		\varepsilon \sqrt{r}
		\leq \limsup_{\mathbb{N}\ni n\to+\infty}
		\frac{\re(\mathcal{Z}_{n\sigma}^{s,r})}{\sqrt{2\log n}}
		\leq \sqrt{\langle \mathrm{Var}_{\delta_\cdot}[\re(\mathcal{H}_{s,r})], \widetilde{\varphi}\rangle W_\infty^\varphi + \varepsilon^2 r}+
		\varepsilon \sqrt{r}.
	\end{align*}
	Taking $\varepsilon\downarrow0$ yields that, for $s\in (r,+\infty)$,
	\begin{align}\label{lim1}
		\limsup_{\mathbb{N}\ni n\to+\infty} \frac{\re (\mathcal{Z}_{n\sigma}^{s,r} )}{\sqrt{2\log n}} =\sqrt{\langle \mathrm{Var}_{\delta_\cdot}[\re(\mathcal{H}_{s,r})], \widetilde{\varphi}\rangle W_\infty^\varphi}\quad\pp_{\delta_{x}}\mbox{-a.s.}
	\end{align}
	Suppose $(\lambda,g)$ is an eigenpair with $c=\re\lambda>\lambda_{1}/2$.
	By \eqref{low2},
	there is some $c_{1}(\sigma)>0$ such that,
	for any $\epsilon>0$ and $l\in\mathbb{N}$,
	$\pp_{\delta_{x}}$-a.s., when $n$ is sufficiently large,
	\begin{equation*}
		\frac{\e^{(c-\frac{\lambda_{1}}{2})n\sigma}}{\sqrt{\log n}}\sum_{k=n+l}^{+\infty}\big(\re W_{(k+1)\sigma}(\lambda,g)-\re W_{k\sigma}(\lambda,g)\big)\le c_{1}(\sigma)(\sqrt{W^{\varphi}_{\infty}}+
		\epsilon)
		\Big(\sum_{j=l}^{+\infty}\e^{-(c-\frac{\lambda_{1}}{2})j\sigma}\sqrt{\frac{\log(n+j)}{\log n}}\Big).
	\end{equation*}
	Recall that $W_{t}(\lambda,g)\to W_{\infty}(\lambda,g)$, $\pp_{\delta_{x}}$-a.s.
	We conclude that
	$$\lim_{\mathbb{N}\ni l\to+\infty}\limsup_{\mathbb{N}\ni n\to+\infty}\frac{\e^{(c-\frac{\lambda_{1}}{2})n\sigma}\big(\re W_{\infty}(\lambda,g)-\re W_{(n+l)\sigma}(\lambda,g)\big)}{\sqrt{\log n}}=0\quad\pp_{\delta_{x}}\mbox{-a.s.}$$
	Moreover,	by taking $g$ to be $-g$ and $-\mathrm{i} g$ in the above equation, it follows that
	\begin{equation}\label{neweq4}
		\lim_{\mathbb{N}\ni l\to+\infty}\limsup_{\mathbb{N}\ni n\to+\infty}\frac{\e^{(c-\frac{\lambda_{1}}{2})n\sigma}\big| W_{\infty}(\lambda,g)-W_{(n+l)\sigma}(\lambda,g)\big|}{\sqrt{\log n}}=0\quad\pp_{\delta_{x}}\mbox{-a.s.}
	\end{equation}
	We note that, for $l\in\mathbb{N}$,
	$$\frac{\re(\mathcal{Z}^{\infty,r}_{n\sigma})}{\sqrt{2\log n}}=\frac{\re(\mathcal{Z}^{l\sigma,r}_{n\sigma})}{\sqrt{2\log n}}+\sum_{\re(\gamma_{k})>\lambda_{1}/2}\frac{\re\big(\e^{\gamma_{k}r}\cdot\e^{(\gamma_{k}-\frac{\lambda_{1}}{2})n\sigma}\big(W_{(n+l)\sigma}(\gamma_{k},g_{k})-W_{\infty}(\gamma_{k},g_{k})\big)
		\big)}{\sqrt{2\log n}}.$$
	We also note that $\langle \mathrm{Var}_{\delta_{\cdot}}[\re(\mathcal{H}_{s,r})],\widetilde{\varphi}\rangle\to \langle\mathrm{Var}_{\delta_{\cdot}}[\re(\mathcal{H}_{\infty,r})],\widetilde{\varphi}\rangle$ as $s\to+\infty$.
	In view of \eqref{lim1} and \eqref{neweq4}, by taking, first $n\to+\infty$ and then $l\to +\infty$ in the above equation, we get the desired result for a $(P_{t},\psi)$-superprocess.
\end{proof}

\begin{proof}[Proof of Lemma \ref{lem4.1} for a $(P_{t},G)$-BMP]
	Fix $r,\sigma>0$.
	Let  $Z:=(Z_{t})_{t\ge 0}$
	be
	an independent continuous-time binary branching process with branching rate $1$. For each $n\in\mathbb{N}$ and $t\ge 0$,
	let $\mathcal{G}_{n}:=
	\mathcal{F}_{n\sigma}
	\vee \sigma(Z_{s}:\ s\le n\sigma)$, $W_t(Z):=e^{-t}Z_t$, and for 	$s\in (r,+\infty)$
	and $\varepsilon>0$, let
	\begin{align*}
		(V^{s,r}_{n\sigma,\varepsilon})^{2}
		&:=  \mathrm{Var}_{\delta_x} \Big[ \re(
		\mathcal{Z}^{s,r}_{n\sigma}
		+ \varepsilon e^{n\sigma/2}(W_{(n+1)\sigma}(Z)-W_{n\sigma}(Z)) )  \Big|
		\mathcal{G}_{n}
		\Big] \nonumber\\
		&=
		\e^{-\lambda_{1}n\sigma}
		\langle \mathrm{Var}_{\delta_{\cdot}}[	
		\re(\mathcal{H}_{s,r})
		], 	X_{n\sigma}
		\rangle + \varepsilon^2 W_{n\sigma}(Z)\mbox{Var}[W_\sigma(Z)].
	\end{align*}
	Applying a similar argument as the one yielding \eqref{neweq1'},
	we can prove that
	\begin{align*}
		\limsup_{\mathbb{N}\ni n\to+\infty} \frac{\re (\mathcal{Z}^{s,r}_{n\sigma} + \varepsilon \e^{n\sigma/2}(W_{(n+1)\sigma}(Z)-W_{n\sigma}(Z)) )}{\sqrt{2\log n} \, V^{s,r}_{n\sigma,\varepsilon}}=1\quad\pp_{\delta_{x}}\mbox{-a.s.}
	\end{align*}
	It follows from
	Lemma \ref{lem2.2}
	that $(	V^{s,r}_{n\sigma,\varepsilon})^{2}	\to \langle\mathrm{Var}_{\delta_\cdot}[\re(\mathcal{H}_{s,r})], \widetilde{\varphi}\rangle W_\infty^\varphi + \varepsilon^2 W_{\infty}(Z)\mbox{Var}[W_\sigma(Z)]$ as $n\to +\infty$, $\mathbb{P}_{\delta_x}$-a.s.  Therefore, we get that
	\begin{align*}
		& \limsup_{\mathbb{N}\ni n\to+\infty} \frac{\re (\mathcal{Z}^{s,r}_{n\sigma} +  \varepsilon \e^{n\sigma/2}(W_{(n+1)\sigma}(Z)-W_{n\sigma}(Z)) )}{\sqrt{2\log n} }\nonumber\\
		&=\sqrt{\langle \mathrm{Var}_{\delta_\cdot}[\re(\mathcal{H}_{s,r})], \widetilde{\varphi}\rangle W_\infty^\varphi + \varepsilon^2 W_{\infty}(Z)\mbox{Var}[W_\sigma(Z)]}\quad\pp_{\delta_{x}}\mbox{-a.s.}
	\end{align*}
	By \eqref{LIL-GW}, it holds that 	
	\begin{align*}
		& \sqrt{\langle \mathrm{Var}_{\delta_\cdot}[\re(\mathcal{H}_{s,r})], \widetilde{\varphi}\rangle W_\infty^\varphi + \varepsilon^2 W_{\infty}(Z)\mbox{Var}[W_\sigma(Z)]} -\sqrt{\varepsilon^2 W_{\infty}(Z)\mbox{Var}[W_\sigma(Z)]}
		\nonumber\\& \leq  \limsup_{\mathbb{N}\ni n\to+\infty} \frac{\re (\mathcal{Z}^{s,r}_{n\sigma}  )}{\sqrt{2\log n} } \nonumber\\
		&\leq \langle \mathrm{Var}_{\delta_\cdot}[\re(\mathcal{H}_{s,r})], \widetilde{\varphi}\rangle W_\infty^\varphi + \varepsilon^2 W_{\infty}(Z)\mbox{Var}[W_\sigma(Z)] +\sqrt{\varepsilon^2 W_{\infty}(Z)\mbox{Var}[W_\sigma(Z)]}\quad\pp_{\delta_{x}}\mbox{-a.s.}
	\end{align*}
	Taking $\varepsilon\downarrow0$, we get \eqref{lim1}.
	The rest of the proof is identical to the superprocess case, thus establishing the desired result.
\end{proof}

\begin{prop}\label{lem4.2}
	For each $h\in \widehat{\mathcal{G}}$ and $r>0$, $\pp_{\delta_{x}}$-a.s.,
	\begin{align*}
		\limsup_{t\to+\infty}\frac{\re\big(\mathcal{Z}^{\infty,r}_{t}\big)}{\sqrt{2\log t}}			
		&=\sqrt{\langle \mathrm{Var}_{\delta_{\cdot}}[\re(\mathcal{H}_{\infty,r})],\widetilde{\varphi}\rangle W^{\varphi}_{\infty}}.
	\end{align*}
	Moreover,
	\begin{align*}
		\langle \mathrm{Var}_{\delta_{\cdot}}[\re(\mathcal{H}_{\infty,r})],\widetilde{\varphi}\rangle=&\e^{\lambda_{1}r}\Big(
		\langle \cc[\re(h_{sm})]-\cc[\e^{-\frac{\lambda_{1}}{2}r}T_{r}(\re(h_{sm}))],\widetilde{\varphi}\rangle+
		\int_0^r e^{-\lambda_1 s}\langle \vartheta[T_{s}(\re(h_{sm}))], \widetilde{\varphi} \rangle \mathrm{d}s\\
		&+ \langle\mathrm{Var}_{\delta_{\cdot}}\Big[\sum_{\re(\gamma_{k})> \lambda_1/2} \re(W_\infty (\gamma_k, g_k))\Big],\widetilde{\varphi}\rangle\Big).
	\end{align*}
\end{prop}

\begin{proof}
	Since $\limsup_{t\to+\infty}\ge
	\limsup_{t=n\sigma\to+\infty}$,
	by Lemma \ref{lem4.1},
	to prove the first conclusion, we only need to show that
	\begin{align}\label{prop4.3.1}
		\limsup_{t\to+\infty}\frac{\re\big(\mathcal{Z}^{\infty,r}_{t}\big)}{\sqrt{2\log t}}			
		&\le  \sqrt{\langle \mathrm{Var}_{\delta_\cdot}[\re(\mathcal{H}_{\infty,r})], \widetilde{\varphi} \rangle W_\infty^\varphi}\quad\pp_{\delta_{x}}\mbox{-a.s.}
	\end{align}
	Taking $(\lambda,g)$ to be $(\gamma_{k},\pm g_k)$ and $(\gamma_{k},\pm \mathrm{i} g_k)$
	in \eqref{lem3.5.3}, we get that, $\pp_{\delta_{x}}$-a.s.,
	\begin{align}\label{lem4.2.1}
		0&= \lim_{\sigma\to 0}\limsup_{n\to+\infty}\sup_{t\in [n\sigma, (n+1)\sigma)} \frac{\e^{(\re(\gamma_{k})-\frac{\lambda_{1}}{2})t}\big| W_{t}(\gamma_k,g_k)- W_{n\sigma}(\gamma_k,g_k)\big|}{\sqrt{\log t}} \nonumber\\
		&=\lim_{\sigma\to 0}\limsup_{n\to+\infty}\sup_{t\in [n\sigma, (n+1)\sigma)} \frac{\e^{-\frac{\lambda_{1} t}{2}}\big| \langle g_k, X_t\rangle - \langle T_{t-n\sigma}g_k, X_{n\sigma}\rangle\big|}{\sqrt{\log t}}.
	\end{align}
	Thus, if $\re(\gamma_{k})>\frac{\lambda_1}{2}$, then
	\begin{align}\label{eq4.1}
		& \limsup_{\sigma\to 0} \limsup_{n\to+\infty}\sup_{t\in [n\sigma, (n+1)\sigma)} \frac{\e^{-\frac{\lambda_{1} t}{2}}\big| \langle g_k, X_t\rangle -e^{\gamma_k t} W_\infty(\gamma_k, g_k) - \Big(
			\langle g_k, X_{n\sigma}\rangle -e^{\gamma_k n\sigma} W_\infty(\gamma_k, g_k)\Big)\big|}{\sqrt{\log t}}\nonumber\\
		&\leq \limsup_{\sigma\to 0} \limsup_{n\to+\infty}\sup_{t\in [n\sigma, (n+1)\sigma)} \frac{\e^{-\frac{\lambda_{1} t}{2}}\big| \langle T_{t-n\sigma}g_k, X_{n\sigma}\rangle -e^{\gamma_k t} W_\infty(\gamma_k, g_k) - \Big(
			\langle g_k, X_{n\sigma}\rangle -e^{\gamma_k n\sigma} W_\infty(\gamma_k, g_k)\Big)\big|}{\sqrt{\log t}}\nonumber\\
		& \leq \limsup_{\sigma\to 0} \limsup_{n\to+\infty}\sup_{t\in [n\sigma, (n+1)\sigma)} \e^{|\gamma_k|\sigma }|\gamma_k| \sigma \frac{\e^{-\frac{\lambda_{1} t}{2}}\big| \langle g_k, X_{n\sigma}\rangle -e^{\gamma_k n\sigma} W_\infty(\gamma_k, g_k) \big|}{\sqrt{\log t}}
		=0\quad\pp_{\delta_{x}}\mbox{-a.s.},
	\end{align}
	where the  final equality is from Theorem \ref{them1}(iii).
	
	On the other hand, if $\re(\gamma_{k}) \leq \frac{\lambda_1}{2}$, then by \eqref{lem4.2.1}, for each $r>0$,
	\begin{align}\label{eq4.2}
		& \limsup_{\sigma\to 0} \limsup_{n\to+\infty}\sup_{t\in [n\sigma, (n+1)\sigma)} \frac{\e^{-\frac{\lambda_{1} t}{2}}\big| \langle g_k, X_t\rangle -e^{-\gamma_k r} \langle g_k, X_{t+r}\rangle - \Big(
			\langle g_k, X_{n\sigma}\rangle -e^{-\gamma_k r}  \langle g_k, X_{n\sigma +r}\rangle\Big)\big|}{\sqrt{\log t}}\nonumber\\
		&\leq \limsup_{\sigma\to 0} \limsup_{n\to+\infty}
		\nonumber\\
		&\qquad  \sup_{t\in [n\sigma, (n+1)\sigma)}
		\frac{\e^{-\frac{\lambda_{1} t}{2}}\big| \langle T_{t-n\sigma}g_k, X_{n\sigma}\rangle -e^{-\gamma_k r} \langle T_{t-n\sigma}g_k, X_{n\sigma+r}\rangle - \Big(
			\langle g_k, X_{n\sigma}\rangle -e^{-\gamma_k r}  \langle g_k, X_{n\sigma +r}\rangle\Big)\big|}{\sqrt{\log t}}\nonumber\\
		& \leq  \limsup_{\sigma\to 0} \limsup_{n\to+\infty}\sup_{t\in [n\sigma, (n+1)\sigma)} \e^{|\gamma_k|\sigma }|\gamma_k| \sigma \frac{\e^{-\frac{\lambda_{1} t}{2}}\big|
			\langle g_k, X_{n\sigma}\rangle -e^{-\gamma_k r}  \langle g_k, X_{n\sigma +r}\rangle\big|}{\sqrt{\log t}}=0\quad\pp_{\delta_{x}}\mbox{-a.s.}
	\end{align}
	Here the final equality is from Proposition \ref{lem3.5}.
	Combining Lemma \ref{lem4.1},  \eqref{eq4.1} and \eqref{eq4.2},  we conclude that, $\pp_{\delta_{x}}$-a.s.,
	\begin{align*}
		\limsup_{t\to+\infty}\frac{\re\big(\mathcal{Z}^{\infty,r}_{t}\big)}{\sqrt{2\log t}}=&\lim_{\sigma\to 0}\limsup_{\mathbb{N}\ni n\to+\infty}\sup_{t\in [n\sigma,(n+1)\sigma)}\frac{\re\big(\mathcal{Z}^{\infty,r}_{t}\big)}{\sqrt{2\log t}}\\
		\le&\lim_{\sigma\to 0}\limsup_{\mathbb{N}\ni n\to+\infty}\sup_{t\in [n\sigma,(n+1)\sigma)}\frac{\e^{-\lambda_{1}(t-n\sigma)/2}\re\big(\mathcal{Z}^{\infty,r}_{n\sigma}\big)}{\sqrt{2\log t}}
		=\sqrt{\langle \mathrm{Var}_{\delta_\cdot}[\re(\mathcal{H}_{\infty,r})], \widetilde{\varphi} \rangle W_\infty^\varphi}.
	\end{align*}
	Hence we prove \eqref{prop4.3.1}.
	
	According to the definition of $\mathcal{H}_{\infty, r}$ in \eqref{Def-of-H}, we have
	\[
	\re({\mathcal{H}_{\infty, r}}) =\re(\langle h_{sm}, X_r\rangle -
	\langle T_{r}h_{sm},X_{0}\rangle
	) - \re\Big( \sum_{\re(\gamma_{k})> \frac{\lambda_1}{2}} e^{\gamma_k r}W_{\infty}(\gamma_k, g_k) - \langle g_k, X_r\rangle  \Big)= : I_r+ II_r.
	\]
	Since $\mathbb P_{\delta_x}\left[ \re({\mathcal{H}_{\infty, r}}) \right]=0$, it holds that
	\begin{align*}
		\mathrm{Var}_{\delta_x}[\re(\mathcal{H}_{\infty,r})] =\mathbb{P}_{\delta_x}\left[ \Big(I_r +II_r\Big)^2 \right] = \mathbb{P}_{\delta_x} [I_r^2] + 2\mathbb{P}_{\delta_x} [I_r II_r]+ \mathbb{P}_{\delta_x} [II_r^2].
	\end{align*}
	Noticing that by the Markov property, $$\mathbb{P}_{\delta_x} [I_r II_r]=
	\pp_{\delta_{x}}\Big[I_{r}\pp_{X_{r}}\Big[\sum_{\re(\gamma_{k})>\lambda_{1}/2}\re\big(\langle g_{k},X_{0}\rangle -W_{\infty}(\gamma_{k},g_{k})\big)\Big]\Big]
	=0,$$
	we conclude that
	\begin{align}\label{step4}
		\mathrm{Var}_{\delta_x}[\re(\mathcal{H}_{\infty,r})] = \mathbb{P}_{\delta_x} [I_r^2] +  \mathbb{P}_{\delta_x} [II_r^2].
	\end{align}
	Combining (H1) and \eqref{eq:variance}, it holds that
	\begin{align}\label{step5}
		\langle \mathbb{P}_{\delta_\cdot} [I_r^2], \widetilde{\varphi} \rangle & = \langle \mathrm{Var}_{\delta_\cdot}[\re(\langle h_{sm}, X_r\rangle)], \widetilde{\varphi} \rangle\nonumber\\
		=&\langle T_{r}(\cc[\re(h_{sm})])-\cc[T_{r}(\re(h_{sm}))],\widetilde{\varphi}\rangle+
		\int_0^r \langle T_{r-s}(\vartheta[T_{s}(\re(h_{sm}))])(\cdot), \widetilde{\varphi} \rangle \mathrm{d}s\nonumber\\
		& = e^{\lambda_1 r}\Big(
		\langle \cc[\re(h_{sm})]-\cc[\e^{-\frac{\lambda_{1}}{2}r}T_{r}(\re(h_{sm}))],\widetilde{\varphi}\rangle+
		\int_0^r e^{-\lambda_1 s}\langle \vartheta[T_{s}(\re(h_{sm}))], \widetilde{\varphi} \rangle \mathrm{d}s\Big).
	\end{align}
	For the second term on the right hand side of \eqref{step4},
	by the Markov property and \eqref{eq:variance},
	\begin{align*}
		\mathbb{P}_{\delta_x} [II_r^2]  & = \mathbb{P}_{\delta_x}\Big[\mathrm{Var}_{X_r}\Big[\sum_{\re(\gamma_{k})> \lambda_1/2} \re(W_\infty (\gamma_k, g_k))\Big]  \Big]  = \mathbb{P}_{\delta_x}\Big[ \langle \mathrm{Var}_{\delta_{\cdot}}\Big[\sum_{\re(\gamma_{k})> \lambda_1/2} \re(W_\infty (\gamma_k, g_k))\Big], X_r\rangle   \Big] \nonumber\\
		& = T_r\left( \mathrm{Var}_{\delta_{\cdot}}\Big[\sum_{\re(\gamma_{k})> \lambda_1/2} \re(W_\infty (\gamma_k, g_k))\Big]\right)(x).
	\end{align*}
	Hence,
	\begin{align}\label{step6}
		\langle   \mathbb{P}_{\delta_\cdot } [II_r^2], \widetilde{\varphi} \rangle = e^{\lambda_1 r} \langle \mathrm{Var}_{\delta_{\cdot}}\Big(\sum_{\re(\gamma_{k})> \lambda_1/2} \re(W_\infty (\gamma_k, g_k))\Big), \widetilde{\varphi} \rangle.
	\end{align}
	Combining \eqref{step4}, \eqref{step5} and \eqref{step6}, we get the second conclusion.
\end{proof}

Next we aim at obtaining an LIL
for combinations of martingales in the critical case.
Let $L\in\mathbb{N}$ and $(\gamma_j, g_j), 1\leq j\leq L$ be eigenpairs such that
$\re(\gamma_{j})=\lambda_{1}/2$
and that $|\im(\gamma_{j})|\neq |\im(\gamma_{k})|$ for $j\neq k$. Define a martingale
\begin{align*}
	\mathcal{W}_t:= \sum_{j=1}^L W_{t}(\gamma_j, g_j).
\end{align*}
A sequence $\{n_{k}:k=0,1,\cdots\}$ of integers is said to be \textit{syndetic}, if $n_{0}<n_{1}<\cdots$ and the gaps $n_{k+1}-n_{k}$ are bounded.
The following lemma is	an extension	of Lemma \ref{lem3.9}.
Though its proof is similar to that of Lemma \ref{lem3.9}, we give the details here for the reader's convenience.

\begin{lemma}\label{lem5.1}
	Suppose $\{n_{k}: k=0,1,\cdots\}$ is a syndetic sequence.
	Then for every $x\in E$ and $\sigma>0$, $\pp_{\delta_{x}}$-a.s.,
	\begin{align*}
		&  \limsup_{\mathbb{N}\ni k\to+\infty}\slash \liminf_{\mathbb{N}\ni k\to+\infty}\frac{\re \mathcal{W}_{n_k\sigma}}{\sqrt{n_k\sigma\log\log n_k}}=(+\slash -)  \sqrt{ \sum_{j=1}^L
		 \Sigma(\gamma_j, g_j)
		 W^{\varphi}_{\infty} }.
	\end{align*}
\end{lemma}
\begin{proof} Fix arbitrary $x\in E$ and $\sigma>0$.
	Let $N:=\sup_{k}(n_{k+1}-n_{k})$ and $\mathcal{G}_{0}:=\mathcal{F}_{n_{0}\sigma}$.
	For $\ell\in\mathbb{N}$,
	define $\mathcal{G}_\ell:= \mathcal{F}_{n_{\ell}\sigma}$ and
	\[
	Y_\ell:=\re \mathcal{W}_{n_\ell \sigma} - \re \mathcal{W}_{n_{\ell-1}\sigma}.
	\]
	Obviously, $\{Y_{\ell}:\ \ell\ge 1\}$ is a martingale difference sequence with respect to the filtration $\{\mathcal{G}_{\ell}:\ \ell\ge 1\}$. Moreover,
	\begin{align*}
		\pp_{\delta_{x}}[Y_{\ell}^{4}|\mathcal{G}_{l-1}]=&
		\pp_{\delta_{x}}\Big[\Big(\sum_{j=1}^{L}\re W_{n_{\ell}\sigma}(\gamma_{j},g_{j})-\re W_{n_{\ell-1}\sigma}(\gamma_{j},g_{j})\Big)^{4}\,|\,\mathcal{G}_{\ell-1}\big]\\ \lesssim&\sum_{j=1}^{L}\sup_{1\leq k\leq N}\pp_{\delta_{x}}\Big[\Big(\re W_{(n_{\ell-1}+k)\sigma}(\gamma_{j},g_{j})-\re W_{n_{\ell-1}\sigma}(\gamma_{j},g_{j})\Big)^{4}\,|\,\mathcal{G}_{\ell-1}\big].
	\end{align*}
	Applying a similar argument (with minor modification) as in the proof of
	Lemma \ref{lem-eq3.20},
	 we can show that
	$$ \sup_{\ell \ge 1, 1\le k\le N} \pp_{\delta_{x}}\Big[\Big(\re W_{(n_{\ell-1}+k)\sigma}(\gamma_{j},g_{j})-\re W_{n_{\ell-1}\sigma} (\gamma_{j},g_{j})\Big)^{4}\,|\,\mathcal{G}_{\ell-1}\big]<+\infty$$
	for each  $1\le j\le L$.
	Consequently, we have $\sup_{\ell\ge 1}\pp_{\delta_{x}}[Y_{\ell}^{4}|\mathcal{G}_{l-1}]<+\infty$, $\pp_{\delta_{x}}$-a.s.
	We note that $n_{k}\asymp k$ as $k\to+\infty$.
	By Lemma \ref{lemA.3.2} in the Appendix,  it suffices to prove that
	\begin{align}\label{Goal}
		\lim_{k\to\infty} \frac{1}{n_k} \sum_{\ell=1}^k \mathbb{P}_{\delta_x}[Y_\ell^2\Big| \mathcal{G}_{\ell-1}] & =  \sigma\sum_{j=1}^L
		 \frac{\Sigma(\gamma_j, g_j) }{2}
		 W^{\varphi}_{\infty} \quad\pp_{\delta_{x}}\mbox{-a.s.}
	\end{align}
	Set $a_\ell:= n_{\ell}-n_{\ell-1}$.  By the Markov property, for each $\ell$,
	\begin{align}\label{step8}
		& \mathbb{P}_{\delta_x}[Y_\ell^2\Big| \mathcal{G}_{\ell-1}]
		\nonumber\\
		& =e^{-\lambda_1n_{\ell-1}\sigma}  \sum_{j=1}^L \mathrm{Var}_{X_{n_{\ell-1}\sigma}} \Big[\re(e^{-\mathrm{i} n_{\ell-1}\sigma\im(\gamma_{j})}W_{a_\ell \sigma}(\gamma_j, g_j))\Big] \nonumber\\
		&\qquad+ e^{-\lambda_1n_{\ell-1}\sigma} \sum_{j\neq k} \mbox{Cov}_{X_{n_{\ell-1}\sigma}} \Big(\re(e^{-\mathrm{i}n_{\ell-1}\sigma\im(\gamma_{j})}W_{a_\ell \sigma}(\gamma_j, g_j) ), \re(e^{-\mathrm{i}n_{\ell-1}\sigma\im(\gamma_{k})}W_{a_\ell \sigma}(\gamma_k, g_k))\Big)\nonumber\\
		&=: I_\ell + II_\ell.
	\end{align}
	By repeating the argument used to establish \eqref{lem3.10.2} and \eqref{lem3.10.3},
	we can show that	
	\begin{align}\label{step8'}
		& \lim_{k\to+\infty} \frac{1}{n_k} \sum_{\ell=1}^k I_\ell = \sigma\sum_{j=1}^L
	 \frac{\Sigma(\gamma_j, g_j)}{2}
	 W^{\varphi}_{\infty}\quad\pp_{\delta_{x}}\mbox{-a.s.}
	\end{align}
	For $II_{\ell}$,
	\eqref{eq:variance}
	and the fact that $\mathrm{Var}_{\mu}[\langle f,X_{t}\rangle]=\langle\mathrm{Var}_{\delta_{\cdot}}[\langle f,X_{t}\rangle],\mu\rangle$
	imply
	that $\mathrm{Cov}_{\mu}(\langle f,X_{t}\rangle,\langle g,X_{t}\rangle)=\langle \mathrm{Cov}_{\delta_{\cdot}}(\langle f,X_{t}\rangle,\langle g,X_{t}\rangle),\mu\rangle$, and
	$$\mathrm{Cov}_{\delta_{x}}(\langle f,X_{t}\rangle,\langle g,X_{t}\rangle)=
	T_{t}(\cc[f,g])(x)-\cc[T_{t}f,T_{t}g](x)+
	\int_{0}^{t}T_{s}\big(\vartheta[T_{t-s}f,T_{t-s}g]\big)(x)\mathrm{d}s.$$
	Hence, we have
	\begin{align*}
		&\mbox{Cov}_{X_{n_{\ell-1}\sigma}} \Big(\re(e^{-\mathrm{i}n_{\ell-1}\sigma\im(\gamma_{j})}W_{a_\ell \sigma}(\gamma_j, g_j) ), \re(e^{-\mathrm{i}n_{\ell-1}\sigma\im(\gamma_{k})}W_{a_\ell \sigma}(\gamma_k, g_k))\Big)\nonumber\\
		&=\langle
		\e^{-\lambda_{1}a_{\ell}\sigma}T_{a_{\ell}\sigma}\big(\cc\big[\re\big(\e^{-\mathrm{i}(n_{\ell-1}+a_{\ell})\sigma\im(\gamma_{j})}g_{j}\big),\re\big(\e^{-\mathrm{i}(n_{\ell-1}+a_{\ell})\sigma\im(\gamma_{k})}g_{k}\big)\big]\big)\nonumber\\
		&\qquad-\cc\big[\re\big(\e^{-\mathrm{i}n_{\ell-1}\sigma\im(\gamma_{j})}g_{j}\big),\re\big(\e^{-\mathrm{i}n_{\ell-1}\sigma\im(\gamma_{k})}g_{k}\big)\big]
		\nonumber\\
		&\qquad+\int_{0}^{a_{\ell}\sigma}\e^{-\lambda_{1}s}T_{s}\big(\vartheta\big[\re\big(\e^{-\mathrm{i}(n_{\ell-1}\sigma+s)\im(\gamma_{j})}g_{j}\big),\re\big(\e^{-\mathrm{i}(n_{\ell-1}\sigma+s)\im(\gamma_{k})}g_{k}\big)\big]\big)\mathrm{d}s,X_{n_{\ell-1}\sigma}\rangle\nonumber\\
		&= \frac{1}{2}\langle
		\re\big(\e^{-\mathrm{i}(n_{\ell-1}+a_{\ell})\sigma(\im(\gamma_{j})-\im(\gamma_{k}))}\cdot\e^{-\lambda_{1}a_{\ell}\sigma}T_{a_{\ell}\sigma}(\cc[g_{j},\bar{g}_{k}])\big)
		\nonumber\\
		&\qquad+\re\big(\e^{-\mathrm{i}(n_{\ell-1}+a_{\ell})\sigma(\im(\gamma_{j})+\im(\gamma_{k}))}\cdot\e^{-\lambda_{1}a_{\ell}\sigma}T_{a_{\ell}\sigma}(\cc[g_{j},g_{k}])\big)
		\nonumber\\
		&\qquad-\re\big(\e^{-\mathrm{i}(n_{\ell-1}+a_{\ell})\sigma(\im(\gamma_{j})-\im(\gamma_{k}))}\cc[g_{j},\bar{g}_{k}])
		-\re\big(\e^{-\mathrm{i}(n_{\ell-1}+a_{\ell})\sigma(\im(\gamma_{j})+\im(\gamma_{k}))}\cc[g_{j},g_{k}]\big)\nonumber\\
		&\qquad+\int_0^{a_\ell\sigma} e^{-\lambda_1 s} T_s \Big(\re\big(e^{-\mathrm{i} (n_{\ell-1}\sigma +s)(\im(\gamma_{j})-\im(\gamma_{k}))} \vartheta[g_j, \bar{g}_k]\big)  \nonumber\\
		&\qquad +   \re\big(e^{-\mathrm{i} (n_{\ell-1}\sigma +s)(\im(\gamma_{j})+\im(\gamma_{k}))} \vartheta[g_j, {g}_k]   \big) \Big)\mathrm{d}s, X_{n_{\ell-1}\sigma}\rangle.
	\end{align*}
	Recall that $|\im(\gamma_{j})|\not= |\im(\gamma_{k})|$ for $j\not=k$.
	By adapting the argument used to establish
	\eqref{neweq12},
	we can prove that
	\begin{align}\label{step9}
		& \lim_{k\to+\infty} \frac{1}{n_k} \sum_{\ell=1}^k I I_\ell = 0\quad\pp_{\delta_{x}}\mbox{-a.s.}
	\end{align}
	Now combining \eqref{step8}, \eqref{step8'} and \eqref{step9}, we get \eqref{Goal},
	which completes the proof of this lemma.
\end{proof}

\begin{prop}\label{prop4.5}
	For every $x\in E$, $\pp_{\delta_{x}}$-a.s.,
	\begin{align*}
		&  \limsup_{t\to+\infty}\slash \liminf_{t\to+\infty}\frac{\re \mathcal{W}_{t}}{\sqrt{t\log\log t}}=(+\slash -)  \sqrt{ \sum_{j=1}^L
		 \Sigma(\gamma_j, g_j)
			W^{\varphi}_{\infty} }
		,
	\end{align*}
	and
	$$\limsup_{t\to+\infty}\frac{ \big|\mathcal{W}_{t}\big|}{\sqrt{t\log\log t}}=(+\slash -)  \sqrt{ \sum_{j=1}^L
	 \Sigma(\gamma_j, g_j)
	 W^{\varphi}_{\infty} }.$$
\end{prop}
\begin{proof} Lemma \ref{lem5.1} yields that, for every $x\in E$ and $\sigma>0$,
	$$\limsup_{\mathbb{N}\ni n\to+\infty}\slash \liminf_{\mathbb{N}\ni n\to+\infty}\frac{\re \mathcal{W}_{n\sigma}}{\sqrt{n\sigma\log\log n}}=(+\slash -)  \sqrt{ \sum_{j=1}^L
		 \Sigma(\gamma_j, g_j)
		 W^{\varphi}_{\infty} }\quad\pp_{\delta_{x}}\mbox{-a.s.}$$
	Then this proposition follows in the same way as Theorem \ref{them1}(ii) (with Lemma \ref{lem3.9} replaced by the above equation).
\end{proof}

We are now ready to prove the LIL
for $\sum_{j=1}^L \langle g_j, X_t \rangle $,
where $\re(\gamma_{j}) =\lambda_1/2$ for all $1\le j\le L$ and $|\im(\gamma_i)|\neq |\im(\gamma_{j})|$ when $i\neq j$.
\begin{prop}\label{prop:critical case}
	For every $x\in E$, $\pp_{\delta_{x}}$-a.s.,
	\begin{align*}
		&  \limsup_{t\to+\infty}\slash \liminf_{t\to+\infty}\frac{ e^{-\lambda_1 t/2}\re \big(\sum_{j=1}^L \langle g_j, X_t \rangle  \big)}{\sqrt{t\log\log t}} =(+\slash -)  \sqrt{ \sum_{j=1}^L
		 \Sigma(\gamma_j, g_j)
		W^{\varphi}_{\infty} }.
	\end{align*}
\end{prop}
\begin{proof}
	We only need to prove the ``$\limsup$" here.
	For each $\theta_1,...,\theta_L\in \mathbb{C}$ with $|\theta_j|=1$, define ${\Theta}:= (\theta_1,...,\theta_L)$ and
	\[
	\Theta\star \mathcal{W}_t:= \sum_{j=1}^L W_t(\gamma_j,\theta_j g_j).
	\]
	An important observation is that,
	\begin{align}\label{step10}
		&  \limsup_{t\to+\infty}
		\frac{\re \big(	\Theta\star \mathcal{W}_t \big)}{\sqrt{t\log\log t}}
		=  \sqrt{ \sum_{j=1}^L
		 \Sigma(\gamma_j, g_j)
		W^{\varphi}_{\infty} }\quad\pp_{\delta_{x}}\mbox{-a.s.},
	\end{align}
	where the limit in the right hand side is independent of $\Theta$.

	\noindent\textit{Upper bound:}
	Let $\mathbb{S}^{1}:=\{z\in\mathbb{C}:\ |z|=1\}$. Fix an arbitrary $\varepsilon>0$. There exists
	a finite set $\{z_{j}:\ 1\le j\le M\}\subset \mathbb{S}^{1}$ satisfying that, for each $z\in \mathbb{S}^{1}$, there is some $1\le j\le M$ such that $|z-z_{j}|\le \varepsilon$. By \eqref{step10}, we have
	\begin{align}\label{step11}
		\limsup_{t\to+\infty} \sup_{\theta_1,...,\theta_L\in \{z_j, 1\leq j\leq M\}}  \frac{ \re \big(\Theta\star \mathcal{W}_t \big) }{\sqrt{t\log\log t}}= \sqrt{ \sum_{j=1}^L
		 \Sigma(\gamma_j, g_j)
			W^{\varphi}_{\infty}}\quad\pp_{\delta_{x}}\mbox{-a.s.}
	\end{align}
	For $1\le k\le L$ and $t>0$, let $z_{k,t}\in \{z_{j}:1\le j\le M\}$ be such that $|\e^{\mathrm{i} t\im(\gamma_{k})}-z_{k,t}|\le \varepsilon$.
	It holds that
	\begin{align*}
		\e^{-\lambda_{1}t/2}\re\big(\sum_{j=1}^{L}\langle g_{j},X_{t}\rangle\big)=&\re\big(\sum_{j=1}^{L}\e^{\mathrm{i} t\im(\gamma_{j})}W_{t}(\gamma_{j},g_{j})\big)\\
		=&\sum_{j=1}^{L}\re W_{t}(\gamma_{j},z_{j,t}g_{j})+\sum_{j=1}^{L}\re\big((\e^{\mathrm{i} t\im(\gamma_{j})}-z_{j,t})W_{t}(\gamma_{j},g_{j})\big).
	\end{align*}
	Then by Theorem \ref{them1}(ii)
	and \eqref{step11}, we have, $\pp_{\delta_{x}}$-a.s.,
	\begin{align*}
		&\limsup_{t\to+\infty}\frac{\e^{-\lambda_{1}t/2}\re\big(\sum_{j=1}^{L}\langle g_{j},X_{t}\rangle\big)}{\sqrt{t\log\log t}}\\
		\le&\limsup_{t\to+\infty} \sup_{\theta_1,...,\theta_L\in \{z_j, 1\leq j\leq M\}}  \frac{ \re \big(\Theta\star \mathcal{W}_t \big) }{\sqrt{t\log\log t}}+\varepsilon\sum_{j=1}^{L}\limsup_{t\to+\infty}\frac{|W_{t}(\gamma_{j},g_{j})|}{\sqrt{t\log\log t}}\\
		=&\sqrt{ \sum_{j=1}^L
			 \Sigma(\gamma_j, g_j)
			W^{\varphi}_{\infty}}+\varepsilon\sum_{j=1}^{L}\sqrt{
			 \Sigma(\gamma_j, g_j)
			W^{\varphi}_{\infty}}.
	\end{align*}
	Letting $\varepsilon\to 0$, we get that
	$$\limsup_{t\to+\infty}\frac{\e^{-\lambda_{1}t/2}\re\big(\sum_{j=1}^{L}\langle g_{j},X_{t}\rangle\big)}{\sqrt{t\log\log t}}\le \sqrt{ \sum_{j=1}^L
	 \Sigma(\gamma_j, g_j)
		W^{\varphi}_{\infty}}\quad\pp_{\delta_{x}}\mbox{-a.s.}$$

	\noindent\textit{Lower bound:}
	According to \cite[Theorem 1.21]{Furstenberg}, for every $\varepsilon>0$, there exists a syndetic sequence $\{n_{k}:k=0,1,\cdots\}$ such that
	\begin{align}\label{claim}
		\sup_{1\leq j\leq L} \Big|e^{\mathrm{i}n_k \im(\gamma_{j})}-1\Big|<\varepsilon.
	\end{align}
	Combining \eqref{claim},
	Theorem \ref{them1}(ii)
	and Lemma \ref{lem5.1} (with $\sigma=1$), it holds that
	\begin{align*}
		&  \limsup_{t\to+\infty}\frac{ e^{-\lambda_1 t/2}\re \big(\sum_{j=1}^L \langle g_j, X_t \rangle \big)}{\sqrt{t\log\log t}}\geq  \limsup_{\mathbb{N}\ni k\to+\infty} \frac{ e^{-\lambda_1 n_k/2}\re \big(\sum_{j=1}^L \langle g_j, X_{n_k}\rangle  \big)}{\sqrt{n_k\log\log n_k}}\nonumber\\
		& = \limsup_{\mathbb{N}\ni k\to+\infty} \frac{ \re \big(\sum_{j=1}^Le^{\mathrm{i} n_k \im(\gamma_{j})} W_{n_k}(\gamma_j, g_j)  \big)}{\sqrt{n_k\log\log n_k}}\nonumber\\
		&\geq  \limsup_{\mathbb{N}\ni k\to+\infty} \frac{ \re \big(\sum_{j=1}^LW_{n_k}(\gamma_j, g_j)  \big)}{\sqrt{n_k\log\log n_k}} - \varepsilon\sum_{j=1}^L \limsup_{t\to+\infty}\frac{\big|  W_{t}(\gamma_j,g_j)\big| }{\sqrt{t\log \log t}} \nonumber\\
		&=\sqrt{ \sum_{j=1}^L
		 \Sigma(\gamma_j, g_j)
			W^{\varphi}_{\infty}}-\varepsilon\sum_{j=1}^{L}\sqrt{
		 \Sigma(\gamma_j, g_j)
			W^{\varphi}_{\infty}}.
	\end{align*}
	Letting $\varepsilon\to 0$, we get that
	$$\limsup_{t\to+\infty}\frac{\e^{-\lambda_{1}t/2}\re\big(\sum_{j=1}^{L}\langle g_{j},X_{t}\rangle\big)}{\sqrt{t\log\log t}}\ge \sqrt{ \sum_{j=1}^L
	 \Sigma(\gamma_j, g_j)
		W^{\varphi}_{\infty}}\quad\pp_{\delta_{x}}\mbox{-a.s.}$$
	Thus we arrive at the desired result.
\end{proof}

\begin{proof}[Proof of Theorem \ref{them4}]
	(i) Suppose $\{\gamma_{k}:\re(\gamma_{k})=\lambda_{1}/2\}=\emptyset$. For any $t,r>0$, we have
	\begin{align}\label{neweq6}
		&\frac{\e^{-\frac{\lambda_{1}}{2}(t+r)}\re\big(\langle h,X_{t+r}\rangle-\sum_{\re(\gamma_{k})>\lambda_{1}/2}\e^{\gamma_{k}(t+r)}W_{\infty}(\gamma_{k},g_{k})\big)}{\sqrt{2\log (t+r)}}\nonumber\\
		=&\frac{\e^{-\frac{\lambda_{1}}{2}(t+r)}\re\big(\langle T_{r}h_{sm},X_{t}\rangle\big)}{\sqrt{2\log(t+r)}}
		+\sqrt{\frac{\log t}{\log(t+r)}}\e^{-\frac{\lambda_{1}r}{2}}\,\frac{\re\big(\mathcal{Z}^{\infty,r}_{t}\big)}{\sqrt{2\log t}}.
	\end{align}
	Note that
	\begin{align*}
		\big|\e^{-\frac{\lambda_{1}}{2}(t+r)}\re\big(\langle T_{r}h_{sm},X_{t}\rangle\big)\big|=&\big|\sum_{\re(\gamma_{k})<\lambda_{1}/2}\re\big(\e^{(\gamma_{k}-\frac{\lambda_{1}}{2})(t+r)}W_{t}(\gamma_{k},g_{k})\big)\big|\\
		\le&\sum_{\re(\gamma_{k})<\lambda_{1}/2}
		\e^{-(\frac{\lambda_{1}}{2}-\re(\gamma_{k}))r}\cdot\e^{(\re(\gamma_{k})-\frac{\lambda_{1}}{2})t}\big|W_{t}(\gamma_{k},g_{k})\big|.
	\end{align*}
	It then follows from Theorem \ref{them1}(i) that
	\begin{equation}\label{neweq5}
		\lim_{r\to+\infty}\limsup_{t\to+\infty}\frac{\big|\e^{-\frac{\lambda_{1}}{2}(t+r)}\re\big(\langle T_{r}h_{sm},X_{t}\rangle\big)\big|}{\sqrt{2\log(t+r)}}=0\quad\pp_{\delta_{x}}\mbox{-a.s.}
	\end{equation}
	We also note that,
	by the second conclusion of Proposition \ref{lem4.2},
	\begin{align*}
		&\e^{-\lambda_{1}r}\langle \mathrm{Var}_{\delta_{\cdot}}[\re(\mathcal{H}_{\infty,r})],\widetilde{\varphi}\rangle\\
		=&\langle \cc[\re(h_{sm})]-\cc[\e^{-\frac{\lambda_{1}}{2}r}T_{r}(\re(h_{sm}))],\widetilde{\varphi}\rangle+
		\int_{0}^{r}\e^{-\lambda_{1}s}\langle\vartheta\big[\sum_{\re(\gamma_{k})<\lambda_{1}/2}\re\big(\e^{\gamma_{k}s}g_{k}\big)\big],\widetilde{\varphi}\rangle \mathrm{d}s\\
		&+\langle \mathrm{Var}_{\delta_{\cdot}}\big[\sum_{\re(\gamma_{k})>\lambda_{1}/2}\re W_{\infty}(\gamma_{k},g_{k})\big],\widetilde{\varphi}\rangle.
	\end{align*}
	It is easy to verify that $\cc[\e^{-\frac{\lambda_{1}}{2}r}T_{r}(\re(h_{sm}))]=\cc\big[\sum_{\re(\gamma_{k})<\lambda_{1}/2}\e^{-(\frac{\lambda_{1}}{2}-\re(\gamma_{k}))r}\re\big(\e^{\mathrm{i}r\im(\gamma_{k})}g_{k}\big)\big]\stackrel{b}{\to}0$
	as $r\to+\infty$. Hence by the dominated theorem, we have
	\begin{align*}
		\e^{-\lambda_{1}r}\langle \mathrm{Var}_{\delta_{\cdot}}[\re(\mathcal{H}_{\infty,r})],\widetilde{\varphi}\rangle\stackrel{r\to+\infty}{\longrightarrow}
		&\langle\cc[\re(h_{sm})],\widetilde{\varphi}\rangle+
		\int_{0}^{+\infty}\e^{-\lambda_{1}s}\langle\vartheta\big[\sum_{\re(\gamma_{k})<\lambda_{1}/2}\re\big(\e^{\gamma_{k}s}g_{k}\big)\big],\widetilde{\varphi}\rangle \mathrm{d}s\\
		&+\langle \mathrm{Var}_{\delta_{\cdot}}\big[\sum_{\re(\gamma_{k})>\lambda_{1}/2}\re W_{\infty}(\gamma_{k},g_{k})\big],\widetilde{\varphi}\rangle\\
		&=H^{sm}_{\infty}(h)+H^{la}_{\infty}(h).
	\end{align*}
	In view of this, \eqref{neweq5} and Proposition \ref{lem4.2}, by letting first $t\to+\infty$ and then $r\to+\infty$ in \eqref{neweq6}, we get that
	\begin{equation}\label{new1}
		\limsup_{t\to+\infty}\frac{\e^{-\frac{\lambda_{1}}{2}t}\re\big(\langle h,X_{t}\rangle-\sum_{\re(\gamma_{k})>\lambda_{1}/2}\e^{\gamma_{k}t}W_{\infty}(\gamma_{k},g_{k})\big)}{\sqrt{2\log t}}
		=\sqrt{\big(H^{sm}_{\infty}(h)+H^{la}_{\infty}(h)\big)W^{\varphi}_{\infty}}
		\quad\pp_{\delta_{x}}\mbox{-a.s.}
	\end{equation}

	\noindent(ii) Suppose $\{\gamma_{k}:\re(\gamma_{k})=\lambda_{1}/2\}\not=\emptyset$. By (i), we have
	$$
	\lim_{t\to+\infty}
	\frac{\e^{-\frac{\lambda_{1}}{2}t}\re\big(\langle \sum_{\re(\gamma_{k})\not=\lambda_{1}/2}g_{k},X_{t}\rangle-\sum_{\re(\gamma_{k})>\lambda_{1}/2}\e^{\gamma_{k}t}W_{\infty}(\gamma_{k},g_{k})\big)}{\sqrt{t\log \log t}}=0\quad\pp_{\delta_{x}}\mbox{-a.s.}$$
	This together with
	Proposition
	\ref{prop:critical case} yields
	 the first conclusion. The second conclusion is a direct result of the first conclusion and Lemma \ref{lem:modula}.
\end{proof}

\appendix
\section{Appendix}
\subsection{Long time behaviour for a class of integrals}
\begin{lemma}\label{lemn1}
	Suppose $\alpha,\theta\in\R$ and $f$ is a $\mathbb{C}$-valued bounded measurable function on $E$.
	\begin{description}
		\item{(i)} If $\alpha<\lambda_{1}/2$, then
		$$\e^{-(\lambda_{1}-2\alpha)t}\int_{0}^{t}\e^{-2\alpha s+\mathrm{i}\theta s}T_{s}f(x)\mathrm{d}s-\varphi(x)\langle f,\widetilde{\varphi}\rangle\frac{\e^{\mathrm{i}\theta t}}{\lambda_{1}-2\alpha+\mathrm{i}\theta}\stackrel{b}{\to}0,\mbox{ as }t\to+\infty.$$
		\item{(ii)} If $\alpha=\lambda_{1}/2$, then
		$$\frac{1}{t}\int_{0}^{t}\e^{-2\alpha s+\mathrm{i}\theta s}T_{s}f(x)\mathrm{d}s\stackrel{b}{\to}
		\begin{cases}
			\varphi(x)\langle f,\widetilde{\varphi}\rangle &\mbox{ if }\theta=0,\\
			0&\mbox{ else}
		\end{cases}\mbox{ as }t\to+\infty.$$
		\item{(iii)} If $\alpha>\lambda_{1}/2$, then
		$$\int_{0}^{t}\e^{-2\alpha s+\mathrm{i}\theta s}T_{s}f(x)\mathrm{d}s\stackrel{b}{\to}\int_{0}^{+\infty}\e^{-2\alpha s+\mathrm{i}\theta s}T_{s}f(x)\mathrm{d}s\mbox{ as }t\to+\infty.$$
	\end{description}
\end{lemma}
\begin{proof} By linearity, we only need to prove these convergence results hold for $f\in\mathcal{B}^{+}_{b}(E)$.
	
	\noindent(i) We have
	$$\e^{-(\lambda_{1}-2\alpha)t}\int_{0}^{t}\e^{-2\alpha s +\mathrm{i}\theta s}T_{s}f(x)\mathrm{d}s=\e^{\mathrm{i}\theta t}\int_{0}^{t}\e^{-(\lambda_{1}-2\alpha+\mathrm{i}\theta)s}\cdot \e^{-\lambda_{1}(t-s)}T_{t-s}f(x)\mathrm{d}s.$$
	Note that $\int_{0}^{t}\e^{-(\lambda_{1}-2\alpha+\mathrm{i}\theta)s}\mathrm{d}s=\big(1-\e^{-(\lambda_{1}-2\alpha+\mathrm{i}\theta)t}\big)/(\lambda_{1}-2\alpha+\mathrm{i}\theta)$. Thus we have
	\begin{align*}
		&\e^{-(\lambda_{1}-2\alpha)t}\int_{0}^{t}\e^{-2\alpha s +\mathrm{i}\theta s}T_{s}f(x)\mathrm{d}s-\varphi(x)\langle f,\widetilde{\varphi}\rangle\frac{\e^{\mathrm{i}\theta t}}{\lambda_{1}-2\alpha+\mathrm{i}\theta}\left(1-\e^{-(\lambda_{1}-2\alpha+\mathrm{i}\theta)t}\right)\\
		=&\e^{\mathrm{i}\theta t}\int_{0}^{t}\e^{-(\lambda_{1}-2\alpha+\mathrm{i}\theta)s}\Big(\e^{-\lambda_{1}(t-s)}T_{t-s}f(x)-\varphi(x)\langle f,\widetilde{\varphi}\rangle\big)\mathrm{d}s.
	\end{align*}
	By \eqref{eq1.9},
	for every $x\in E$ and $s>0$,
	$$\left|\e^{-(\lambda_{1}-2\alpha+\mathrm{i}\theta)s}\Big(\e^{-\lambda_{1}(t-s)}T_{t-s}f(x)-\varphi(x)\langle f,\widetilde{\varphi}\rangle\big)\right|\le c_{1}\e^{-(\lambda_{1}-2\alpha)s}
	\|f\|_{\infty},$$
	and
	$$\e^{-(\lambda_{1}-2\alpha+\mathrm{i}\theta)s}\Big(\e^{-\lambda_{1}(t-s)}T_{t-s}f(x)-\varphi(x)\langle f,\widetilde{\varphi}\rangle\big)\to 0\mbox{ as }t\to+\infty.$$
	Since $\alpha<\lambda_{1}/2$, it follows by the dominated convergence theorem that
	$$\e^{-(\lambda_{1}-2\alpha)t}\int_{0}^{t}\e^{-2\alpha s +\mathrm{i}\theta s}T_{s}f(x)\mathrm{d}s-\varphi(x)\langle f,\widetilde{\varphi}\rangle\frac{\e^{\mathrm{i}\theta t}}{\lambda_{1}-2\alpha+\mathrm{i}\theta}\left(1-\e^{-(\lambda_{1}-2\alpha+\mathrm{i}\theta)t}\right)\stackrel{b}{\to }0\mbox{ as }t\to+\infty.$$
	This together with the fact that $\e^{-(\lambda_{1}-2\alpha)t}/(\lambda_{1}-2\alpha+\mathrm{i}\theta)\stackrel{b}{\to }0$ as $t\to+\infty$ yields (i).

	\noindent(ii) For $\alpha=\lambda_{1}/2$, we have
	\begin{align*}
		\frac{1}{t}\int_{0}^{t}\e^{-2\alpha s+\mathrm{i}\theta s}T_{s}f(x)\mathrm{d}s-\frac{1}{t}\int_{0}^{t}\e^{\mathrm{i}\theta s}\mathrm{d}s \cdot \varphi(x)\langle f,\widetilde{\varphi}\rangle
		=&\frac{1}{t}\int_{0}^{t}\e^{\mathrm{i}\theta s}\Big(\e^{-\lambda_{1}s}T_{s}f(x)-\varphi(x)\langle f,\widetilde{\varphi}\rangle\big)\mathrm{d}s\\
		=&\int_{0}^{1}\e^{\mathrm{i}\theta t r}\Big(\e^{-\lambda_{1}t r}T_{t r}f(x)-\varphi(x)\langle f,\widetilde{\varphi}\rangle\big)\mathrm{d}r.
	\end{align*}
	Again, using \eqref{eq1.9} and the dominated convergence theorem, we can show that
	$$\frac{1}{t}\int_{0}^{t}\e^{-2\alpha s+\mathrm{i}\theta s}T_{s}f(x)\mathrm{d}s-\frac{1}{t}\int_{0}^{t}\e^{\mathrm{i}\theta s}\mathrm{d}s \cdot \varphi(x)\langle f,\widetilde{\varphi}\rangle\stackrel{b}{\to}0\mbox{ as }t\to+\infty.$$
	Note that $\frac{1}{t}\int_{0}^{t}\e^{\mathrm{i}\theta s}\mathrm{d}s=1_{\{\theta=0\}}+\frac{\e^{\mathrm{i}\theta t}-1}{\mathrm{i}\theta t}1_{\{\theta\not=0\}}$. Moreover, for $\theta\not=0$, $\frac{\e^{\mathrm{i}\theta t}-1}{\mathrm{i}\theta t}\stackrel{b}{\to}0$ as $t\to+\infty$. Hence we conclude (ii).

	\noindent (iii) We have
	$$\int_{0}^{t}\e^{-2\alpha s+\mathrm{i}\theta s}T_{s}f(x)\mathrm{d}s=\int_{0}^{+\infty}1_{\{s<t\}}\e^{-(2\alpha-\lambda_{1})s+\mathrm{i}\theta s}\cdot\e^{-\lambda_{1}s}T_{s}f(x)\mathrm{d}s.$$
	By \eqref{eq1.9},
	$|\e^{-\lambda_{1}s}T_{s}f(x)|\le  c_{1}\|f\|_{\infty} $ for all $x\in E$ and $s>0$.
	Since $\alpha>\lambda_{1}/2$, (iii) follows directly from the dominated convergence theorem.
\end{proof}

\subsection{Law of the iterated logarithm for a martingale difference sequence}
Let $\{Y_{n}:\ n\ge 1\}$ be a sequence of random variables defined on a probability space $(\Omega,\mathcal{H},\mathrm{P})$, and $\{\mathcal{H}_{n}:n\ge 0\}$ be an increasing sequence of $\sigma$-fields such that $\mathcal{H}_{0}=\{\emptyset,\Omega\}$, $Y_{n}$ is $\mathcal{H}_{n}$-measurable and $\mathrm{E}[Y_{n}|\mathcal{H}_{n-1}]=0$ for all $n\ge 1$. Let $s^{2}_{n}:=\sum_{k=1}^{n}\mathrm{E}[Y^{2}_{k}|\mathcal{H}_{k-1}]$ and $u_{n}:=\sqrt{2\log\log s^{2}_{n}}$.

\begin{lemma}\label{lemA.3.1}
	If $s^{2}_{n}\to +\infty$ $\mathrm{P}$-a.s. and
	$$\sum_{n=1}^{+\infty}a_{n}\mathrm{E}\left[Y^{2}_{n}1_{\{a_{n}Y^{2}_{n}>1\}}|\mathcal{H}_{n-1}\right]<+\infty\quad\mathrm{P}\mbox{-a.s.,}$$
	where $a_{n}$ are nonnegative $\mathcal{H}_{n-1}$-measurable random variables such that $\frac{u^{2}_{n}}{a_{n}s^{2}_{n}}\to 0$ $\mathrm{P}$-a.s., then
	$$\limsup_{\mathbb{N}\ni n\to+\infty}\frac{\sum_{k=1}^{n}Y_{k}}{s_{n}u_{n}}=1\quad\mathrm{P}\mbox{-a.s.}$$
\end{lemma}

\begin{proof} This lemma is a direct result of \cite[Theorem 3]{Stout} by setting $K^{2}_{n}=\frac{u^{2}_{n}}{a_{n}s^{2}_{n}}$.
\end{proof}

\begin{lemma}\label{lemA.3.2}
	If $\sup_{n}\mathrm{E}\left[Y^{4}_{n}|\mathcal{H}_{n-1}\right]<+\infty$ $\mathrm{P}$-a.s.,
	and $s^{2}_{n}/b(n)\to W$ $\mathrm{P}$-a.s., where $W$ is a nonnegative random variable, and $\{b(n):n\ge 1\}$ is an increasing sequence of numbers such that $b(n)\asymp n$ as $n\to+\infty$, then
	$$\limsup_{\mathbb{N}\ni n\to+\infty}\frac{\sum_{k=1}^{n}Y_{k}}{\sqrt{2b(n)\log\log b(n)}}=\sqrt{W}\quad\mathrm{P}\mbox{-a.s.}$$
\end{lemma}

\begin{proof} Let $\{X_{n}: n\ge 1\}$ be a sequence of i.i.d. random variables with $\mathrm{P}(X=\pm 1)=1/2$ such that $\{X_{n}:n\ge 1\}$ is independent of $\{\mathcal{H}_{n}:n\ge 1\}$. Let $\varepsilon>0$ be an arbitrary constant. Let $\mathcal{G}_{0}:=\mathcal{H}_{0}$, $\mathcal{G}_{n}:=\mathcal{H}_{n}\vee \sigma\{X_{1},\cdots,X_{n}\}$, and $Z_{n}(\varepsilon):=\varepsilon X_{n}+Y_{n}$ for $n\ge 1$. Define $s^{2}_{n}(\varepsilon):=\sum_{k=1}^{n}\mathrm{E}\left[Z^{2}_{k}(\varepsilon)|\mathcal{G}_{k-1}\right]$ and $u_{n}(\varepsilon):=\sqrt{2\log\log s^{2}_{n}(\varepsilon)}$.
	It is easy to verify that $\mathrm{E}\left[Z_{n}(\varepsilon)|\mathcal{H}_{n-1}\right]=0$ for $n\ge 1$, and that $s^{2}_{n}(\varepsilon)=n\epsilon^{2}+s^{2}_{n}$.
	We note that
	$$\mathrm{E}\left[Z^{4}_{n}(\varepsilon)|\mathcal{G}_{n-1}\right]=\mathrm{E}\left[(
	\varepsilon
	X_{n}+Y_{n})^{4}|\mathcal{G}_{n-1}\right]\le 8\left(\varepsilon^{4}\mathrm{E}[X^{4}_{n}|\mathcal{G}_{n-1}]+\mathrm{E}[Y^{4}_{n}|\mathcal{G}_{n-1}]\right)=8\varepsilon^{4}+8\mathrm{E}[Y^{4}_{n}|\mathcal{H}_{n-1}].$$
	If we set $a_{n}=\log n/n$ for $n\ge 1$, then we have
	\begin{eqnarray*}
		\sum_{n=1}^{+\infty}a_{n}\mathrm{E}[Z^{2}_{n}(\varepsilon)1_{\{a_{n}Z^{2}_{n}(\varepsilon)>1\}}|\mathcal{G}_{n-1}]&\le&\sum_{n=1}^{+\infty}a^{2}_{n}\mathrm{E}[Z^{4}_{n}(\varepsilon)|\mathcal{G}_{n-1}]\\
		&\le&8\sum_{n=1}^{+\infty}a^{2}_{n}\big(\varepsilon^{4}+\sup_{k}\mathrm{E}[Y^{4}_{k}\,|\,\mathcal{H}_{k-1}]\big)<+\infty\quad\mathrm{P}\mbox{-a.s.}
	\end{eqnarray*}
	It is easy to see the conditions of Lemma \ref{lemA.3.1} are satisfied by the sequence $\{Z_{n}(\varepsilon):n\ge 1\}$.
	Thus we have
	$$1=\limsup_{\mathbb{N}\ni n\to+\infty}\frac{\sum_{k=1}^{n}Z_{k}(\varepsilon)}{s_{n}(\varepsilon)u_{n}(\varepsilon)}=\limsup_{\mathbb{N}\ni n\to+\infty}\frac{\sum_{k=1}^{n}Y_{k}+\varepsilon\sum_{k=1}^{n}X_{k}}{\sqrt{(n\varepsilon^{2}+s^{2}_{n})\log\log(n\varepsilon^{2}+s^{2}_{n})}}\quad\mathrm{P}\mbox{-a.s.}$$
	We note that
	\begin{align}\label{eqA.3.2}
		\frac{\sum_{k=1}^{n}Y_{k}}{\sqrt{2b(n)\log\log b(n)}}& =\frac{\sum_{k=1}^{n}Z_{k}(\varepsilon)}{s_{n}(\varepsilon)u_{n}(\varepsilon)}\cdot\sqrt{\frac{(n\varepsilon^{2}+s^{2}_{n})\log\log(n\varepsilon^{2}+s^{2}_{n})}{b(n)\log\log b(n)}}
		\nonumber\\
		&\qquad -\frac{\varepsilon\sum_{k=1}^{n}X_{k}}{\sqrt{2n\log \log n}}\cdot\sqrt{\frac{n\log\log n}{b(n)\log\log b(n)}}.
	\end{align}
	Under our assumptions, there are positive constants $c_{1}$ and $c_{2}$ such that $c_{2}b(n)\le n\le c_{1}b(n)$ for $n$ sufficiently large. Using this and the fact that $s^{2}_{n}/b(n)\to W$ a.s., we can easily show that, $\mathrm{P}$-a.s.,
	$$\sqrt{c_{2}\varepsilon^{2}+W}\le \liminf_{n\to+\infty}\sqrt{\frac{(n\varepsilon^{2}+s^{2}_{n})\log\log(n\varepsilon^{2}+s^{2}_{n})}{b(n)\log\log b(n)}}\le\limsup_{n\to+\infty}\sqrt{\frac{(n\varepsilon^{2}+s^{2}_{n})\log\log(n\varepsilon^{2}+s^{2}_{n})}{b(n)\log\log b(n)}}\le \sqrt{c_{1}\varepsilon^{2}+W},$$
	and
	$$\sqrt{c_{2}}\le \liminf_{n\to+\infty}\sqrt{\frac{n\log\log n}{b(n)\log\log b(n)}}\le \limsup_{n\to+\infty}\sqrt{\frac{n\log\log n}{b(n)\log\log b(n)}}\le \sqrt{c_{1}}.$$	
	By the law of iterated logarithm for a simple random walk, we have
	$$\liminf_{\mathbb{N}\ni n\to+\infty}\frac{\sum_{k=1}^{n}X_{k}}{\sqrt{2n\log\log n}}=-1
	\quad\mbox{and}\quad
	\limsup_{\mathbb{N}\ni n\to+\infty}\frac{\sum_{k=1}^{n}X_{k}}{\sqrt{2n\log\log n}}=1\quad\mathrm{P}\mbox{-a.s.}$$
	Hence, by letting $n\to+\infty$ in \eqref{eqA.3.2}, we get that
	$$\sqrt{c_{2}\varepsilon^{2}+W}-\varepsilon\sqrt{c_{1}}\le \limsup_{\mathbb{N}\ni n\to+\infty}
	\frac{\sum_{k=1}^{n}Y_{k}}{\sqrt{2b(n)\log\log b(n)}}\le \sqrt{c_{1}\varepsilon^{2}+W}+\varepsilon\sqrt{c_{1}} \quad \mathrm{P}\mbox{-a.s.}$$
	The conclusion follows by letting $\varepsilon\to 0$.
\end{proof}

	\subsection{Proof of Remark \ref{Rem}(ii)}\label{Appen3}

We first claim that, for any $a,b,c\in\mathbb{R}$,
	\begin{equation}\label{new2}
		\sup_{x,y\in\R,\,x^{2}+y^{2}=1}\big(a x^{2}-2bxy+cy^{2}\big)=\frac{a+c+\sqrt{(a-c)^{2}+4b^{2}}}{2}.
	\end{equation}
	To see this, we let $x=\cos \theta$ and $y=\sin \theta$ for some $\theta\in\R$, then
	\begin{align*}
		ax^{2}-2bxy+cy^{2}=&\frac{a+c}{2}+\frac{a-c}{2}\cos 2\theta-b\sin 2\theta\\
		=&\frac{a+c}{2}+\sqrt{\big(\frac{a-c}{2}\big)^{2}+b^{2}}1_{\{\big(\frac{a-c}{2}\big)^{2}+b^{2}>0\}}\cos(2\theta-\phi),
	\end{align*}
	where $\phi$ satisfies that $\cos \phi=\frac{a-c}{\sqrt{(a-c)^{2}+4b^{2}}}$ and $\sin\phi=\frac{2b}{\sqrt{(a-c)^{2}+4b^{2}}}$. This yields \eqref{new2}.

Now suppose $h=\sum_{k=1}^{m}g_{k}\in\widehat{\mathbb{T}}$ and $\{\gamma_{k}:\ \re(\gamma_{k})=\lambda_{1}/2\}=\emptyset$.
	Noting that $$\im \big(\langle h,X_{t}\rangle-\sum_{\re(\gamma_{k})>\lambda_{1}/2}\e^{\gamma_{k}t}W_{\infty}(\gamma_{k},g_{k})\big)= \re\big(\langle -\mathrm{i}h,X_{t}\rangle-\sum_{\re(\gamma_{k})>\lambda_{1}/2}\e^{\gamma_{k}t}W_{\infty}(\gamma_{k},-\mathrm{i}g_{k})\big),$$
	we have
	\begin{align*}
		&\limsup_{t\to+\infty}\frac{\e^{-\lambda_{1}t}/2}{\sqrt{2\log t}}\Big|\langle h,X_{t}\rangle-\sum_{\re(\gamma_{k})>\lambda_{1}/2}\e^{\gamma_{k}t}W_{\infty}(\gamma_{k},g_{k})\Big|\\
		=&\limsup_{t\to+\infty}\frac{\e^{-\lambda_{1}t}/2}{\sqrt{2\log t}}\nonumber\\
		&\quad \times \sqrt{\big(\re\big(\langle h,X_{t}\rangle-\sum_{\re(\gamma_{k})>\lambda_{1}/2}\e^{\gamma_{k}t}W_{\infty}(\gamma_{k},g_{k})\big)\big)^{2}+\big(\re\big(\langle -\mathrm{i}h,X_{t}\rangle-\sum_{\re(\gamma_{k})>\lambda_{1}/2}\e^{\gamma_{k}t}W_{\infty}(\gamma_{k},-\mathrm{i}g_{k})\big)\big)^{2}}.
	\end{align*}
	Thus, the upper bound in \eqref{new5} follows immediately from the above equation and the first conclusion of Theorem \ref{them4}(i). It remains to prove the lower bound.
	Let $\mathbb{S}^{1}:=\{z\in\mathbb{C}:\ |z|=1\}$. Since $|y|=\sup_{z\in\mathbb{S}^{1}}\re(z\cdot y)$ for any $y\in\mathbb{C}$, we have, for every $z\in\mathbb{S}^{1}$,
	\begin{align}\label{new3}
		&\limsup_{t\to+\infty}\frac{\e^{-\lambda_{1}t}/2}{\sqrt{2\log t}}\Big|\langle h,X_{t}\rangle-\sum_{\re(\gamma_{k})>\lambda_{1}/2}\e^{\gamma_{k}t}W_{\infty}(\gamma_{k},g_{k})\Big|\nonumber \\
		\ge&\limsup_{t\to+\infty}\frac{\e^{-\lambda_{1}t}/2}{\sqrt{2\log t}}\re\big(z\cdot\big(\langle h,X_{t}\rangle-\sum_{\re(\gamma_{k})>\lambda_{1}/2}\e^{\gamma_{k}t}W_{\infty}(\gamma_{k},g_{k})\big)\big)\nonumber\\
		\ge&\sqrt{H^{sm}_{\infty}(z\cdot h)+H^{la}_{\infty}(z\cdot h)}\cdot\sqrt{W^{\varphi}_{\infty}}\quad\pp_{\delta_{x}}\mbox{-a.s.}
	\end{align}
	Simple computation using the fact that $\re(z\cdot y)=\re z\cdot \re y-\im z\cdot \im y$ for $y,z\in\mathbb{C}$ gives that
	$$H^{sm}_{\infty}(z\cdot h)+H^{la}_{\infty}(z\cdot h)=a_{1}(h)(\re z)^{2}-2a_{2}(h)\re z\im z +a_{3}(h)(\im z)^{2},$$
	where the constants $a_{i}(h)$, $i=1,2,3$ are defined above \eqref{new5}.
	Applying \eqref{new2}, we get
	$$\sup_{z\in\mathbb{S}^{1}}\big(H^{sm}_{\infty}(z\cdot h)+H^{la}_{\infty}(z\cdot h)\big)=\frac{a_{1}(h)+a_{3}(h)+\sqrt{(a_{1}(h)-a_{3}(h))^{2}+4a_{2}(h)^{2}}}{2}.$$
	This together with \eqref{new3} yields the lower bound in \eqref{new5}.

\bigskip
\noindent
{\bf Acknowledgements:} The authors would like to thank Jianhua Zhang for helpful discussions.
 We also thank the Associate Editor and two anonymous referees for
their valuable comments and suggestions that have led to the present improved version of the original manuscript.

\end{document}